\documentclass[12pt]{article}
\usepackage{makeidx}
\usepackage{amssymb,amsthm,amsmath,amsfonts,enumerate}
\usepackage{setspace,graphicx,caption,subcaption,float,relsize}
\usepackage{url}
\usepackage{hyperref}
\usepackage{rotating}
\usepackage{natbib}
\usepackage{array}

\hypersetup{
  colorlinks = true, urlcolor = cyan, linkcolor = blue, citecolor = red }
\renewenvironment{abstract}
 {\small
  \begin{center}
  \bfseries \abstractname\vspace{-.0em}\vspace{0pt}
  \end{center}
  \list{}{    \setlength{\leftmargin}{0mm}
    \setlength{\rightmargin}{\leftmargin}  }  \item\relax}
 {\endlist}
\makeatletter
\def\maketag@@@#1{\hbox{\m@th\normalfont\normalsize#1}}
\makeatother
\DeclareMathOperator{\plim}{plim}

\allowdisplaybreaks
\newtheorem {theorem}{Theorem}[section]
\newtheorem {assumption}{Assumption}

\newtheorem{corollary}[theorem]{Corollary}

\newtheorem{lemma}[theorem]{Lemma}

\newtheorem{remark}{Remark}[section]

\allowdisplaybreaks
\input{tcilatex}
\begin{document}

\title{Testing (Infinitely) Many Zero Restrictions}
\author{Jonathan B. Hill\thanks{%
Department of Economics, University of North Carolina, Chapel Hill. E-mail: 
\texttt{jbhill@email.unc.edu}; \texttt{https://jbhill.web.unc.edu}.. We
thank two anonymous referees, an associate editor, and Co-Editor Xiaohong
Chen for many helpful remarks that lead to significant improvements to this
article}\medskip \\
Dept. of Economics, University of North Carolina}
\date{{\large This draft:} \today
}
\maketitle

\begin{abstract}
This paper proposes a max-test for testing (possibly infinitely) many zero
parameter restrictions in an extremum estimation framework. The test
statistic is formed by estimating key parameters one at a time based on many
empirical loss functions that map from a low dimension parameter space, and
choosing the largest in absolute value from these individually estimated
parameters. The parsimoniously parametrized loss identify whether the
original parameter of interest is or is not zero. Estimating fixed low
dimension sub-parameters ensures greater estimator accuracy, does not
require a sparsity assumption, and using only the largest in a sequence of
weighted estimators reduces test statistic complexity and therefore
estimation error, ensuring sharper size and greater power in practice.
Weights allow for standardization in order to control for estimator
dispersion. In a nonlinear parametric regression framework we provide a
parametric wild bootstrap for p-value computation without directly requiring
the max-statistic's limit distribution. A simulation experiment shows the
max-test dominates a conventional bootstrapped test.\bigskip \newline
\textbf{Key words and phrases}: high dimension, max-test, white noise test.
\smallskip \newline
\textbf{AMS classifications} : 62G10, 62M99, 62F35. \smallskip \newline
\textbf{JEL classifications} : C12, C55.
\end{abstract}

\emergencystretch 3em

\setstretch{1.2}

\section{Introduction\label{sec:intro}}

Let $\mathcal{L}$ $:$ $\mathcal{B}$ $\rightarrow $ $[0,\infty )$ be a loss
function with (pseudo) metric parameter space $\mathcal{B}$ $\equiv $ $%
\mathcal{D}$ $\times $ $\Theta $, where $\mathcal{D}$ is a compact subset of 
$\mathbb{R}^{k_{\delta }}$, $k_{\delta }$ $\in $ $\mathbb{N}$, and $\Theta $ 
$\subset $ $\mathbb{R}^{k_{\theta }}$ is a $k_{\theta }$ dimensional linear
space with $k_{\theta }$ $\in $ $\mathbb{N}\cup \infty $. In particular $%
\Theta $ $=$ $\times _{i=1}^{k_{\theta }}\Theta _{i}$ where $\Theta _{i}$ $%
\subset $ $\mathbb{R}$ are compact.\footnote{%
In the finite dimensional case $k_{\theta }$ $<$ $\infty $, $\Theta $ is a
compact subset of the Euclidean space $\mathbb{R}^{k_{\theta }}$. In the
infinite dimensional case $\Theta $ $\subset $ $\mathbb{R}^{\infty }$, $%
(\Theta ,d)$ is a (pseudo) metric space, where $\Theta $ need not be compact
with respect to the (pseudo) metric $d$. In this paper, that is irrelevant
for estimation purposes because we only estimate finite dimensional
parameters on the compact Euclidean subspace $\mathcal{D}$ $\times $ $\Theta
_{i}.$ We assume throughout $\Theta $ $\subset $ $\{\theta $ $\in $ $\mathbb{%
R}^{k_{\theta }}$ $:$ $|\theta |$ $<$ $\infty \}$.}

Assume $\beta _{0}$ $\equiv $ $[\delta _{0}^{\prime },\theta _{0}^{\prime
}]^{\prime }$ is the unique minimizer of $\mathcal{L}(\beta )$ on $\mathcal{B%
}$. This paper presents a test of the sub-parameter $\theta _{0}$ zero
restriction%
\begin{equation*}
H_{0}:\theta _{0}=0\text{ vs. }H_{1}:\theta _{0,i}\neq 0\text{ for at least
one }i.
\end{equation*}%
Although the dimension $k_{\theta }$ of $\theta _{0}$ $=$ $[\theta
_{0,i}]_{i=1}^{k_{\theta }}$ may be infinite, or finite and $k_{\theta }$ $>$
$n$ for a given sample size $n$, we do not require a sparsity condition
under $H_{1}$. It is straightforward to allow for general hypotheses $%
H_{0}:\theta _{0}=\tilde{\theta}$ vs. $H_{1}:\theta _{0}\neq \tilde{\theta}$%
, but testing $H_{0}$ $:$ $\theta _{0}$ $=$ $0$ naturally saves notation.

The main theory in Sections \ref{sec:ident}-\ref{sec:max_test_dist} is
developed at a high level based on a general loss $\mathcal{L}(\beta )$. In
order to facilitate an operational test, a parametric wild bootstrap method
is presented in Section \ref{sec:param_boot} based on the regression model 
\begin{equation}
y_{t}=f(x_{t},\beta _{0})+\epsilon _{t}\text{ with }\beta _{0}=\left[ \delta
_{0}^{\prime },\theta _{0}^{\prime }\right] ^{\prime }\in \mathcal{B},
\label{model}
\end{equation}%
with scalar $y_{t}$, known response $f$ $:$ $\mathbb{R}^{k_{x}}$ $\times $ $%
\mathcal{B}$ $\rightarrow $ $\mathbb{R}$, and $x_{t}$ $\in $ $\mathbb{R}%
^{k_{x}}$ with $k_{x}$ $\in $ $\mathbb{N}\cup \infty $\ are covariates or
design points. We use squared error loss $\mathcal{L}(\beta )$ $\equiv $ $%
.5E[(y_{t}$ $-$ $f(x_{t},\beta ))^{2}]$ in this setting to ease notation. We
assume $(\epsilon _{t},x_{t})$ are independent over $t$ due to technical
challenges with verifying conditions required for a high dimensional
non-Gaussian first order approximation, although $(\epsilon _{t},x_{t})$ may
be heterogeneous, including heteroscedastic. See \cite{Belloni_et_al_2017}
for examples and references concerning high dimensional nonlinear models.

A large (or infinite) dimension for $\beta _{0}$ is a natural possibility in
time series settings due to lags and potentially very high frequency or
tick-by-tick data. Some examples include (infinite order) vector
autoregressions and related multivariate distributed lag models %
\citep[e.g.][]{LutkepohlPoskitt1996,IngWei2003}, forecasting with a large
number of predictors \citep{StockWatson2006,DeMolGiannoneReichlin2008}, and
mixed frequency models \citep{GhyselsHillMotegi2016,GhyselsHillMotegi2017}.
In cross sections, panels and spatial settings an enormous amount of
information is available due to the size of surveys (e.g. U.S. Census,
Current Population Survey, National Longitudinal Survey of Youth),
increasingly sophisticated survey techniques (e.g. text records and word
counts, household transaction data), the unit of observation (e.g. county
wide house prices are dependent on neighbor county prices) and use of dummy
variables for groups (e.g. age, race, region, etc.). See \cite{FanLvQi2011}
and \cite{BelloniChernozhukovHansen2014} for recent surveys.

Inference in this case is typically done post estimation on a regression
model with an imposed sparsity or regularization condition, or related
structural assumption. Sparsity in practice is achieved by a (shrinkage)
penalized estimator, while recent theory allows for an increasing number of
candidate and relevant covariates as the sample size increases %
\citep{FanPeng2004,HuangHorowitzMa2008,MedeirosMendes2015}. Valid inference,
however, typically exists only for the non-zero valued parameter subset,
although there is a nascent literature for gaining inference on the (number
of) zero valued parameters %
\citep[e.g.][]{Carpentier_Verzelen_2019,Carpentier_Verzelen_2021}. Imposing
sparsity at all, however, may be too restrictive, in particular without
justification via a pre-estimation test procedure.

We approach inference on a high dimension parameter \textit{pre}-estimation,
in the sense that the original model, like (\ref{model}), is not estimated.
We exploit Ghysels, Hill and Motegi's (\citeyear{GhyselsHillMotegi2017})
max-test method which operates on many low dimension regression models, each
with one key parameter $\theta _{0,i}$ from the original set of parameters
to be tested, and any other (nuisance) parameters $\delta _{0}$ to be
estimated. Parsimony improves estimation of each $\theta _{0,i}$ when the
dimension $k_{\theta }$ of $\theta _{0}$\ is large, and allows us to
sidestep sparsity considerations. \cite{GhyselsHillMotegi2017} only consider
linear models, a fixed dimension $k_{\theta }$ $\in $ $\mathbb{N}$ with $%
k_{\theta }$ $<$ $n$, and least squares estimation. We generalize their
approach by developing a general theory of hypothesis identification and
estimation based on a general differentiable loss; allowing for an infinite
dimensional $\theta _{0}$ and therefore an increasing number of estimated
parameters; and delivering parametric bootstrap inference for the broad
class of nonlinear regression models (\ref{model}). An increasing number of
estimated parameters leads to both greater generality, and a substantially
different asymptotic theory approach detailed below.

We estimate pseudo true versions $[\delta _{(i)}^{\ast \prime },\theta
_{i}^{\ast }]^{\prime }$ of $[\delta _{0}^{\prime },\theta _{0,i}]^{\prime }$
based on parsimoniously parameterized loss functions $\mathcal{L}%
_{(i)}(\cdot )$ $:$ $\mathcal{D}$ $\times $ $\Theta _{i}$ $\rightarrow $ $%
[0,\infty )$, hence the $i^{th}$ parameterization contains only $\delta $
and $\theta _{i}$. In the case of (\ref{model}), $i$ indexes a parsimonious
version of $f(x_{t},\beta _{0})$. A key result shows that under mild
conditions $[\theta _{i}^{\ast }]_{i=1}^{k_{\theta }}$ $=$ $0$ \textit{if
and only if} $\theta _{0}$ $=$ $0$, and therefore each nuisance parameter $%
\delta _{(i)}^{\ast }$ is identically the (pseudo) true $\delta _{0}$.
Obviously we cannot generally identify $\theta _{0}$ (or $\delta _{0}$)
under the alternative, but we can identify whether $\theta _{0}$ $=$ $0$ is
true or not.

Our chosen parsimonious parameterizations are the \emph{least parameterized}
with nuisance term $\delta _{(i)}^{\ast }$ $\in $ $\mathbb{R}^{k_{\delta }}$
and just one $\theta _{i}^{\ast }$, while in a general setting also allowing
for identification of the null hypothesis: $\theta ^{\ast }$ $=$ $0$ \textit{%
if and only if} $\theta _{0}$ $=$ $0$ (and therefore $\delta _{(i)}^{\ast }$ 
$=$ $\delta _{0}$). Any further parametric reduction, e.g. using $\theta
_{i}^{\ast }$ and a subvector of $\delta _{0}$, cannot generally lead to
identification of the hypotheses. The least viable parameterization is
helpful for estimation purposes (especially under the null): higher degrees
of freedom generally lead to sharper estimates in practice.

Let $\hat{\theta}_{i}$ be an extremum estimator of $\theta _{i}^{\ast }$, $i$
$=$ $1,...,k_{\theta ,n}$, where $k_{\theta ,n}$ $\leq $ $k_{\theta }$ is
the number of parameterizations used in practice. The test statistic is $%
\max_{1\leq i\leq k_{\theta ,n}}|\sqrt{n}\mathcal{W}_{n,i}\hat{\theta}_{i}|$
where $\mathcal{W}_{n,i}$ $>$ $0$ $a.s.$ are possibly stochastic weights.
The weights $\mathcal{W}_{n,i}$\ allow for a variety of test statistics,
including $\max_{1\leq i\leq k_{\theta ,n}}|\sqrt{n}\hat{\theta}_{i}|$, or a
max-t-statistic. Bootstrap inference for the max-statistic does not require
a covariance matrix inversion like Wald and LM tests, it does not require a
specific likelihood function like an LR test, and operates like a shrinkage
estimator by using only the most relevant (weighted) estimator. The
max-tests derived here are asymptotically correctly sized, and consistent
because we let $k_{\theta ,n}$ $\rightarrow $ $k_{\theta }$ $\in $ $\mathbb{N%
}\cup \infty $. Bootstrapping in high dimension has a rich history with
current developments: see, e.g., \cite{BickelFreedman1983}, \cite{Mammen1993}%
, \citet{Chernozhukov_etal2013,Chernozhukov_etal2016}, and \cite%
{JentchPolitis2015}.

The dimension $k_{\theta }$ of $\theta _{0}$ may be larger than the current $%
n$, or infinitely large. The theory developed here applies whether $%
k_{\theta }$ is fixed as we assume (including $k_{\theta }$ $=$ $\infty $),
or increases with the sample size which is common in the shrinkage
literature \citep[e.g][]{vandeGreer2008} and has a deep history %
\citep[e.g.][]{Huber1973}. See \cite{CattaneoJanssonNewey2018} and the
reference therein. In the latter case we would just take $k_{\theta }$ $=$ $%
k_{\theta ,n}$ $\rightarrow $ $\infty $. The dimension $k_{\delta }$ of $%
\delta _{0}$\ is assumed fixed and finite to focus ideas, but can in
principle vary with $n$, or be unbounded (including being a function).%
\footnote{%
In a more general framework, if $\delta _{0}$ is sparse with an increasing
dimension, then a penalized estimator can be used to estimate $[\delta
^{\ast \prime }\theta _{i}^{\ast }]^{\prime }$, like lasso or ridge. If $%
\delta _{0}$ is a (infinite dimensional) function, then $[\delta
_{(i)}^{\ast },\theta _{i}^{\ast }]$ can be estimated by sieves. In this
case inference is on the high dimensional parameter $\theta _{0}$, contrary
to \cite{CattaneoJanssonNewey2018} who propose inference for a fixed low
dimension parameter in a linear setting. The theory presented here does not
cover these cases in order to focus ideas.} Examples broadly speaking
include nonlinear ARX models, random volatility and related processes. In a
linear model $y_{t}$ $=$ $\delta _{0}^{\prime }x_{\delta ,t}$ $+$ $\theta
_{0}^{\prime }x_{\theta ,t}$ $+$ $\epsilon _{t}$, the null aligns with a
test of sub-covariate $x_{\theta ,t}$ inclusion, but in nonlinear models
such a correspondence need not hold. Another application is a test of a
sparsity assumption based on a target set of parameters %
\citep{Tibshirani1996,FanLi2001,CandesTao2007,ZhangHuang2008}. Another is a
(semi)nonparametric regression model based on possibly infinitely many
components from flexible functional forms (via basis expansions), Hermite or
related polynomials, splines, and deep neural networks or machine learning %
\citep[e.g.][]{Andrews1991,Belloni_et_al_2017,CattaneoJanssonNewey2018}.
There is a large literature on high dimensional linear regression with an
increasing number of covariates; consult %
\citet{CattaneoJanssonNewey2018_ET,CattaneoJanssonNewey2018} and \cite%
{LiMuller2021}, amongst others. This literature concerns estimation of high
dimensional models, while we perform inference pre-estimation on fixed low
dimensional models.

Conventional test statistics for testing many zero restrictions typically
exhibit size distortions due both to poor estimator sample properties, and
from inverting a high dimension covariance matrix estimator that may be a
poor approximation of the small or asymptotic variance. The latter is
especially relevant in time series settings when long run covariance matrix
estimation requires a nonparametric approximation. Such estimators are
highly sensitive to tuning parameter choices like bandwidth, and data-driven
bandwidth selection does not alleviate such sensitivity %
\citep[e.g.][]{Shao2010,Politis2011}. The challenge of high dimensional
covariance matrix estimation is well documented %
\citep[consider][]{ChenXuWu2015}. A bootstrap method is therefore typically
applied, but bootstrap tests may only have size corrected power equal to the
original test, which may be low under parameter proliferation %
\citep{DavidsonMacKinnon2006,GhyselsHillMotegi2017}. Indeed, the
bootstrapped Wald test is undersized with low power in linear iid models
when $k_{\delta }$ is even mildly large (e.g. $k_{\delta }$ $=$ $10)$%
\footnote{%
This may carry over to (semi)nonparametric settings with fixed low dimension 
$\delta _{{}}$ and infinite dimensional function $\theta _{0}$. See \cite%
{ChenPouzo2015} for simulation evidence of sieve based bootstrapped t-tests.}%
, and potentially significantly under-sized with low power when $k_{\theta
,n}$ is large, becoming acute when there are nuisance parameters. Max-tests
yield better (typically sharp) size, while power improvements over the Wald
test can be sizable when deviations from the null are small.

Parametric models may be too restrictive in practice, while fully
nonparametric models do not exploit potentially useful parametric
information. See \cite{Chen2007} and \cite{ChenPouzo2015}. Our methods in
principle allow for semi-nonparametric models, like a partially (non)linear
model $y_{t}$ $=$ $f(x_{t},\beta _{0})$ $+$ $h(z_{t})$ $+$ $\epsilon _{t}$
with known response $f$, unknown measurable function $h$ on some infinite
dimensional (pseudo) metric space, and additional covariates $z_{t}$; a
single index model $y_{t}$ $=$ $h(\beta _{0}^{\prime }x_{t})$ $+$ $\epsilon
_{t}$; or a generalized additive partial linear model 
\citep[see][for theory
and references]{CattaneoJanssonNewey2018_ET}. Such models are popularly
estimated by the method of (profile) sieves, which involves estimating
models with increasing dimension/complexity 
\citep[see,
e.g.,][]{Chen2007,ChenLiaoSun2014,ChenPouzo2015}. In order to compute an
approximate p-value for our test, however, we focus on a \textit{parametric}
bootstrap based on reduced dimension versions of (\ref{model}). See \cite%
{HardleHuetMemmen2004} for a wild bootstrap method for semi-nonparametric
models. Thus, we do not treat the infinite dimensional case where $\theta
_{0}$ lies in a general function space, e.g. $y_{t}$ $=$ $\theta _{0}(\delta
_{0}^{\prime }x_{t})$ $+$ $\epsilon _{t}$ for unknown square integrable $%
\theta _{0}$ $:$ $\mathbb{R}$ $\rightarrow $ $\mathbb{R}$.

We discuss below the related literature and highlight our contributions. In
Section \ref{sec:ident} we prove that a sequence of parsimonious loss
functions lead to correct identification of null and alternative hypotheses.
Section \ref{sec:dist} develops the asymptotic theory for a sequence of
extremum estimators based on empirical parsimonious loss. The max-test based
on a bootstrapped p-value approximation is presented in Section \ref%
{sec:pvalue} under the assumption of model (\ref{model}) for independent
data. We verify all assumptions for a linear regression model in Section \ref%
{sec:example}. A Monte Carlo study is presented in Section \ref{sec:sim},
and brief parting comments are left for Section \ref{sec:conclude}. The
supplemental material \cite{supp_mat_testmanyzeros} contains omitted proofs
and an additional logistic regression example (Appendices C and D), and all
simulation results (Appendix E).

\subsection{Related Literature}

The proposed max-test is related to a variety of inference methods, and
exploits key results in statistics. We cannot be exhaustive due to space,
and will instead focus on the most relevant theory and methods: $i$. methods
similar to our reduced dimension approach, including adaptive tests in the
(semi)nonparametrics literature where sieves are used; $ii$. max-tests in
possibly high dimensional settings; and $iii$. high dimensional (Gaussian)
approximations.

\subsubsection{Reduced Dimension Methods}

Reduced dimension regression models in the high dimensional parametric
statistics and machine learning literatures are variously called \textit{%
marginal regression}, \textit{correlation learning}, and \textit{sure
screening} \citep[e.g.][]{Fan_Lv_2008,Genovese_etal_2012}. \cite%
{McKeague_Qian_2015}, for example, regress $y$ against a single covariate $%
x_{i}$ one at a time, $i$ $=$ $1,...,k$ $<$ $\infty $, and then compute the
most relevant index $\hat{k}$ $\equiv $ $\arg \max_{1\leq i\leq k}|\hat{%
\theta}_{i}|$ where $\hat{\theta}_{i}$ $\equiv $ $\widehat{cov}(y,x_{i})/%
\widehat{var}(x_{i})$, and test $\tilde{H}_{0}$ $:$ $\theta _{0}$ $=$ $%
[cov(y,x_{i})/var(x_{i})]_{i=1}^{k}$ $=$ $0$ using $\hat{\theta}_{\hat{k}}$.
We also work pre-estimation, but with key differences. We allow for general
extremum estimators, nonlinear models, and nuisance parameters. Thus we are
not necessarily selecting covariates, using least squares, or testing that 
\textit{no} covariate can be used to predict $y$ linearly. In a linear
setting with squared error loss and no nuisance parameters, however, we
identically test $\tilde{H}_{0}$ with $k$ $=$ $k_{\theta }$, and allow $%
k_{\theta }$ = $\infty $.

Adaptive minimax methods in the (semi)nonparametrics literature often
utilize sequences of sieves for specification tests. In the functional
linear regression $y_{i}$ $=$ $\int_{\mathcal{T}}x_{i}(t)\theta _{0}(t)dt$ $%
+ $ $\epsilon _{i}$, for example, tests of $\theta _{0}$ $=$ $0$ may rely on
projecting $y_{i}$ onto linear subspaces based on $x_{t}$ 
\citep[e.g.][]{HallHorowitz2007,CaiYuan2012,
    HilgertMasVerzelen2013}. Other techniques involve functional principle
components \citep[e.g.,][]{RamsaySilverman1997,HallHorowitz2007}, or
reproducing kernel Hilbert spaces \citep{CaiYuan2012}.

The goal of minimax testing against sparse alternatives is to seek optimal
or sharp convergence rates of (type I and II) error probabilities %
\citep[cf.][]{Ingster1997}. See, e.g., \cite{Baraud_Huet_Laurent_2003} and 
\cite{Baraud_Giraud_Huet_2009} for treatment of mean models without
covariates. \cite{Carpentier_et_al_2019} study minimax rates for tests of $%
\theta _{0}$ $=$ $0$ in the linear model $y_{t}$ $=$ $\theta _{0}^{\prime
}x_{t}$ $+$ $\epsilon _{t}$ against an $l_{2}$ $s$-sparse alternative $0$ $<$
$\tau $ $\leq $ $||\theta _{0}||_{2}$, and $||\theta _{0}||_{0}$ $\leq $ $s$
with $||\theta _{0}||_{0}$ $\equiv $ $\sum_{i=1}^{k_{\theta }}I(\theta
_{0,i} $ $\neq $ $0)$\ for some $s$ $\in $ $\{1,...,k_{\theta }\}$. \cite%
{IngsterSapatinas2009} treat $y_{i}=\theta _{0}(x_{i})+\epsilon _{i}$\ with
square integrable $\theta _{0}$, iid uniform design points $x_{i}$, and
Gaussian error $\epsilon _{i}$. They present an adaptive minimax test of $%
H_{0}$ $:$ $\theta _{0}$ $=$ $\tilde{\theta}$ for given $\tilde{\theta}$,
against a sequence of local alternatives $H_{1}$ $:$ $||\theta _{0}$ $-$ $%
\theta ||$ $\geq $ $r_{n}$ $\rightarrow $ $0$. See, also, \cite%
{Donoho_Jin_2004}, \cite{Arias-Castro_Candes_Plan_2011} and \cite%
{Zhong_Chen_Xu_2013}.

In each case above where covariates appear, minimax rates are sought for
testing whether covariates have predictive ability. Our initial environment
works with a general loss and extremum estimator for testing the parameter $%
\theta _{0}$ $=$ \ $0$, and then a nonlinear parametric model for
bootstrapping a p-value. Thus we allow for greater generality than merely
testing for covariate importance. As noted above, moreover, we do not
require sparsity under $H_{1}$, and we do not estimate the original model.

\cite{ChenPouzo2015} provide a comprehensive theory of inference for
conditional moment (semi)nonparametric models estimated by sieves. They
study Wald and QLR tests of the functional $\phi (\alpha _{0})$ $=$ $\tilde{%
\phi}$ $\in $ $\mathbb{R}^{d}$ for a given mapping $\tilde{\phi}$ of $\alpha
_{0}$ $=$ $[\theta _{0}^{\prime },h_{0}],$ where $\theta _{0}$ is finite
dimensional and $h_{0}$ is an unknown (infinite dimensional) function. They
cover tests of $\theta _{0}$ $=$ $\tilde{\theta}$ for known $\tilde{\theta}$%
, $h_{0}(x)=$ $\tilde{h}(x)$ for a given point $x$ and function $\tilde{h}$,
and much more. See also \cite{ChenLiaoSun2014}.

We focus on testing $\theta _{0}=0$, while $\theta _{0}=\tilde{\theta}$ is
similar. In principle we can allow for the partially nonlinear class $y_{t}$ 
$=$ $f(x_{t},\beta _{0})$ $+$ $h(z_{t})$ $+$ $\epsilon _{t}$, but this
requires a broader bootstrap theory than used here after partially out $%
h(z_{t})$. The major differences, however, are $i$. we allow for infinite
dimensional $\beta _{0}$ $\equiv $ $[\delta _{0}^{\prime },\theta
_{0}^{\prime }]^{\prime }$ where $\theta _{0}$ $\in $ $\times
_{i=1}^{k_{\theta }}\Theta _{i}$ is not an element of a function space; and $%
ii$. we do not use sieves. Indeed, in terms of (\ref{model}) we do not
estimate $\theta _{0}$, and our parsimonious parameterizations all have the
same non-increasing dimension $k_{\delta }$ $+$ $1$, but the number of such
parameterizations increases with $n$. Conversely, sieves spaces increase in
complexity with $n$. In the linear model $y_{t}$ $=$ $\delta _{0}^{\prime
}x_{\delta ,t}$ $+$ $\theta _{0}^{\prime }x_{\theta ,t}$ $+$ $\epsilon _{t}$%
, $\delta $ $\in $ $\mathbb{R}^{k_{\delta }}$, for example, the parsimonious
models are $y_{t}$ $=$ $\delta _{(i)}^{\ast \prime }x_{\delta ,t}$ $+$ $%
\theta _{(i)}^{\ast }x_{\theta ,i,t}$ $+$ $v_{(i),t}$ for $i$ $=$ $%
1,...,k_{\theta ,n}$ each with a $k_{\delta }$ $+$ $1$ dimensional $[\delta
_{(i)}^{\ast \prime },\theta _{(i)}^{\ast }]^{\prime }$, while a\ simple
linear sieve would be $y_{t}$ $=$ $\delta ^{\ast \prime }x_{\delta ,t}$ $+$ $%
\sum_{j=1}^{k_{\theta ,n}}\theta _{n,j}x_{\theta ,t,j}$ $+$ $v_{n,t}$ with
increasing dimension $k_{\delta }$ $+$ $k_{\theta ,n}$. This has been
considered in disparate settings, from AR($\infty $) and VAR($\infty $), to
(semi)nonparametric models %
\citep{LewisReinsel1985,LutkepohlPoskitt1996,CattaneoJanssonNewey2018,CattaneoJanssonNewey2018_ET}%
. Ultimately, our inference approach is not a substitute for sieves, because
in an expanded setting sieves could be used for estimating partially
(non)linear models $y_{t}$ $=$ $f(x_{t},\beta _{0})$ $+$ $h(z_{t})$ $+$ $%
\epsilon _{t}$.

\subsubsection{Max-Statistics}

Interest in extremes for the purpose of model selection arises in diverse
fields , including frontier curve analysis for panel data %
\citep{SchmidtSickles1984}, and the minimum eigenvalue of an income effects
matrix in order to determine the law of demand %
\citep{HardleHildebrandJansen1991}. The use of a max-statistic across
econometric models is used in White's (\citeyear{White2000}) predictive
model selection criterion. Our approach, however, is decidedly different. In
the case of (\ref{model}) we have one model and one hypothesis, from which
many equally low dimensioned models are generated solely for testing the
hypothesis. We also allow the number of low dimensional models to diverge as
the sample size increases.

Max-statistics over subsets of hypotheses appear in different forms: e.g.
over an increasing bandwidth set controlling smoothness of parametric and
nonparametric estimators \citep{Horowitz_Spokoiny_2001}; in mean function
tests \citep{Baraud_Huet_Laurent_2003}; and in (sequential) tests for linear
models against sparse (and faint) alternatives %
\citep{Arias-Castro_Candes_Plan_2011,Zhong_Chen_Xu_2013}. By comparison, we
do not require a linear model, nor approach max-test inference via
sequential testing (requiring size adjustment). Our max-test statistic has a
limit law under the null that is easily bootstrapped, and the test does not
require post sequential testing adjustment.

\subsubsection{High Dimensional Gaussian Approximations}

The limit theory presented here has deep ties to the high dimensional
Gaussian approximation literature for sample means of observed (i.e.
unfiltered) random variables $\{X_{i,t}$ $:$ $i$ $=$ $1,...,k_{n}\}$. \cite%
{Chernozhukov_etal2013}\ profoundly improve on the allowed dimension rate $%
k_{n}$ $\rightarrow $ $\infty $ for independent observed sequences of sample
means. \citet[Appendix
B]{Chernozhukov_etal2019_sm} allow for \textit{almost surely} bounded
stationary $\beta $-mixing processes, and %
\citet{Chernozhukov_etal2019,Chernozhukov_etal2019_sm} consider tests of
(infinitely) many moment inequalities $H_{0}$ $:$ $E[X_{i,t}]$ $\leq $ $0$,
with limited details when a plug-in estimator is involved. Stationary
dependent sequences are covered in \cite{ZhangWu2017} and \cite%
{ZhangCheng2018} among others.

A \textit{Gaussian} approximation theory, however, obviously neglects
approximations for (possibly pre-Gaussian) random functions. The latter is
necessary for handling an expansion of $\hat{\theta}_{i}$ around the pseudo
true $\theta _{i}^{\ast }$. A Gaussian theory\ cannot handle that scenario
because neither object may have a Gaussian distribution.\footnote{%
Consider a linear model without nuisance terms $y_{t}$ $=$ $\theta
_{0}^{\prime }x_{t}$ $+$ $\epsilon _{t}$ with iid square integrable $%
(x_{t},\epsilon _{t})$ and parsimonious models $y_{t}$ $=$ $\theta
_{i}^{\ast }x_{i,t}$ $+$ $v_{i,t}$. The least squares estimator of $\theta
_{i}^{\ast }$ under $H_{0}$ satisfies $\sqrt{n}\hat{\theta}_{i}$ $=$ $%
(E[x_{i,t}^{2}])^{-1}\sum_{t=1}^{n}x_{i,t}\epsilon _{t}/\sqrt{n}$ $+$ $%
\mathcal{\hat{R}}_{i}$, where $\mathcal{\hat{R}}_{i}$ is $o_{p}(1)$ for each 
$i$. Even if $x_{i,t}$ and $\epsilon _{t}$ are mutually independent iid
Gaussian, generally $x_{i,t}\epsilon _{t}$ is not Gaussian. Hence a
different approach is required for proving $\max_{1\leq i\leq k_{,n}}|\sqrt{n%
}\hat{\theta}_{i}$ $-$ $(E[x_{i,t}^{2}])^{-1}\sum_{t=1}^{n}x_{i,t}\epsilon
_{t}/\sqrt{n}|$ $\overset{p}{\rightarrow }$ $0$.} We therefore provide a
general approximation theory for handling the difference between $\hat{\theta%
}_{i}$ and its expansion.\footnote{%
See \citet[Theorem B.2, Comment B.2]{Chernozhukov_etal2019_sm} for the use
of such an expansion, and some discussion on ensuring it holds.} We use a
smooth-max (log-sum-exp) technique exploited in the literature on spin
glasses \citep[cf.][]{Talagrand2011}, machine learning and convex
optimization \citep[see, e.g.,][]{BoydVanderberghe2004}, and high
dimensional Gaussian approximation theory %
\citep[e.g.][]{Chernozhukov_etal2013,Chernozhukov_etal2015,ChangChenWu2021}.
This ultimately forces the existence of a moment generating function for key
first order approximation summands (see Assumption \ref{assum:expand}%
.b(ii),c). Since verification of these properties is greatly expedited under
independence, we work in the setting of \cite{Chernozhukov_etal2013} for
asymptotic distribution and bootstrap limit theory (cf. Section \ref%
{sec:pvalue}).

\cite{HillMotegi2020} propose a max-correlation white noise test using a
fixed dimension parametric filter. The present paper differs significantly.
First, we test zero restrictions in a general extremum estimator framework
that, under high level assumptions, trivially includes a white noise test
couched in an AR framework. Second, they use a general theory to establish
both a high dimensional first order approximation \textit{and} Gaussian
approximation theory for \textit{some} $\{k_{\theta ,n}\}$, $k_{\theta ,n}$ $%
\rightarrow $ $\infty $, but do not establish an upper bound on the rate $%
k_{\theta ,n}$ $\rightarrow $ $\infty $. We exploit a smooth-max
approximation in order to characterize sequences $\{k_{\theta ,n}\}$ that
satisfy a high dimensional first order expansion, and then extant high
dimensional Gaussian approximation theory. Together, we can use any $%
\{k_{\theta ,n}\}$ provided $k_{\theta ,n}$ $=$ $o(\sqrt{n}/(\ln (n))^{4})$.
This slow rate compared to the current Gaussian approximation literature is
due to the (non-Gaussian) first order expansion. The latter generally does
not permit tools like Slepian and Sudakov-Fernique methods, and is not
improved by working in simple settings, like iid linear
regressions.\medskip\ 

We assume all random variables exist on a complete measure space \footnote{%
Completeness of the measure space ensures majorants and integrals over
uncountable families of measurable functions are measurable, and
probabilities where applicable are outer probability measures. Cf. 
\citet[Appendix C:
permissibility criteria]{Pollard1984} and 
\citet[p. 101: admissible
Suslin]{Dudley1984}.}. $|x|$ $=$ $\sum_{i,j}|x_{i,j}|$ is the $l_{1}$-norm, $%
|x|_{2}$ $=$ $(\sum_{i,j}x_{i,j}^{2})^{1/2}$ is the Euclidean or $l_{2}$
norm, and $||A||$ $=$ $\max_{|\lambda |_{2}}\{|A\lambda |_{2}/|\lambda
|_{2}\}$ is the spectral norm for finite dimensional square matrices $A$
(and the Euclidean norm for vectors). $||\cdot ||_{p}$ denotes the $L_{p}$%
-norm. $a.s.$ is \textit{almost surely}. $\boldsymbol{0}_{k}$ denotes a zero
vector with dimension $k$ $\geq $ $1$. Write $r$-vectors as $x$ $\equiv $ $%
[x_{i}]_{i=1}^{r}$. $[\cdot ]$ rounds to the nearest integer. $K$ $>$ $0$ is
non-random and finite, and may take different values in different places. $%
a_{n}$ $\varpropto b_{n}$ implies $a_{n}/b_{n}$ $\rightarrow $ $K$.
Derivatives of functions $f$ $:$ $\mathbb{X}\rightarrow $ $\mathbb{R}$ with
infinite dimensional $\mathbb{X}$ $\subset $ $\mathbb{R}^{\infty }$, denoted 
$(\partial /\partial x)f(x)$, are partial derivatives: $(\partial /\partial
x)f(x)$ $\equiv $ $[(\partial /\partial x_{i})f(x)]_{i=1}^{\infty }$. $awp1$
= asymptotically with probability approaching one.

\section{Identification of Hypotheses\label{sec:ident}}

We first show that low dimension parameterizations leads to correct
identification of the hypotheses in an extremum estimator setting. We do not
require a specific regression model like (\ref{model}), and instead operate
solely on loss $\mathcal{L}(\beta )$ and a parsimonious version $\mathcal{L}%
_{(i)}(\cdot )$ defined below.

\subsection{Parsimonious Loss\label{sec:parsim_model}}

Recall $\mathcal{D}$ and $\Theta _{i}$ are compact subsets of $\mathbb{R}%
^{k_{\delta }}$ and $\mathbb{R}$, $\Theta $ $=$ $\times _{i=1}^{k_{\theta
}}\Theta _{i}$, and denote by $\beta _{0}$ $\equiv $ $[\delta _{0}^{\prime
},\theta _{0}^{\prime }]^{\prime }$ the unique minimizer of loss $\mathcal{L}
$ $:$ $\mathcal{B}$ $\rightarrow $ $[0,\infty )$ on $\mathcal{B}$ $\equiv $ $%
\mathcal{D}$ $\times $ $\Theta $. $\beta _{0}$ may be a true or pseudo-true
value depending on an underlying model specification and choice of loss %
\citep[see][]{Sawa1978,White1982}. In practice $\mathcal{L}(\beta )$ may be
an M-, GMM, or (generalized) empirical likelihood criterion, among others,
where identification under smoothness properties is typically achieved by
first order moment conditions

Now define $k_{\theta ,n}$ $\leq $ $k_{\theta }$ low dimension compact
parameter spaces:%
\begin{equation*}
\mathcal{B}_{(i)}\equiv \mathcal{D}\times \Theta _{i}\text{ hence }\beta
_{(i)}\equiv \left[ \delta ^{\prime },\theta _{i}\right] ^{\prime }\in 
\mathcal{B}_{(i)}\text{, }i=1,...,k_{\theta ,n},
\end{equation*}%
and \textit{parsimonious loss} $\mathcal{L}_{(i)}$ $:$ $\mathcal{B}_{(i)}$ $%
\rightarrow $ $[0,\infty )$. The (pseudo) true value $\beta _{(i)}^{\ast }$
is the unique minimizer of $\mathcal{L}_{(i)}(\beta _{(i)})$ on $\mathcal{B}%
_{(i)}$. Define the set of $\theta _{i}^{\ast \prime }s$: 
\begin{equation*}
\theta ^{\ast }\equiv \left[ \theta _{i}^{\ast }\right] _{i=1}^{k_{\theta
,n}}.
\end{equation*}

If the dimension of $\theta _{0}$ is finite and $k_{\theta }$ $<$ $n$ then
we use $k_{\theta ,n}$ $=$ $k_{\theta }$, and conventional theory applies.
In general we assume%
\begin{equation*}
k_{\theta ,n}\rightarrow k_{\theta }\in \mathbb{N}\cup \infty \text{ \ (}%
k_{\theta ,n}\text{ is monotonically increasing)}.
\end{equation*}

Assume $\mathcal{L}(\cdot )$ and $\mathcal{L}_{(i)}(\cdot )$ are
differentiable, a standard convention in the literature %
\citep[cf.][]{Huber1967,PakesPollard1989}. In order to achieve valid
asymptotics for our max-statistic, however, that \textit{also} yields a
class of characterizable sequences $\{k_{\theta ,n}\}$ with $k_{\theta ,n}$ $%
\rightarrow $ $\infty $ and easy to verify assumptions, we assume with some
sacrifice of generality that the sample version $\mathcal{\hat{L}}%
_{(i)}(\cdot )$ of $\mathcal{L}_{(i)}(\cdot )$\ is twice differentiable with
a sufficiently smooth Hessian. This neglects non-smooth criteria used in
certain penalized estimators (e.g. lasso), quantile, trimmed, and $L_{p}$%
-estimation, non-smooth moment indicators, and non-smooth underlying
regression models (e.g. binary switching models). In regression setting (\ref%
{model}), $\mathcal{\hat{L}}_{(i)}(\cdot )$ is a measurable function of the
sample $\{x_{t},y_{t}\}_{t=1}^{n}$ with twice differentiable response $%
f(\cdot ,\beta )$.

A link between the pseudo-true $\beta _{(i)}^{\ast }$ and $\beta _{0}$ is
required such that $|\theta ^{\ast }$ $-$ $\boldsymbol{0}_{k_{\theta ,n}}|$ $%
\rightarrow $ $0$ as $n$ $\rightarrow $ $\infty $ \textit{if and only if} $%
\theta _{0}$ $=$ $\boldsymbol{0}_{k_{\theta }}$. Hence the nuisance
parameter satisfies $\delta _{(i)}^{\ast }$ $=$ $\delta _{0}$ under the
null, a desirable property. A useful link arises when the parsimonious loss $%
\mathcal{L}_{(i)}(\beta _{(i)})$ operates like a restricted $\mathcal{L}%
(\beta )$, at the null value $\theta $ $=$ $0$:%
\begin{equation}
\frac{\partial }{\partial \beta _{(i)}}\mathcal{L}\left( \delta ,\boldsymbol{%
0}_{k_{\theta }}\right) =\frac{\partial }{\partial \beta _{(i)}}\mathcal{L}%
_{(i)}\left( \delta ,0\right) \text{ for all }i=1,...,k_{\theta }\text{ and
all }\delta \in \mathcal{D},  \label{p_pi_equiv}
\end{equation}%
where $\mathcal{L}(\delta ,\theta )$ $\equiv $ $\mathcal{L}(\beta )$ and $%
\mathcal{L}_{(i)}(\delta ,\theta _{i})$ $\equiv $ $\mathcal{L}_{(i)}(\beta
_{(i)})$. This implies each $\mathcal{L}_{(i)}$ identifies the same $\delta $
under the null. A natural mechanism in practice that promotes (\ref%
{p_pi_equiv}) is setting:%
\begin{equation}
\mathcal{L}_{(i)}(\beta _{(i)})=\mathcal{L}(\delta ,[0,...,\theta
_{i},0,,,]^{\prime })\text{ }\forall \beta _{(i)},  \label{LiL}
\end{equation}%
where $[0,...,\theta _{i},0,...]^{\prime }$ is a zero vector with $i^{th}$
element $\theta _{i}$. Then $(\partial /\partial \beta _{(i)})\mathcal{L}%
_{(i)}(\beta _{(i)})$ $=$ $(\partial /\partial \beta _{(i)})\mathcal{L}%
(\delta ,[0,...,\theta _{i},0,...]^{\prime })$ for each $\beta _{(i)}$,
hence (\ref{p_pi_equiv}) holds.

A simple example is the linear regression: 
\begin{equation}
y_{t}=\delta _{0}^{\prime }x_{\delta ,t}+\theta _{0}^{\prime }x_{\theta
,t}+\epsilon _{t}=\beta _{0}^{\prime }x_{t}+\epsilon _{t}\text{, \ \ }%
E[\epsilon _{t}x_{t}]=0  \label{linear}
\end{equation}%
with parsimonious models \citep[cf.][]{GhyselsHillMotegi2017}:%
\begin{equation}
y_{t}=\delta _{(i)}^{\ast \prime }x_{\delta ,t}+\theta _{i}^{\ast }x_{\theta
,i,t}+v_{(i),t}=\beta _{(i)}^{\ast \prime }x_{(i),t}+v_{(i),t},\text{ }%
i=1,...,k_{\theta ,n}.  \label{linear_parsim}
\end{equation}%
Squared loss and its parsimonious version are $\mathcal{L}(\beta )$ $=$ $%
.5E[(y_{t}$ $-$ $\beta ^{\prime }x_{t})^{2}]$ and $\mathcal{L}_{(i)}(\beta
_{(i)})$ $=$ $.5E[(y_{t}$ $-$ $\beta _{(i)}^{\prime }x_{(i),t})^{2}]$, hence
link (\ref{LiL}) holds. Intuitively when $\theta _{0}$ $=$ $0$ then each $%
\theta _{i}^{\ast }$ $=$ $0$ given orthogonality $E[\epsilon _{t}x_{t}]$ $=$ 
$0$ and squared error loss, hence each $\delta _{(i)}^{\ast }$ $=$ $\delta
_{0}$.

The estimator $\hat{\beta}_{(i)}$ $=$ $[\hat{\delta}_{(i)}^{\prime },\hat{%
\theta}_{i}]^{\prime }$ is the unique minimizer of $\mathcal{\hat{L}}%
_{(i)}(\beta _{(i)})$: $\mathcal{\hat{L}}_{(i)}(\hat{\beta}_{(i)})$ $<$ $%
\mathcal{\hat{L}}_{(i)}(\beta _{(i)})$ $\forall \beta _{(i)}$ $\in $ $\{%
\mathcal{B}_{(i)}$ $:$ $||\beta _{(i)}$ $-$ $\hat{\beta}_{(i)}||$ $>$ $0\}$.
The estimated $\theta ^{\ast }$ $=$ $[\theta _{i}^{\ast }]_{i=1}^{k_{\theta
,n}}$ across parameterizations is:%
\begin{equation*}
\hat{\theta}_{n}\equiv \lbrack \hat{\theta}_{i}]_{i=1}^{k_{\theta ,n}}.
\end{equation*}%
We assume $k_{\theta ,n}$ is small enough for each $n$ to ensure $\hat{\theta%
}_{n}$ exists, and $\hat{\theta}_{i}$ $\overset{p}{\rightarrow }$ $\theta
_{i}^{\ast }$ for each $1$ $\leq $ $i$ $\leq $ $\mathring{k}$ and any fixed $%
\mathring{k}$ $\in $ $\mathbb{N}$, $\mathring{k}$ $\leq $ $%
\lim_{n\rightarrow \infty }k_{\theta ,n}$. In (\ref{linear_parsim}), for
example, if scalar $x_{\delta ,t}$ $=$ $1$ and $x_{\theta ,i,t}$ $=$ $%
y_{t-i} $ then $\hat{\theta}_{i}$ $\overset{p}{\rightarrow }$ $\theta
_{i}^{\ast }$ for each $1$ $\leq $ $i$ $\leq $ $k_{\theta ,n}$ only when $%
k_{\theta ,n}\in $ $\{1,...,n$ $-$ $1\}$ and $k_{\theta ,n}$ $=$ $o(n)$. A
more restrictive bound on $k_{\theta ,n}$ $\rightarrow $ $\infty $, however,
is generally required for a non-Gaussian first order approximation theory
and our bootstrapped p-value (Sections \ref{sec:first_order} and \ref%
{sec:param_boot}).

The max-statistic is

\begin{equation*}
\mathcal{T}_{n}=\max_{1\leq i\leq k_{\theta ,n}}\left\vert \sqrt{n}\mathcal{W%
}_{n,i}\hat{\theta}_{i}\right\vert ,
\end{equation*}%
where $\{\mathcal{W}_{n,i}\}$ are sequences of possibly stochastic weights, $%
\mathcal{W}_{n,i}$ $>$ $0$ $a.s.$ for each $i$, with non-random
(probability) limits $\mathcal{W}_{i}$ $\in $ $(0,\infty )$. The max-test
rejects $H_{0}$ at level $\alpha $ when an appropriate p-value approximation 
$\hat{p}_{n}$ based on $\mathcal{T}_{n}$\ satisfies $\hat{p}_{n}$ $<$ $%
\alpha $ (see Section \ref{sec:param_boot}).

\subsection{Identification: Link between $\protect\theta _{0}$ and $\protect%
\theta ^{\ast }$}

Assumption \ref{assum:ident} ensures $\theta _{0}$ $=$ $\boldsymbol{0}%
_{k_{\theta }}$ \textit{if and only if} $|\theta ^{\ast }$ $-$ $\boldsymbol{0%
}_{k_{\theta ,n}}|$ $\rightarrow $ $0$ as $n$ $\rightarrow $ $\infty $.

\begin{assumption}[Identification]
\label{assum:ident} $\ \ \medskip $\newline
$a.$ The number of parameterizations $k_{\theta ,n}$ $\rightarrow $ $%
k_{\theta }$ $\in $ $\mathbb{N}\cup \infty $.$\medskip $\newline
$b.$ Loss $\mathcal{L}$ $:$ $\mathcal{D}$ $\times $ $\Theta $ $\rightarrow $ 
$[0,\infty )$ is continuous and differentiable on $\mathcal{B}$ $\equiv $ $%
\mathcal{D}\times \Theta $, where $\mathcal{D}$ is a compact subset of $%
\mathbb{R}^{k_{\delta }}$, $k_{\delta }$ $\in $ $\mathbb{N}$, and $\Theta $ $%
=$ $\times _{i=1}^{k_{\theta }}\Theta _{i}$, with non-empty compact $\Theta
_{i}$ $\subset $ $\mathbb{R}$. Moreover:%
\begin{equation}
\frac{\partial }{\partial \beta }\mathcal{L}(\beta )=0\text{ \emph{if and
only if} }\beta =\beta _{0}=\left[ \delta _{0}^{\prime },\theta _{0}^{\prime
}\right] ^{\prime }\in \mathcal{B},  \label{p_orthog}
\end{equation}%
where $\beta _{0}$ is unique, and $\delta _{0}$ and $\theta _{0,i}$ are
interior points of $\mathcal{D}$ and $\Theta _{i}$.$\medskip $\newline
$c.$ Parsimonious loss $\mathcal{L}_{(i)}$ $:$ $\mathcal{B}_{(i)}$ $%
\rightarrow $ $[0,\infty )$ are continuous and differentiable on compact $%
\mathcal{B}_{(i)}$ $\equiv $ $\mathcal{D}\times \Theta _{i}$. $\beta
_{(i)}^{\ast }$ are the unique interior points of $\mathcal{D}$ $\times $ $%
\Theta _{i}$ that satisfy 
\begin{equation}
\frac{\partial }{\partial \beta _{(i)}}\mathcal{L}_{(i)}(\beta _{(i)})=0%
\text{ \emph{if and only if} }\beta _{(i)}=\beta _{(i)}^{\ast }=[\delta
_{(i)}^{\ast \prime },\theta _{i}]^{\prime }\in \mathcal{D}\times \Theta
_{i}.  \label{pi_orthog}
\end{equation}%
$d$. The loss functions $\mathcal{L}(\beta )$ and $\mathcal{L}_{(i)}(\beta
_{(i)})$\ are linked by (\ref{LiL}): $\mathcal{L}_{(i)}(\beta _{(i)})=%
\mathcal{L}(\delta ,[0,...,\theta _{i},0,...]^{\prime })$ $\forall \beta
_{(i)}.$
\end{assumption}

\begin{remark}
\normalfont($b$) and ($c$) are standard assumptions (for smooth criteria),
and neither requires $\beta _{0}$ or $\beta _{(i)}^{\ast }$ to be part of
the true data generating process.
\end{remark}

In our first main result of the paper we prove our class of low dimension
parameterizations can be used to identify whether $H_{0}$ is true or not,
when exactly $k_{\theta }$ parameterizations are used. The result is
important compared to the extant statistics and econometrics literatures
because in our framework we do not require a (linear) regression model, nor
specified (quadratic) loss. The workhorse condition is gradient linkage (\ref%
{p_pi_equiv}) which holds under loss equivalence (\ref{LiL}).

The proof itself is instructive, and differs significantly from arguments in %
\citet[proof of Theorem 2.1]{GhyselsHillMotegi2017}, hence we present it
below.

\begin{theorem}
\label{th:ident}Let $(\theta _{0},\theta ^{\ast })$ $\in $ $\mathbb{R}%
^{k_{\theta }}$. Under Assumption \ref{assum:ident}.b-d, $\theta _{0}$ $=$ $%
\boldsymbol{0}_{k_{\theta }}$ \emph{if and only if} $\theta ^{\ast }$ $=$ $%
\boldsymbol{0}_{k_{\theta }}$. Under $H_{0}$ $:$ $\theta _{0}$ $=$ $%
\boldsymbol{0}_{k_{\theta }}$ it therefore follows that $\delta ^{\ast }$ $=$
$\delta _{0}$, and under $H_{1}$ $:$ $\theta _{0}$ $\neq $ $\boldsymbol{0}%
_{k_{\theta }}$ there exists $i$, $1$ $\leq $ $i$ $\leq $ $k_{\theta }$,
such that $\theta _{i}^{\ast }$ $\neq $ $0$.
\end{theorem}

\textbf{Proof.}\qquad If $\theta _{0}$ $=$ $\boldsymbol{0}_{k_{\theta }}$
then by orthogonality (\ref{p_orthog}) it follows $(\partial /\partial \beta
)\mathcal{L}([\delta _{0},\boldsymbol{0}_{k_{\theta }}])$ $=$ $\boldsymbol{0}%
_{k_{\delta }+k_{\theta }}$, hence $(\partial /\partial \beta _{(i)})%
\mathcal{L}([\delta _{0},\boldsymbol{0}_{k_{\theta }}])$ $=$ $\boldsymbol{0}%
_{k_{\delta }+1}$. Loss gradient equivalence (\ref{p_pi_equiv}) follows from
loss equivalence (\ref{LiL}) under Assumption \ref{assum:ident}.d, and
therefore yields $(\partial /\partial \beta _{(i)})\mathcal{L}_{(i)}([\delta
_{0},0])$ $=$ $\boldsymbol{0}_{k_{\delta }+1}$ for each $i$. Hence, $\theta
_{i}^{\ast }$ $=$ $0$ for each $i$ and $\delta ^{\ast }$ $=$ $\delta _{0}$\
by the construction and uniqueness of $\beta _{(i)}^{\ast }$ $=$ $[\delta
^{\ast \prime },\theta _{i}^{\ast }]^{\prime }$.

Conversely, if $\theta ^{\ast }$ $=$ $\boldsymbol{0}_{k_{\theta }}$ then $%
(\partial /\partial \beta _{(i)})\mathcal{L}_{(i)}([\delta _{(i)}^{\ast
},0]) $ $=$ $0$ for each $i$ by identification condition (\ref{pi_orthog})
for $\theta _{i}^{\ast }$ under Assumption \ref{assum:ident}.c. Apply (\ref%
{p_pi_equiv}) to deduce $(\partial /\partial \beta _{(i)})\mathcal{L}%
([\delta _{(i)}^{\ast },\boldsymbol{0}_{k_{\theta }}])$ $=$ $0$ for each $i$%
, hence the sub-gradient $(\partial /\partial \delta )\mathcal{L}([\delta
_{(i)}^{\ast },\boldsymbol{0}_{k_{\theta }}])$ $=$ $0$. But this implies $%
\delta _{(i)}^{\ast }$ is for each $i$ the unique solution to the
optimization problem $\inf_{\delta \in \mathcal{D}}\mathcal{L}(\delta ,%
\boldsymbol{0}_{k_{\theta }})$, hence the $\delta _{(i)}^{\ast \prime }s$
are equivalent: there exists a unique $\delta ^{\ast }$ $\in $ $\mathcal{D}$
such that $\delta _{(i)}^{\ast }$ $=$ $\delta ^{\ast }$ for each $i$. This
yields each sub-gradient $(\partial /\partial \beta _{(i)})\mathcal{L}%
(\delta ^{\ast },\boldsymbol{0}_{k_{\theta }})$ $=$ $0$ for each $i,$ hence
the gradient $(\partial /\partial \beta )\mathcal{L}(\delta ^{\ast },%
\boldsymbol{0}_{k_{\theta }})$ $=$ $0$. Therefore $[\delta _{0},\theta _{0}]$
$=$ $[\delta ^{\ast },\boldsymbol{0}_{k_{\theta }}]$ by orthogonality and
uniqueness Assumption \ref{assum:ident}.b. $\mathcal{QED}$.\bigskip

Theorem \ref{th:ident} and its proof imply identification of the hypotheses
asymptotically when $k_{\theta ,n}$ parameterizations are used, provided $%
k_{\theta ,n}$ $\rightarrow $ $k_{\theta }$ cf. Assumption \ref{assum:ident}%
.a.

\begin{corollary}
\label{cor:ident_kn}Under Assumption \ref{assum:ident}, $\theta _{0}$ $=$ $%
\boldsymbol{0}_{k_{\theta }}$ \emph{if and only if} $|\theta ^{\ast }$ $-$ $%
\boldsymbol{0}_{k_{\theta ,n}}|$ $\rightarrow $ $0$ as $n$ $\rightarrow $ $%
\infty $, hence $\delta ^{\ast }$ $=$ $\delta _{0}$ as $n$ $\rightarrow $ $%
\infty $ under $H_{0}$ $:$ $\theta _{0}$ $=$ $\boldsymbol{0}_{k_{\theta }}$.
Under $H_{1}$ $:$ $\theta _{0}$ $\neq $ $\boldsymbol{0}_{k_{\theta }}$ it
follows $\lim_{n\rightarrow \infty }\max_{1\leq i\leq k_{\theta ,n}}|\theta
_{i}^{\ast }|$ $>$ $0$.
\end{corollary}

\section{Limit Theory for Parsimonious Model Estimators\label{sec:dist}}

This section presents core asymptotic theory for the parsimonious estimators 
$\hat{\beta}_{(i)}$. We prove under general conditions that $\hat{\beta}%
_{(i)}$ is consistent, and derive a required first order approximation.
Gaussian approximation and bootstrap theory are then presented in the
following section.

\subsection{Consistency\label{sec:consist}}

The following are standard conditions for consistency when the data are
non-trending \citep[see][]{PakesPollard1989}. In the case of stochastic or
deterministic trend, consistency can be proven by case using specific model
information and assumptions on the data generating process. Recall the
sample criterion for the $i^{th}$ parameterization is $\mathcal{\hat{L}}%
_{(i)}(\cdot )$, $k_{\theta ,n}$ is the number of parameterizations, and $%
k_{\theta }$ $\in $ $\mathbb{N}\cup \infty $ is the true dimension of $%
\theta _{0}$.

\begin{assumption}[Consistency]
\label{assum:consist_suff} \ \ \medskip \newline
$a$. Each estimator $\hat{\beta}_{(i)}$ satisfies $\mathcal{\hat{L}}_{(i)}(%
\hat{\beta}_{(i)})$ $=$ $\inf_{\beta _{(i)}\in \mathcal{B}_{(i)}}\{\mathcal{%
\hat{L}}_{(i)}(\beta _{(i)})\}$.$\medskip $\newline
$b$. $\sup_{\beta _{(i)}\in \mathcal{B}_{(i)}}|\mathcal{\hat{L}}_{(i)}(\beta
_{(i)})/n$ $-$ $\mathcal{L}_{(i)}(\beta _{(i)})|$ $\overset{p}{\rightarrow }$
$0$ for each $i$.
\end{assumption}

\begin{remark}
\normalfont Consistency can be proven under $\mathcal{\hat{L}}_{(i)}(\hat{%
\beta}_{(i)})$ $=$ $\inf_{\beta _{(i)}\in \mathcal{B}_{(i)}}\{\mathcal{\hat{L%
}}_{(i)}(\beta _{(i)})\}$ $+$ $o_{p}(1)$, a weaker version of $(a)$,
accounting for non-differentiable loss $\mathcal{\hat{L}}_{(i)}(\cdot )$,
and iterative numerical optimization procedures %
\citep[e.g.][]{NeweyMcFadden1994}. In order to characterize the allowed rate 
$k_{\theta ,n}$ $\rightarrow $ $\infty $ for our first order approximation,
however, we need to know the rate of approximation (under differentiability
we need to known the rate of $(\partial /\partial \beta _{(i)})\mathcal{\hat{%
L}}_{(i)}(\hat{\beta}_{(i)})$ $\overset{p}{\rightarrow }$ $0$). This will
depend on the numerical method, degree of loss non-smoothness, and method of
gradient approximation. Even with this information, the approximation rate
may still be unknowable beyond $o_{p}(1)$ 
\citep[see,
e.g.,][]{HongMahajanNekipelov2015}.
\end{remark}

\begin{remark}
\normalfont($a$) defines $\hat{\beta}_{(i)}$, while existence of $\hat{\beta}%
_{(i)}$ is shown in, e.g., \citet[Lemma
2]{Jennrich1969}. ($b$) can be verified when $\mathcal{\hat{L}}_{(i)}(\cdot
) $\ is a sufficiently smooth sample mean, under various dependence
properties.
\end{remark}

The proof of consistency is well known, hence we present it in %
\citet[Appendix C]{supp_mat_testmanyzeros}.

\begin{theorem}
\label{th:consist}Under Assumptions \ref{assum:ident}.c and \ref%
{assum:consist_suff}, $\hat{\beta}_{(i)}$ $\overset{p}{\rightarrow }$ $\beta
_{(i)}^{\ast }$ for each $1$ $\leq $ $i$ $\leq $ $\mathring{k}$ and any $%
\mathring{k}$ $\in $ $\mathbb{N}$, $\mathring{k}$ $\leq $ $k_{\theta }$.
Hence $\hat{\theta}_{i}$ $\overset{p}{\rightarrow }$ $0$ for each $1$ $\leq $
$i$ $\leq $ $\mathring{k}$, and any $\mathring{k}$ $\leq $ $k_{\theta }$,\ 
\emph{if and only if} $H_{0}$ is true.
\end{theorem}

\begin{remark}
\normalfont The theorem implies $\hat{\theta}_{i}$ $\overset{p}{\nrightarrow 
}$ $0$ for some $i$\ \emph{if and only if} $H_{0}$ is false, ensuring the
max-test is consistent.
\end{remark}

\begin{remark}
\normalfont In order to prove $\max_{1\leq i\leq k_{\theta ,n}}||\hat{\beta}%
_{(i)}$ $-$ $\beta _{(i)}^{\ast }||$ $\overset{p}{\rightarrow }$ $0$ we need
more information on loss and $\{k_{\theta ,n}\}$. See Lemma \ref%
{lm:beta_hat_rate} in Appendix \ref{app:proofs} where we show $\max_{1\leq
i\leq k_{\theta ,n}}||\hat{\beta}_{(i)}$ $-$ $\beta _{(i)}^{\ast }||$ $%
=O_{p}(\ln (k_{\theta ,n})/\sqrt{n})$ provided $k_{\theta ,n}=O(\sqrt{n}%
/(\ln (n))^{2})$, under the first order expansion setting laid out below. In
that setting, therefore, $\max_{1\leq i\leq k_{\theta ,n}}||\hat{\beta}%
_{(i)} $ $-$ $\beta _{(i)}^{\ast }||$ $\overset{p}{\rightarrow }$ $0$..
\end{remark}

\subsection{First Order Expansion\label{sec:first_order}}

Considering $\hat{\beta}_{(i)}$ $\overset{p}{\rightarrow }$ $\beta
_{(i)}^{\ast }$ and identification Corollary \ref{cor:ident_kn}\ ensure
asymptotic test consistency under $H_{1}$, we now work exclusively under $%
H_{0}$. Define gradient and Hessian functions:%
\begin{eqnarray}
&&\widehat{\mathcal{G}}_{(i)}(\beta _{(i)})\equiv \frac{\partial }{\partial
\beta _{(i)}}\mathcal{\hat{L}}_{(i)}(\beta _{(i)}),\text{ and set\ }\widehat{%
\mathcal{G}}_{(i)}=\widehat{\mathcal{G}}_{(i)}(\beta _{(i)}^{\ast })
\label{GH} \\
&&\widehat{\mathcal{H}}_{(i)}(\beta _{(i)})\equiv \frac{\partial ^{2}}{%
\partial \beta _{(i)}\partial \beta _{(i)}^{\prime }}\mathcal{\hat{L}}%
_{(i)}(\beta _{(i)}).  \notag
\end{eqnarray}%
The estimators $\hat{\beta}_{(i)}$ satisfy $\widehat{\mathcal{G}}_{(i)}(\hat{%
\beta}_{(i)})$ $=$ $0$ for each $n$ by construction. The mean value theorem
yields for for some sequence $\{\ddot{\beta}_{(i)}\}$, $||\ddot{\beta}_{(i)}$
$-$ $\beta _{(i)}^{\ast }||$ $\leq $ $||\hat{\beta}_{(i)}$ $-$ $\beta
_{(i)}^{\ast }||,$%
\begin{equation}
0=\widehat{\mathcal{G}}_{(i)}+\widehat{\mathcal{H}}_{(i)}(\ddot{\beta}%
_{(i)})\left( \hat{\beta}_{(i)}-\beta _{(i)}^{\ast }\right) .
\label{1st_expand}
\end{equation}%
Now define:%
\begin{equation}
\widehat{\mathfrak{H}}_{(i)}(\beta _{(i)})\equiv \frac{1}{n}\widehat{%
\mathcal{H}}_{(i)}(\beta _{(i)})\text{ and }\mathfrak{H}_{(i)}(\beta
_{(i)})\equiv \plim_{n\rightarrow \infty }\frac{1}{n}\widehat{\mathcal{H}}%
_{(i)}(\beta _{(i)})\text{ and }\mathfrak{H}_{(i)}=\mathfrak{H}_{(i)}(\beta
_{(i)}^{\ast }).  \label{H(b)}
\end{equation}%
Hence%
\begin{eqnarray}
\sqrt{n}\left( \hat{\beta}_{(i)}-\beta _{(i)}^{\ast }\right) &=&-\mathfrak{H}%
_{(i)}^{-1}\frac{1}{\sqrt{n}}\widehat{\mathcal{G}}_{(i)}  \label{b_expand} \\
&&-\left\{ \widehat{\mathfrak{H}}_{(i)}^{-1}(\ddot{\beta}_{(i)})-\mathfrak{H}%
_{(i)}^{-1}(\ddot{\beta}_{(i)})\right\} \frac{1}{\sqrt{n}}\widehat{\mathcal{G%
}}_{(i)}  \notag \\
&&-\left\{ \mathfrak{H}_{(i)}^{-1}(\ddot{\beta}_{(i)})-\mathfrak{H}%
_{(i)}^{-1}(\beta _{(i)}^{\ast })\right\} \frac{1}{\sqrt{n}}\widehat{%
\mathcal{G}}_{(i)}  \notag \\
&=&\mathcal{\hat{Z}}_{(i)}+\mathcal{\hat{R}}_{1,i}(\ddot{\beta}_{(i)})+%
\mathcal{\hat{R}}_{2,i}(\ddot{\beta}_{(i)}),  \notag
\end{eqnarray}%
say. $\widehat{\mathfrak{H}}_{(i)}(\cdot )$ and its limit $\mathfrak{H}%
_{(i)}(\cdot )$\ are assumed below to be positive definite in a neighborhood
of $\beta _{(i)}^{\ast }$, and $||\ddot{\beta}_{(i)}-\beta _{(i)}^{\ast }||$ 
$\overset{p}{\rightarrow }$ $0$ by construction and Theorem \ref{th:consist}%
. This ensures $\mathcal{\hat{Z}}_{(i)}$ and each $\mathcal{\hat{R}}_{j,i}$\
are well defined $awp1$.

In order to yield a class of sequences $\{k_{\theta ,n}\}$ that satisfies $%
\max_{1\leq i\leq k_{\theta ,n}}||\sqrt{n}(\hat{\beta}_{(i)}$ $-$ $\beta
_{(i)}^{\ast })$ $-$ $\mathcal{\hat{Z}}_{(i)}||\overset{p}{\rightarrow }0$
we need rates of convergence for $\mathcal{\hat{R}}_{\cdot ,i}(b)$ uniformly
in $i$ and $b$. The following provides a suitable setting. Let $\{\varpi
_{n}\}_{n\in \mathbb{N}}$ be a sequence of positive real numbers, $\varpi
_{n}$ $\rightarrow $ $0$, that may be different in different places; and $%
\mathcal{B}_{n,(i)}$ $\equiv $ $\{\mathcal{B}_{(i)}:||\beta _{(i)}-\beta
_{(i)}^{\ast }||$ $\leq $ $\varpi _{n}\},$ thus $\mathcal{B}_{n,(i)}$\
contains $\beta _{(i)}$ close to $\beta _{(i)}^{\ast }$. Let $\{k_{\theta
,n}\}$ be arbitrary positive integers.

\begin{assumption}[Asymptotic Expansion]
\label{assum:expand}Let $H_{0}$ hold.\medskip \newline
$a$. $||\widehat{\mathfrak{H}}_{(i)}(\beta _{(i)})$ $-$ $\widehat{\mathfrak{H%
}}_{(i)}(\beta _{(i)}^{\ast })||$ $\leq $ $\widehat{\mathcal{C}}%
_{(i)}||\beta _{(i)}$ $-$ $\beta _{(i)}^{\ast }||$ $\forall \beta _{(i)}$ $%
\in $ $\mathcal{B}_{n,(i)}$ for every $\{\varpi _{n}\}$ and some positive
stochastic $\widehat{\mathcal{C}}_{(i)}$ with $\max_{1\leq i\leq k_{\theta
,n}}\widehat{\mathcal{C}}_{(i)}$ $=$ $O_{p}(1)$. $||\mathfrak{H}_{(i)}(\beta
_{(i)})$ $-$ $\mathfrak{H}_{(i)}(\beta _{(i)}^{\ast })||$ $\leq $ $\mathcal{C%
}_{(i)}||\beta _{(i)}$ $-$ $\beta _{(i)}^{\ast }||$ $\forall \beta _{(i)}$ $%
\in $ $\{\mathcal{B}_{(i)}$ $:$ $||\beta _{(i)}$ $-$ $\beta _{(i)}^{\ast }||$
$\leq $ $\varepsilon \}$, some $\varepsilon $ $>$ $0$, and nonstochastic $%
\mathcal{C}_{(i)}$ $>$ $0$ with $\max_{i\in \mathbb{N}}\mathcal{C}_{(i)}$ $<$
$\infty $.$\medskip $\newline
$b$.

$(i)$ $\widehat{\mathfrak{H}}_{(i)}(\beta _{(i)})$ is symmetric, and
positive definite and bounded uniformly $awp1$. Specifically:%
\begin{eqnarray*}
&&\inf_{\lambda ^{\prime }\lambda =1}\min_{1\leq i\leq k_{\theta ,n}}\left\{
\inf_{\beta _{(i)}\in \mathcal{B}_{n,(i)}}\left\{ \lambda ^{\prime }\widehat{%
\mathfrak{H}}_{(i)}(\beta _{(i)})\lambda \right\} \right\} >0\text{ }awp1. \\
&&\max_{1\leq i\leq k_{\theta ,n}}\left\{ \sup_{\beta _{(i)}\in \mathcal{B}%
_{n,(i)}}\left\Vert \widehat{\mathfrak{H}}_{(i)}^{-1}(\beta
_{(i)})\right\Vert \right\} =O_{p}(k_{\theta ,n}).
\end{eqnarray*}

$(ii)$ For some $\zeta ,\xi >0$, 
\begin{equation*}
\max_{1\leq i\leq k_{\theta ,n}}E\left[ \exp \left\{ \zeta \sup_{\beta
_{(i)}\in \mathcal{B}_{n,(i)}}\sqrt{n}\left\vert \widehat{\mathfrak{H}}%
_{(i)}(\beta _{(i)})-\mathfrak{H}_{(i)}(\beta _{(i)})\right\vert \right\} %
\right] =O\left( n^{\xi \ln (k_{n})}\right) ,
\end{equation*}%
where $\mathfrak{H}_{(i)}(\beta _{(i)})$ is non-random, symmetric, and
uniformly positive definite on $\{\beta _{(i)}$ $:$ $||\beta _{(i)}$ $-$ $%
\beta _{(i)}^{\ast }||$ $\leq $ $\varepsilon \}$ for some $\varepsilon $ $>$ 
$0$. In particular $\max_{i\in \mathbb{N}}\{\sup_{\beta _{(i)}:||\beta
_{(i)}-\beta _{(i)}^{\ast }||\leq \varepsilon }||\mathfrak{H}%
_{(i)}^{-1}(\beta _{(i)})||\}$ $<$ $\infty .\medskip $\newline
$c$. $\max_{1\leq i\leq k_{\theta ,n}}E[\exp \{\zeta |\widehat{\mathcal{G}}%
_{(i)}/\sqrt{n}|\}]$ $=$ $O(n^{\xi \ln (k_{n})})$ for some $\zeta ,\xi >0.$
\end{assumption}

\begin{remark}
\normalfont Local Lipschitz ($a$) simplifies higher order asymptotics. Under
($b$) the asymptotic Hessian $\mathfrak{H}_{(i)}(\cdot )$ could be random
(permitting, e.g., stochastic trend). We assume it is non-random because in
Section (\ref{sec:max_test_dist}) we exploit a Gaussian approximation theory
on the expansion term $\mathcal{\hat{Z}}_{(i)}$ $=$ $-\mathfrak{H}_{(i)}^{-1}%
\widehat{\mathcal{G}}_{(i)}/\sqrt{n}$ in (\ref{b_expand}), where $\widehat{%
\mathcal{G}}_{(i)}/\sqrt{n}$ is asymptotically normal (hence $\mathfrak{H}%
_{(i)}$ generally cannot be random).
\end{remark}

\begin{remark}
\normalfont The Hessian and gradient exponential moment bounds ($b.ii$) and (%
$c$) are used to yield a non-Gaussian first order approximation theory based
on a \textit{smooth-max} technique. See Lemma \ref{lm:max_p-converg} in
Appendix \ref{app:proofs} and its proof in 
\citet[Appendix
C.2]{supp_mat_testmanyzeros}. Sub-Gaussian and sub-exponential distribution
classes provide a popular realm for developing high dimensional probability
theory like concentration bounds required here %
\citep[see][]{Vershynin2018,KuchibhotlaChakrabortty2020}.
\end{remark}

\begin{remark}
\normalfont The ($b.i$) bound $\mathcal{X}_{n}$ $\equiv $ $\max_{1\leq i\leq
k_{\theta ,n}}\{\sup_{\beta _{(i)}\in \mathcal{B}_{n,(i)}}||\widehat{%
\mathfrak{H}}_{(i)}^{-1}(\beta _{(i)})||\}$ $=$ $O_{p}(k_{\theta ,n})$ is
mild since $\mathcal{X}_{n}$ $\leq $ $\sum_{i=1}^{k_{\theta ,n}}\sup_{\beta
_{(i)}\in \mathcal{B}_{n,(i)}}||\widehat{\mathfrak{H}}_{(i)}^{-1}(\beta
_{(i)})||$ and $\sup_{\beta _{(i)}\in \mathcal{B}_{n,(i)}}||\widehat{%
\mathfrak{H}}_{(i)}^{-1}(\beta _{(i)})||$ $=$ $O_{p}(1)$ under suitable
dependence, heterogeneity and smoothness constraints. See also Lemma \ref%
{lm:max_p-converg2} in Appendix \ref{app:proofs}.
\end{remark}

The following provide low level sufficient conditions for the ($b.ii$) and ($%
c$) exponential bounds. See \citet[Appendix B.2]{supp_mat_testmanyzeros} for
expanded details and references. In\ brief, we assume $\widehat{\mathcal{G}}%
_{(i)}$ and $\widehat{\mathcal{H}}_{(i)}(\beta _{(i)})$ are partial sums of
independent random variables with sub-exponential tails, and therefore have
a moment generating function (covering bounded support, Gaussian, Laplace,
exponential, Poisson, and Gamma random variables). This ensures the summands
are locally sub-Gaussian, which ultimately suffices %
\citep[cf.][]{Chareka_etal_2006}.\medskip \newline
\textbf{Assumption \ref{assum:expand}.b(ii)}$^{\ast }$ \textit{Let }$%
\widehat{\mathfrak{H}}_{(i)}$\textit{\ }$\equiv $\textit{\ }$%
1/n\sum_{t=1}^{n}h_{i,t}(\beta _{(i)})$\textit{\ and }$\mathfrak{H}%
_{(i)}(\beta _{(i)})$ $\equiv $ $\lim_{n\rightarrow \infty
}1/n\sum_{t=1}^{n}E[h_{i,t}(\beta _{(i)})]$\textit{, where }$\{h_{i,t}(\cdot
)$ $=$ $[h_{l,m,i,t}(\cdot )]_{l,m=1}^{k_{\delta }+1}\}$\textit{\ are
measurable random variables\ on a common probability space, independent
across }$t$\textit{. For some finite universal constants }$(\delta ,\mathcal{%
K})$ $>$ $0$:\textit{\ }%
\begin{equation*}
\max_{1\leq l,m\leq k_{\delta }+1}\max_{i,t\in N}P\left( \sup_{\beta
_{(i)}\in \mathcal{B}_{n,(i)}}\left\vert h_{l,m,i,t}(\beta _{(i)})-\mathfrak{%
H}_{l,m,(i)}(\beta _{(i)})\right\vert >c\right) \leq \mathcal{K}e^{-\delta c}%
\text{ }\forall c>0.
\end{equation*}%
\textbf{Assumption \ref{assum:expand}.c}$^{\ast }$ \textit{Let }$\widehat{%
\mathcal{G}}_{(i)}$\textit{\ }$\equiv $\textit{\ }$\sum_{t=1}^{n}\epsilon
_{i,t}w_{i,t}$\textit{, }$w_{i,t}$ $=$ $[w_{l,i,t}]_{l=1}^{k_{\delta }+1}$, 
\textit{where }$(\epsilon _{i,t},w_{i,t})$\textit{\ are random variables\ on
a common probability space, independent across }$t$\textit{, }$E[\epsilon
_{i,t}]$ $=$ $0$\textit{, and }$\epsilon _{i,s}$\textit{\ and }$w_{l,i,t}$%
\textit{\ are mutually independent }$\forall i,l,s,t$\textit{. Further, }$%
\max_{i,t\in N}P(|\epsilon _{i,t}|$ $>$ $c)$ $\leq $ $\mathcal{K}e^{-\delta
c}$ \textit{and} $\max_{i,t\in N}P(|w_{i,t}$ $-$ $E[w_{i,t}]|$ $>$ $c)$ $%
\leq $ $\mathcal{K}e^{-\delta c}\ \forall c$ $>$ $0$ \textit{for finite
constants }$(\delta ,\mathcal{K})$ $>$ $0$ \textit{that may be different in
different places}.

\begin{remark}
\normalfont The ($c^{\ast }$) form $\widehat{\mathcal{G}}_{(i)}$\textit{\ }$%
= $\textit{\ }$\sum_{t=1}^{n}\epsilon _{i,t}w_{i,t}$\ includes M-estimator
loss with, e.g., regression model error $\epsilon _{i,t}$ $=$ $\epsilon _{t}$
and response gradient $w_{i,t}$; and moment-based estimators like GMM, with
identification \ condition $E[m_{(i),t}(\beta _{(i)}^{\ast })]$ $=$ $0$,
hence $\epsilon _{i,t}$ $=$ $m_{(i),t}(\beta _{(i)}^{\ast })$ and $w_{i,t}$
is a multiple of $(\partial /\partial \beta _{(i)})m_{(i),t}(\beta
_{(i)}^{\ast })$).
\end{remark}

\begin{remark}
\normalfont Verification of Assumptions \ref{assum:expand}.b(ii),c for
finite dependent random variables is analogous. Expanding the result to
mixing (or similar) sequences, and therefore demonstrating a partial sum of
dependent sub-exponential random variables has a moment generating function,
is left for future work.
\end{remark}

\begin{lemma}
\label{lm:Assum3bc*}Assumptions \ref{assum:expand}.b(ii)$^{\ast }$,c$^{\ast
} $ suffice for Assumptions \ref{assum:expand}.b(ii),c for any sequence of
positive integers $\{k_{\theta ,n}\}$.
\end{lemma}

\noindent\textbf{Proof.} Corollaries B.2 and B.3 and the surrounding
discussion in \citet[Appendix B.2]{supp_mat_testmanyzeros} prove the claim. $%
\mathcal{QED}$.\medskip

\begin{remark}
\normalfont Corollary B.2 in \cite{supp_mat_testmanyzeros} works with $%
\widehat{\mathcal{G}}_{(i)}$\textit{\ }$\equiv $\textit{\ }$%
\sum_{t=1}^{n}\epsilon _{i,t}w_{i,t}$\textit{. }It does not require
independence for $w_{i,t}$ over $t$, hence the result covers some time
series settings where $w_{i,t}$ is strictly exogenous (e.g. distributed lag
models).
\end{remark}

We now have the following key approximation result.

\begin{lemma}
\label{lm:expansion}Let $H_{0}$, and Assumptions \ref{assum:consist_suff}-%
\ref{assum:expand} hold. Let $\{k_{\theta ,n}\}$ be any monotonic sequence
of positive integers satisfying $k_{\theta ,n}$ $=$ $o(\sqrt{n}/(\ln
(n))^{4})$, and let the weight sequences $\{\mathcal{W}_{n,i}\}_{n\in 
\mathbb{N}}$ satisfy $\max_{1\leq i\leq k_{\theta ,n}}|\mathcal{W}_{n,i}$ $-$
$\mathcal{W}_{i}|$ $=$ $o_{p}(1/(\ln (n))^{2})$ for non-stochastic $\mathcal{%
W}_{i}$ $\in $ $(0,\infty )$. Then: 
\begin{equation}
\left\vert \max_{1\leq i\leq k_{\theta ,n}}\left\vert \sqrt{n}\mathcal{W}%
_{n,i}\hat{\theta}_{i}\right\vert -\max_{1\leq i\leq k_{\theta
,n}}\left\vert \mathcal{W}_{i}[\boldsymbol{0}_{k_{\delta }}^{\prime },1]%
\mathcal{\hat{Z}}_{(i)}\right\vert \right\vert \overset{p}{\rightarrow }0.
\label{theta_approx}
\end{equation}
\end{lemma}

\begin{remark}
\normalfont$\{\mathcal{W}_{n,i}\}_{n\in \mathbb{N}}$ must converge fast
enough in high dimension in order to show (\ref{theta_approx}) holds when $%
\max_{1\leq i\leq k_{\theta ,n}}|\sqrt{n}\hat{\theta}_{i}$ $-$ $[\boldsymbol{%
0}_{k_{\delta }}^{\prime },1]\mathcal{\hat{Z}}_{(i)}|$ $\overset{p}{%
\rightarrow }$ $0$ by a standard triangle inequality argument. The proof
shows $k_{\theta ,n}$ $=$ $o(\sqrt{n}/(\ln (n))^{4})$ is not improved if we
merely assume $\mathcal{W}_{n,i}$ $=$ $\mathcal{W}_{i}$ (e.g. $\mathcal{W}%
_{n,i}$ $=$ $1$). See arguments following the non-Gaussian approximation
requirement (\ref{maxtheta_Z}).
\end{remark}

\begin{remark}
\normalfont We require $k_{\theta ,n}$ $=$ $o(\sqrt{n}/(\ln (n))^{4})$.
Since test consistency requires $k_{\theta ,n}$ $\rightarrow $ $\infty $
when $k_{\theta }$ $=$ $\infty $, without much loss of generality, and to
ease notation,\ we use in some cases:%
\begin{equation}
k_{\theta ,n}\varpropto n^{\kappa }\text{ for any }\kappa \in (0,1/2)
\label{Kn}
\end{equation}
\end{remark}

\section{Max-Test and p-Value Computation\label{sec:pvalue}}

Our goal is to use the max-statistic $\mathcal{T}_{n}$ $=$ $\max_{1\leq
i\leq k_{\theta ,n}}|\sqrt{n}\mathcal{W}_{n,i}\hat{\theta}_{i}|$ in order to
compute a p-value approximation $\hat{p}_{n}$ $\in $ $[0,1]$ that leads to
an asymptotically correctly sized test: $P(\hat{p}_{n}$ $<$ $\alpha )$ $%
\rightarrow $ $\alpha $ under $H_{0}$ for any level $\alpha $ $\in $ $(0,1)$%
. The max-test is therefore:%
\begin{equation}
\text{reject }H_{0}\text{ at level }\alpha \in \left( 0,1\right) \text{ when 
}\hat{p}_{n}<\alpha \text{.}  \label{max_test}
\end{equation}

\subsection{Max-Test Gaussian Approximation\label{sec:max_test_dist}}

Recall $\mathcal{\hat{Z}}_{(i)}$ $=$ $-\mathfrak{H}_{(i)}^{-1}\widehat{%
\mathcal{G}}_{(i)}/\sqrt{n}$ with non-random $\mathfrak{H}_{(i)}(\beta
_{(i)})$ $\equiv $ $\plim_{n\rightarrow \infty }\widehat{\mathcal{H}}%
_{(i)}(\beta _{(i)})/n$. Lemma \ref{lm:expansion} implies that for
max-statistic asymptotics under the null we need a Gaussian approximation
result for $\max_{1\leq i\leq k_{\theta ,n}}|[\boldsymbol{0}_{k_{\delta
}}^{\prime },1]\mathcal{\hat{Z}}_{(i)}|$.

Define a long run covariance kernel (Assumption \ref{assum:max_dist}.a below
implies existence):%
\begin{equation*}
\mathcal{S}(i,j)=\lim_{n\rightarrow \infty }\mathcal{S}_{n}(i,j)\text{ where 
}\mathcal{S}_{n}(i,j)\equiv \frac{1}{n}E\left[ \widehat{\mathcal{G}}_{(i)}%
\widehat{\mathcal{G}}_{(j)}^{\prime }\right] .
\end{equation*}%
Define for arbitrary $\lambda $ $\in $ $\mathbb{R}^{k_{\delta }+1}$, $%
\lambda ^{\prime }\lambda $ $=$ $1$: 
\begin{eqnarray}
&&\mathcal{V}_{n,(i)}\equiv \mathfrak{H}_{(i)}^{-1}\mathcal{S}_{n}(i,i)%
\mathfrak{H}_{(i)}^{-1}\text{ \ and \ }\mathcal{V}_{(i)}\equiv \mathfrak{H}%
_{(i)}^{-1}\mathcal{S}(i,i)\mathfrak{H}_{(i)}^{-1}  \label{ZVsig} \\
&&\sigma _{n,(i)}^{2}(\lambda )\equiv \lambda ^{\prime }\mathcal{V}%
_{n,(i)}\lambda \text{ \ and }\sigma _{i}^{2}(\lambda )\equiv \lambda
^{\prime }\mathcal{V}_{(i)}\lambda .  \notag
\end{eqnarray}%
By construction $\sigma _{n,(i)}^{2}(\lambda )$ $=$ $E[(\lambda ^{\prime }%
\mathcal{\hat{Z}}_{(i)})^{2}]$. We specifically care about the case $\lambda 
$ $=$ $[\boldsymbol{0}_{k_{\delta }}^{\prime },1]^{\prime }$ in view of
expansion (\ref{theta_approx}), hence we only require a pointwise limit
theory vis-\`{a}-vis $\lambda $.

The main result of this section delivers a class of sequences $\{k_{\theta
,n}\}$, and an array of scalar normal random variables $\{\boldsymbol{Z}%
_{n,(i)}(\lambda )$ $:$ $n$ $\in $ $\mathbb{N}\}_{i=1}^{\infty }$, $%
\boldsymbol{Z}_{n,(i)}(\lambda )$ $\sim $ $N(0,\sigma _{n,(i)}^{2}(\lambda
)) $, such that for each $\lambda $ $\in $ $\mathbb{R}^{k_{\delta }+1}$, $%
\lambda ^{\prime }\lambda $ $=$ $1$, the Kolmogorov distance 
\begin{equation}
\rho _{n}(\lambda )\equiv \sup_{z\geq 0}\left\vert P\left( \max_{1\leq i\leq
k_{\theta ,n}}\left\vert \lambda ^{\prime }\mathcal{\hat{Z}}%
_{(i)}\right\vert \leq z\right) -P\left( \max_{1\leq i\leq k_{\theta
,n}}\left\vert \boldsymbol{Z}_{n,(i)}(\lambda )\right\vert \leq z\right)
\right\vert \rightarrow 0.  \label{Gauss_apprx}
\end{equation}%
The approximation does not require standardized $\lambda ^{\prime }\mathcal{%
\hat{Z}}_{(i)}$ and $\boldsymbol{Z}_{n,(i)}(\lambda )$ because we assume
below $\sigma _{n,(i)}^{2}(\lambda )$ lies in a compact subset of $(0,\infty
)$ asymptotically uniformly in $i$ and $\lambda $.

We could work in a broad dependence setting to allow for possibly
non-stationary mixing random variables, as in \cite{ChangChenWu2021}. A
physical dependence setting like \cite{ZhangWu2017} generally will not work
for a general model like (\ref{model}) since the dependence property does
not carry over to any measurable function of underlying random variables.
However, as noted above, we only verify the Assumptions \ref{assum:expand}%
.b(ii),c exponential moment bounds for independent random variables. We
therefore work in the independence setting of %
\citet[(E.1)]{Chernozhukov_etal2013}.

Now assume loss gradient $\widehat{\mathcal{G}}_{(i)}$ has the form%
\begin{equation}
\widehat{\mathcal{G}}_{(i)}=\sum_{t=1}^{n}G_{i,t}  \label{Gg}
\end{equation}%
where $\{G_{i,t}\}_{t=1}^{n}$ for each $i$\ are properly defined zero mean
independent random variables. The long run covariance kernel becomes: 
\begin{equation*}
\mathcal{S}(i,j)=\lim_{n\rightarrow \infty }\frac{1}{n}%
\sum_{t=1}^{n}E[G_{i,t}G_{j,t}^{\prime }],
\end{equation*}%
and the weighted first order process becomes:%
\begin{equation}
\lambda ^{\prime }\mathcal{\hat{Z}}_{(i)}=\frac{1}{\sqrt{n}}%
\sum_{t=1}^{n}z_{(i),t}(\lambda )\text{ where }z_{(i),t}(\lambda )=-\lambda
^{\prime }\mathfrak{H}_{(i)}^{-1}G_{i,t}.  \label{Zhg}
\end{equation}%
Assumption \ref{assum:expand}.b(ii) ensures $\mathfrak{H}_{(i)}^{-1}$
exists, hence $\lambda ^{\prime }\mathcal{\hat{Z}}_{(i)}$ is well defined.
Observe by independence:%
\begin{equation*}
\sigma _{n,(i)}^{2}(\lambda )=E\left[ (\lambda ^{\prime }\mathcal{\hat{Z}}%
_{(i)})^{2}\right] =\frac{1}{n}\sum_{t=1}^{n}\lambda ^{\prime }\mathfrak{H}%
_{(i)}^{-1}E\left[ G_{i,t}G_{i,t}^{\prime }\right] \mathfrak{H}%
_{(i)}^{-1}\lambda =\lambda ^{\prime }\mathcal{V}_{n,(i)}\lambda .
\end{equation*}

\begin{assumption}[Limit Distribution]
\label{assum:max_dist}Let $H_{0}$ hold, and let $\max_{i,t\in \mathbb{N}}||%
\mathfrak{H}_{(i)}^{-1}||$ $<$ $\infty $.\medskip \newline
$a$. $G_{i,t}$ is independent across $t$, and $\max_{i,t\in \mathbb{N}%
}||E[G_{i,t}G_{i,t}^{\prime }]||$ $<$ $\infty $.$\medskip $\newline
$b.$ Let $(i)\lim \inf_{n\rightarrow \infty }\min_{i\in \mathbb{N}%
}\inf_{\lambda ^{\prime }\lambda =1}\sigma _{n,(i)}^{2}(\lambda )$ $\geq $ 
\b{c} and $(ii)$ $\lim \sup_{n\rightarrow \infty }\max_{i\in \mathbb{N}%
}\sup_{\lambda ^{\prime }\lambda =1}\sigma _{n,(i)}^{2}(\lambda )$ $\leq $ $%
\bar{c}$ for some $0$ $<\text{\b{c}}\leq $ $\bar{c}$ $<$ $\infty $.$\medskip 
$\newline
$c.$ $\max_{\gamma =1,2}\{E|z_{(i),t}(\lambda )|^{2+\gamma }/B_{n}^{\gamma
}\}$ $+$ $E[\exp \{|z_{(i),t}(\lambda )|/B_{n}\}]$ $\leq $ $4$ uniformly in $%
i,t$ $\in $ $\mathbb{N}$ and $\lambda ^{\prime }\lambda $ $=$ $1$, for some
sequence of positive non-random numbers $\{B_{n}\}$, $B_{n}$ $\geq $ $1$
where $B_{n}$ $\rightarrow $ $\infty $ is possible.
\end{assumption}

\begin{remark}
\normalfont The Hessian and sub-exponential tail properties in Assumptions %
\ref{assum:expand}.b(ii) and \ref{assum:max_dist}.a together imply $%
\max_{i,t\in \mathbb{N}}||z_{(i),t}(\lambda )||_{2}$ $<$ $\infty $.
\end{remark}

\begin{remark}
\normalfont($b$) and ($c$) include sub-exponential and sub-Gaussian random
variables \citep[Comment 2.2]{Chernozhukov_etal2013}. Indeed, for ease of
notation suppose $G_{i,t}$ are scalar iid, sub-exponential uniformly in $i$
and $t$. Then $\max_{i,t\in \mathbb{N}}E[z_{(i),t}^{4}(\lambda )]$ $<$ $%
\infty $ and $\max_{i,t\in \mathbb{N}}E[\exp \{|z_{(i),t}(\lambda )|/K\}]\}$ 
$\leq $ $2$ for some $K$ $>$ $0$ \citep[Proposition 2.7.1]{Vershynin2018},
hence ($c$) holds for some large $B_{n}$ $=$ $B$.
\end{remark}

\begin{remark}
\normalfont Under Assumption \ref{assum:expand}.c$^{\ast }$ $G_{i,t}$ $=$ $%
\epsilon _{i,t}w_{i,t}$ is independent over $t$, and has uniformly
sub-exponential tails: $\max_{i,t\in \mathbb{N}}P(|\epsilon _{i,t}w_{i,t}|$ $%
>$ $c)$ $\leq $ $\exp \{-\delta c\}$ for some finite $\delta $ $>$ $0$ (see
the proof of Lemma B.1, and the subsequent remarks, in \cite%
{supp_mat_testmanyzeros}). Coupled with independence and $\max_{i,t\in 
\mathbb{N}}||\mathfrak{H}_{(i)}^{-1}||$ $<$ $\infty $ under Assumptions \ref%
{assum:expand}.b(ii), it follows that ($a$) and the upper bound on $\sigma
_{n,(i)}^{2}(\lambda )$ in ($b$) hold, and ($c$) holds by the arguments in
the previous remark. Finally, by mutual independence $\sigma
_{n,(i)}^{2}(\lambda )$ $=$ $1/n\sum_{t=1}^{n}E\left[ \epsilon _{i,t}^{2}%
\right] E[(\lambda ^{\prime }\mathfrak{H}_{(i)}^{-1}w_{i,t})^{2}]$ hence the
lower bound in ($b$) rests on bounding $E[\epsilon _{i,t}^{2}]$ and $%
E[(\lambda ^{\prime }w_{i,t})^{2}]$.
\end{remark}

In view of the preceding remarks, the following provides lower level
sufficient conditions for Assumption \ref{assum:max_dist}.\medskip \newline
\textbf{Assumption \ref{assum:max_dist}}$^{\ast }$.\textit{\ Let }$H_{0}$%
\textit{\ hold, let }$\lambda $\textit{\ }$\neq $\textit{\ }$0$\textit{, and
let }$\max_{i,t\in \mathbb{N}}||\mathfrak{H}_{(i)}^{-1}||$ $<$ $\infty $.$%
\medskip $\newline
$a.$ $G_{i,t}$\textit{\ }$=$\textit{\ }$\epsilon _{i,t}w_{i,t}$\textit{, }$%
w_{i,t}$ $=$ $[w_{l,i,t}]_{l=1}^{k_{\delta }+1}$\textit{, where }$(\epsilon
_{i,t},w_{i,t})$\textit{\ are independent over }$t$\textit{, }$E[\epsilon
_{i,t}]$\textit{\ }$=$\textit{\ }$0$\textit{, and }$\epsilon _{i,s}$\textit{%
\ and }$w_{l,i,t}$\textit{\ are mutually independent }$\forall i,l,s,t$%
\textit{. Further, }$\max_{t\in N}P(|\epsilon _{i,t}|$\textit{\ }$>$\textit{%
\ }$c)$\textit{\ }$\leq $\textit{\ }$Ke^{-\delta c}$\textit{\ and }$%
\max_{i,t\in N}P(|w_{i,t}$\textit{\ }$-$\textit{\ }$E[w_{i,t}]|$\textit{\ }$%
> $\textit{\ }$c)$\textit{\ }$\leq $\textit{\ }$\mathcal{K}e^{-\delta c}\
\forall c$\textit{\ }$>$\textit{\ }$0$\textit{\ for finite constants }$%
(\delta ,\mathcal{K})$\textit{\ }$>$\textit{\ }$0$ \textit{that may be
different in different places}.$\medskip $\newline
$b.$ \textit{Let }$\min_{i,t\in \mathbb{N}}\min \{E[\epsilon
_{i,t}^{2}],\inf_{\lambda ^{\prime }\lambda =1}E[(\lambda ^{\prime
}w_{i,t})^{2}]\}$ $\geq $ \b{c} \textit{for some }$\text{\b{c}}$ $\in $ $%
(0,\infty )$.

\begin{lemma}
\label{lm:suff_Assum4abiic}Assumption \ref{assum:max_dist}$^{\ast }$
suffices for Assumption \ref{assum:max_dist}, in particular $B_{n}$ $=$ $B$
for some $B$ $\in $ $(0,\infty )$.
\end{lemma}

Gaussian approximation (\ref{Gauss_apprx}) is now formally stated.

\begin{lemma}
\label{lm:max_dist}Let $H_{0}$ and Assumption \ref{assum:max_dist} hold.
Then $\rho _{n}(\lambda )$ $=$ $O(n^{-c})$ for some $c$ $>$ $0$ and any $%
\{k_{\theta ,n}\}$ that satisfies $B_{n}^{2}(\ln (k_{\theta ,n}n))^{7}$ $=$ $%
O(n^{1-\xi })$ for some $\xi $ $>0$.
\end{lemma}

Under (\ref{Kn}) with $k_{\theta ,n}$ $\varpropto $ $n^{\kappa }$ for any $%
\kappa $ $\in $ $(0,1/2)$, the bound $B_{n}^{2}(\ln (k_{\theta ,n}n))^{7}/n$ 
$=$ $O(n^{-\xi })$ reduces to $B_{n}^{2}(\ln (n))^{7}$ $=$ $O(n^{1-\xi })$,
which automatically holds when $B_{n}$ $=$ $O(n^{a})$ for some $a$ $\in $ $%
[0,1/2)$ and any $\xi $ $<$ $1$ $-$ $2a$. As an example, the latter
condition holds when $z_{(i),t}$ is sub-exponential uniformly in $i$ and $t$
since then $\max_{i,t\in \mathbb{N}}E[z_{(i),t}^{4}(\lambda )]$ $<$ $\infty $
and $E[\exp \{|z_{(i),t}(\lambda )|/K\}]\}$ $\leq $ $2$ for some $K$ $>$ $0$%
, hence $a$ $=$ $0$. This in turn holds under Assumption \ref{assum:max_dist}%
$^{\ast }$ in which case $z_{(i),t}(\lambda )$ $=$ $-\lambda ^{\prime }%
\mathfrak{H}_{(i)}^{-1}\epsilon _{i,t}w_{i,t}$, yielding the following
corollary.

\begin{corollary}
\label{cor:max_dist}Let $H_{0}$ and Assumption \ref{assum:max_dist}$^{\ast }$
hold, and let $k_{\theta ,n}$ $\varpropto $ $n^{\kappa }$ for any $\kappa $ $%
\in $ $(0,1/2)$. Then $\rho _{n}(\lambda )$ $=$ $O(n^{-c})$ for some $c$ $>$ 
$0$.
\end{corollary}

Expansion (\ref{theta_approx}) and Gaussian approximation (\ref{Gauss_apprx}%
) yield the following main result.

\begin{theorem}
\label{th:max_theta_hat}Let Assumptions \ref{assum:ident}-\ref%
{assum:max_dist} hold, let $\{k_{\theta ,n}\}$ be any monotonic sequence of
positive integers satisfying $k_{\theta ,n}$ $\varpropto $ $n^{\kappa }$
with $\kappa $ $\in $ $(0,1/2)$, and assume the weight sequences satisfy $%
\max_{1\leq i\leq k_{\theta ,n}}|\mathcal{W}_{n,i}$ $-$ $\mathcal{W}_{i}|$ $%
= $ $o_{p}(1/(\ln (n))^{2})$ for non-stochastic $\mathcal{W}_{i}$ $\in $ $%
(0,\infty )$ and $\max_{i\in \mathbb{N}}\mathcal{W}_{i}$ $<$ $\infty $.$%
\medskip $\newline
$a.$ Under $H_{0}$, for the Gaussian random variables $\boldsymbol{Z}%
_{n,(i)}(\lambda )$ $\sim $ $N(0,\sigma _{n,(i)}^{2}(\lambda ))$ in (\ref%
{Gauss_apprx}),%
\begin{equation}
\left\vert P\left( \max_{1\leq i\leq k_{\theta ,n}}\left\vert \sqrt{n}%
\mathcal{W}_{n,i}\hat{\theta}_{i}\right\vert \leq z\right) -P\left(
\max_{1\leq i\leq k_{\theta ,n}}\left\vert \mathcal{W}_{i}\boldsymbol{Z}%
_{n,(i)}([\boldsymbol{0}_{k_{\delta }}^{\prime },1])\right\vert \leq
z\right) \right\vert \rightarrow 0\text{ }\forall z\geq 0,  \label{PP0}
\end{equation}%
provided the Assumption \ref{assum:max_dist}.c sequence $\{B_{n}\}$
satisfies $B_{n}$ $=$ $O(n^{a})$ for some $a$ $\in $ $[0,1/2)$. In that case 
$\max_{1\leq i\leq k_{\theta ,n}}|\sqrt{n}\mathcal{W}_{n,i}\hat{\theta}_{i}|$
$\overset{d}{\rightarrow }$ $\max_{i\in \mathbb{N}}|\mathcal{W}_{i}%
\boldsymbol{Z}_{(i)}([\boldsymbol{0}_{k_{\delta }}^{\prime },1])|$ where $%
\boldsymbol{Z}_{(i)}(\lambda )$ $\sim $ $N(0,\sigma _{(i)}^{2}(\lambda ))$
and $\sigma _{(i)}^{2}(\lambda )$ $=$ $\lim_{n\rightarrow \infty }\sigma
_{n,(i)}^{2}(\lambda )$.\medskip \newline
$b.$ Under $H_{1}$, $\max_{1\leq i\leq k_{\theta ,n}}|\sqrt{n}\mathcal{W}%
_{n,i}\hat{\theta}_{i}|$ $\rightarrow $ $\infty $ for any monotonic sequence
of positive integers $\{k_{\theta ,n}\}$, $k_{\theta ,n}$ $\rightarrow $ $%
k_{\theta }$.
\end{theorem}

\begin{remark}
\normalfont Gaussian approximation theory requires $\mathcal{W}_{i}%
\boldsymbol{Z}_{n,(i)}(\cdot )$ to be Gaussian, hence $\mathcal{W}_{i}$\ is
assumed to be non-random.
\end{remark}

\subsection{Parametric Bootstrap\label{sec:param_boot}}

We work with the nonlinear regression model $y_{t}$ $=$ $f(x_{t},\beta _{0})$
$+$ $\epsilon _{t}$ in (\ref{model}) with $\beta _{0}$ $=$ $\left[ \delta
_{0}^{\prime },\theta _{0}^{\prime }\right] ^{\prime }$ for simplicity,
where $\{x_{t},\epsilon _{t}\}$ are independent over $t$ \ and mutually
independent. Squared error loss is used also for notational ease. See
Assumption \ref{assum:param_boot} below for details.

Let the parsimonious response functions $f_{(i)}(x_{t},\beta _{(i)})$ be $%
f(x_{t},\beta )$ $=$ $f(x_{t},\delta ,\theta )$ when the $i^{th}$ element of 
$\theta $ is $\theta _{i}$ and the remaining $\theta _{j}$ $=$ $0$ for $j$ $%
\neq $ $i$:%
\begin{equation}
f_{(i)}(x_{t},\beta _{(i)})=f(x_{t},\delta ,[0,...,\theta _{i},0,,,]^{\prime
}).  \label{f(i)f}
\end{equation}%
This implies link (\ref{LiL}) with squared error loss, and therefore (\ref%
{p_pi_equiv}). Under twice differentiability imposed below, define:%
\begin{equation*}
g_{(i)}(x_{t},\beta _{(i)})\equiv \frac{\partial }{\partial \beta _{(i)}}%
f_{(i)}(x_{t},\beta _{(i)})\text{ and }h_{(i)}(x_{t},\beta _{(i)})\equiv
\left( \frac{\partial }{\partial \beta _{(i)}}\right)
^{2}f_{(i)}(x_{t},\beta _{(i)}).
\end{equation*}%
The parsimonious models are 
\begin{equation*}
y_{t}=f_{(i)}(x_{t},\beta _{(i)}^{\ast })+v_{(i),t}\text{ where }\beta
_{(i)}^{\ast }=[\delta _{(i)}^{\ast \prime },\theta _{i}^{\ast }]^{\prime }%
\text{ uniquely satisfy }E\left[ v_{(i),t}g_{(i)}(x_{t},\beta _{(i)}^{\ast })%
\right] =0.
\end{equation*}

\subsubsection{Bootstrap Method and Theory}

Consider a variant of Gon\c{c}alves and Kilian's (%
\citeyear{Goncalves_Kilian_2004}) fixed design wild (multiplier) bootstrap.
They follow \cite{Kreiss1997} for autoregressions, and the method is related
to a sub-sampling non-wild bootstrap method in \cite{Bose1988}. The method
is asymptotically equivalent to a Gaussian wild bootstrap applied to the
first order approximation process.

Define the restricted estimator $\hat{\beta}^{(0)}$ $\equiv $ $[\hat{\delta}%
^{(0)\prime },\boldsymbol{0}_{k_{\theta }}^{\prime }]^{\prime }$ under $%
H_{0} $ $:$ $\theta _{0}$ $=$ $0$, where $\hat{\delta}^{(0)}$ uniquely
minimizes the restricted criterion $\sum_{t=1}^{n}\{y_{t}$ $-$ $%
f(x_{t},[\delta ,\boldsymbol{0}_{k_{\theta }}])\}^{2}$ on $\mathcal{D}$.
Define residuals $\epsilon _{n,t}^{(0)}$ $\equiv $ $y_{t}$ $-$ $f(x_{t},\hat{%
\beta}^{(0)}).$ Draw iid $\{\eta _{t}\}_{t=1}^{n}$ from $N(0,1)$, and
generate the array $y_{n,t}^{\ast }$ $\equiv $ $f(x_{t},\hat{\beta}^{(0)})$ $%
+$ $\epsilon _{n,t}^{(0)}\eta _{t}$. Construct $k_{\theta ,n}$ regression
models $y_{n,t}^{\ast }$ $=$ $f_{(i)}(x_{t},\beta _{(i)})$ $+$ $v_{n,(i),t}$%
, and let $\widehat{\tilde{\beta}}_{(i)}$ $=$ $[\widehat{\tilde{\delta}}%
_{(i)}^{\prime },\widehat{\tilde{\theta}}_{i}]^{\prime }$ be the least
squares estimator of $\beta _{(i)}$ on $\mathcal{B}_{(i)}$ $\equiv $ $%
\mathcal{D}$ $\times $ $\Theta _{i}$. The bootstrapped max-statistic is 
\begin{equation}
\mathcal{\tilde{T}}_{n}\equiv \max_{1\leq i\leq k_{\theta ,n}}\left\vert 
\sqrt{n}\mathcal{W}_{n,i}\widehat{\tilde{\theta}}_{i}\right\vert .
\label{T_boot}
\end{equation}%
Repeat the above steps $\mathcal{M}$ times, each time drawing a new iid
sequence $\{\eta _{j,t}\}_{t=1}^{n}$, $j$ $=$ $1...\mathcal{M}$. This
results in a sequence of bootstrapped estimators $\{\widehat{\tilde{\theta}}%
_{j}\}_{j=1}^{\mathcal{M}}$ and test statistics $\{\mathcal{\tilde{T}}%
_{n,j}\}_{j=1}^{\mathcal{M}}$ that are iid conditional on the sample $%
\{x_{t},y_{t}\}_{t=1}^{n}$. The approximate p-value is 
\begin{equation}
\tilde{p}_{n,\mathcal{M}}\equiv \frac{1}{\mathcal{M}}\sum_{j=1}^{\mathcal{M}%
}I(\mathcal{\tilde{T}}_{n,j}>\mathcal{T}_{n}).  \label{p_boot}
\end{equation}

The above algorithm varies from \cite{Goncalves_Kilian_2004} since they
operate on the \textit{unrestricted} estimator in order to generate
residuals. We impose the null hypothesis when we estimate $\beta _{0}$: this
reduces dimensionality without affecting asymptotics under the null, \textit{%
and} leads to a consistent test under the alternative.

We require some low level assumptions. Write $\beta ^{(0)}$ $=$ $[\delta
^{(0)\prime },\boldsymbol{0}_{k_{\theta }}^{\prime }]^{\prime }$\ and $\beta
_{(i)}^{(0)}$ $=$ $[\delta ^{(0)\prime },0]^{\prime }$, where $\delta ^{(0)}$
minimizes $E[(y_{t}$ $-$ $f(x_{t}[\delta ,\boldsymbol{0}_{k_{\theta
}}]))^{2}]$ $\forall t$.

\begin{assumption}[Bootstrap]
\label{assum:param_boot} \ \ \ \medskip \newline
$a.$ \emph{Response Function}:$\ y_{t}$ $=$ $f(x_{t},\beta _{0})$ $+$ $%
\epsilon _{t}$, $\{x_{t},\epsilon _{t}\}$ are independent over $t$ \ and
mutually independent$,$ $E[\epsilon _{t}]$ $=$ $0$, $E[\epsilon _{t}^{2}]$ $%
\in $ $(0,\infty )$ uniformly in $t$ $\in $ $\mathbb{N}$; $f$ $:$ $\mathbb{R}%
^{k_{x}}$ $\times $ $\mathcal{B}$ $\rightarrow $ $\mathbb{R}$, where $%
\mathcal{B}$ $\equiv $ $\mathcal{D}$ $\times $ $\Theta $ $\subset $ $\mathbb{%
R}^{k_{\beta }}$ with $\Theta $ $=$ $\times _{i=1}^{k_{\theta }}\Theta _{i}$%
; $\mathcal{D}$ $\subset $ $\mathbb{R}^{k_{\delta }}$ and $\Theta _{i}$ $%
\subset $ $\mathbb{R}$ are compact. $f(x,\cdot )$ is for each $x$ Borel
measurable, and $f(\cdot ,\beta )$ is three times continuously
differentiable on $\mathcal{B}.\medskip $\newline
$b.$ $\emph{Identification}$: $\beta _{0}$ $\equiv $ $[\delta _{0}^{\prime
},\theta _{0}^{\prime }]^{\prime }$ uniquely minimizes $E[(y_{t}-f(x_{t},%
\beta ))^{2}]$ on $\mathcal{B}$ $\forall t$, where $\delta _{0}$ and $\theta
_{0,i}$ are interior points of $\mathcal{D}$ and $\Theta _{i}$. $%
E[(y_{t}-f(x_{t},[\delta ,\boldsymbol{0}_{k_{\theta }}]))^{2}]$ has $\forall
t$ a unique minimum on $\mathcal{D}.\medskip $\newline
$c$. \emph{Sub-Exponential Tails}: Let $w_{i,t}$ denote $\sup_{\beta
_{(i)}\in \mathcal{B}_{(i)}}|g_{(i)}(x_{t},\beta _{(i)})|$, $\sup_{\delta
\in \mathcal{D}}|f(x_{t},\beta _{0})$ $-$ $f(x_{t},[\delta ,\boldsymbol{0}%
_{k_{\theta }}])|$, $\sup_{\beta _{(i)}\in \mathcal{B}_{(i)}}|\{f(x_{t},%
\beta _{0})$ $-$ $f(x_{t},[\delta ,\boldsymbol{0}_{k_{\theta
}}])\}g_{(i)}(x_{t},\beta _{(i)})|$ or $\sup_{\beta _{(i)}\in \mathcal{B}%
_{(i)}}|\{f(x_{t},\beta _{0})-f(x_{t},[\delta ,\boldsymbol{0}_{k_{\theta
}}])\}h_{(i)}(x_{t},\beta _{(i)})|$. Then $\max_{t\in \mathbb{N}}P(|\epsilon
_{t}|$ $>$ $c)$ $\leq $ $C\exp \{-\mathcal{K}c\}$ and $\max_{i,t\in \mathbb{N%
}}P(|w_{i,t}|$ $>$ $c)$ $\leq $ $C\exp \{-\mathcal{K}c\}$ for some finite $C,%
\mathcal{K}$ $>$ $0$ that may be different for different $(\epsilon
_{t},w_{i,t})$.
\end{assumption}

\begin{remark}
\normalfont($c$) imposes uniform sub-exponentiality on key first order
expansion terms. It yields the existence of moment generating functions for
first order components, and therefore promotes exponential moment bounds
Assumptions \ref{assum:expand}.b(ii),c. ($c$) is low level, and seems
necessary ultimately for supporting a non-Gaussian approximation theory in a
parametric regression setting. Additional envelope moment bounds that arise
due to high dimensionality and working with parametric functions are
presented as Assumption \ref{assum:param_boot}.d in 
\citet[Appendix
B.4]{supp_mat_testmanyzeros}.
\end{remark}

Define the least squares gradient process under the null $G_{i,t}^{(0)}$ $%
\equiv $ $(y_{t}$ $-$ $f(x_{t},\beta ^{(0)})g_{(i)}(x_{t},\beta
_{(i)}^{(0)}) $, and the normalized bootstrapped gradient:%
\begin{equation*}
\mathcal{\tilde{Z}}_{(i)}^{(0)}\equiv -\left( \frac{1}{n}\sum_{t=1}^{n}E%
\left[ g_{(i)}(x_{t},\beta _{(i)}^{(0)})g_{(i)}(x_{t},\beta
_{(i)}^{(0)})^{\prime }\right] \right) ^{-1}\frac{1}{\sqrt{n}}%
\sum_{t=1}^{n}\eta _{t}G_{i,t}^{(0)}.
\end{equation*}

We first have an asymptotic bootstrap expansion similar to Lemma \ref%
{lm:expansion}.

\begin{lemma}
\label{lm:b_boot_expansion}Let Assumption \ref{assum:param_boot} hold, and
let the weights satisfy $\max_{1\leq i\leq k_{\theta ,n}}|\mathcal{W}_{n,i}$ 
$-$ $\mathcal{W}_{i}|$ $=$ $o_{p}(1/(\ln (n))^{2})$ for non-stochastic $%
\mathcal{W}_{i}$ $\in $ $(0,\infty )$. Then for any monotonic sequence of
positive integers $\{k_{\theta ,n}\}$, $k_{\theta ,n}$ $=$ $o(\sqrt{n})$:%
\begin{equation*}
\left\vert \max_{1\leq i\leq k_{\theta ,n}}\left\vert \sqrt{n}\mathcal{W}%
_{n,i}\widehat{\tilde{\theta}}_{i}\right\vert -\max_{1\leq i\leq k_{\theta
,n}}\left\vert \mathcal{W}_{i}[\boldsymbol{0}_{k_{\delta }}^{\prime },1]%
\mathcal{\tilde{Z}}_{(i)}^{(0)}\right\vert \right\vert \overset{p}{%
\rightarrow }0.
\end{equation*}
\end{lemma}

\begin{remark}
\normalfont The restriction $k_{\theta ,n}$ $=$ $o(\sqrt{n})$ is slightly
more lenient than $k_{\theta ,n}$ $=$ $o(\sqrt{n}/(\ln (n))^{4})$ for the
Lemma \ref{lm:expansion} expansion, ultimately due to the iid Gaussian draw $%
\eta _{t}$ in the bootstrap case.
\end{remark}

Now define%
\begin{equation*}
\mathcal{\bar{H}}_{n,(i)}^{(0)}\equiv \frac{1}{n}\sum_{t=1}^{n}E\left[
g_{(i)}(x_{t},\beta _{(i)}^{(0)})g_{(i)}(x_{t},\beta _{(i)}^{(0)})^{\prime }%
\right] \text{, hence }\mathcal{\tilde{Z}}_{(i)}^{(0)}=-\frac{1}{\sqrt{n}}%
\sum_{t=1}^{n}\mathcal{\bar{H}}_{n,(i)}^{(0)-1}\eta _{t}G_{i,t}^{(0)},
\end{equation*}%
and%
\begin{eqnarray}
&&\tilde{\sigma}_{n,(i)}^{2}(\lambda )\equiv E\left[ \left( \lambda ^{\prime
}\mathcal{\tilde{Z}}_{(i)}^{(0)}\right) ^{2}\right] =\lambda ^{\prime }%
\mathcal{\bar{H}}_{n,(i)}^{(0)}\frac{1}{n}\sum_{t=1}^{n}E\left[
G_{i,t}^{(0)}G_{i,t}^{(0)\prime }\right] \mathcal{\bar{H}}%
_{n,(i)}^{(0)}\lambda ^{\prime }  \label{Vni0} \\
&&\tilde{\sigma}_{(i)}^{2}(\lambda )=\lim_{n\rightarrow \infty }\tilde{\sigma%
}_{n,(i)}^{2}(\lambda ).  \notag
\end{eqnarray}%
Notice $\tilde{\sigma}_{n,(i)}^{2}(\lambda )$, and $\sigma
_{n,(i)}^{2}(\lambda )$ $=$ $\lambda ^{\prime }\mathfrak{H}_{(i)}^{-1}%
\mathcal{S}_{n}(i,i)\mathfrak{H}_{(i)}^{-1}\lambda $ in (\ref{ZVsig}), are
identical asymptotically under $H_{0}$. Likewise $\tilde{\sigma}%
_{(i)}^{2}(\lambda )$ and $\sigma _{(i)}^{2}(\lambda )$ $=$ $%
\lim_{n\rightarrow \infty }\sigma _{n,(i)}^{2}(\lambda )$ are identical
under $H_{0}$. As usual, we only care about the selection vector $\lambda $ $%
=$ $[\boldsymbol{0}_{k_{\delta }}^{\prime },1]^{\prime }$.

Let $\Rightarrow ^{p}$ denote \textit{weak convergence in probability} (Gin%
\'{e} and Zinn, \citeyear{GineZinn1990}: Section 3; see also van der Vaart
and Wellner, \citeyear{VaartWellner1996}, Chapt. 2.9). This notion of
convergence is convenient for characterizing a multiplier bootstrap p-value
limit theory \citep[e.g.][]{Hansen1996}.

We have the following conditional multiplier bootstrap central limit theorem.

\begin{lemma}
\label{lm_ZZ}Let Assumption \ref{assum:param_boot} hold. Let $\{\boldsymbol{%
\tilde{Z}}_{n,(i)}(\lambda )\}_{n\in \mathbb{N}}$ be sequences of scalar
normal random variables, $\boldsymbol{\tilde{Z}}_{n,(i)}(\lambda )$ $\sim $ $%
N(0,\tilde{\sigma}_{n,(i)}^{2}(\lambda ))$, independent of the sample $%
\mathfrak{S}_{n}$ $\equiv $ $\{x_{t},y_{t}\}_{t=1}^{n}$. Then for any
monotonic sequence of positive integers $\{k_{\theta ,n}\}$ that satisfies $%
k_{\theta ,n}=o(\sqrt{n}/\ln (n)^{2})$: 
\begin{equation}
\sup_{z\geq 0}\left\vert P\left( \max_{1\leq i\leq k_{\theta ,n}}\left\vert
\lambda ^{\prime }\mathcal{\tilde{Z}}_{(i)}^{(0)}\right\vert \leq z|%
\mathfrak{S}_{n}\right) -P\left( \max_{1\leq i\leq k_{\theta ,n}}\left\vert 
\boldsymbol{\tilde{Z}}_{n,(i)}(\lambda )\right\vert \leq z\right)
\right\vert \overset{p}{\rightarrow }0.  \label{PZ|W-PZ}
\end{equation}%
Now let $\{\boldsymbol{\tilde{Z}}_{(i)}(\lambda )\}_{i\in \mathbb{N}}$, $%
\boldsymbol{\tilde{Z}}_{(i)}(\lambda )$ $\sim $ $N(0,\tilde{\sigma}%
_{(i)}^{2}(\lambda ))$, be an independent copy of the Theorem \ref%
{th:max_theta_hat}.a null distribution process $\{\boldsymbol{Z}%
_{(i)}(\lambda )\}_{i\in \mathbb{N}}$, that is independent of the asymptotic
draw $\{x_{t},y_{t}\}_{t=1}^{\infty }$. Then $\max_{1\leq i\leq k_{\theta
,n}}|[\boldsymbol{0}_{k_{\delta }}^{\prime },1]\mathcal{\tilde{Z}}%
_{(i)}^{(0)}|$ $\Rightarrow ^{p}$ $\max_{i\in \mathbb{N}}|\boldsymbol{\tilde{%
Z}}_{(i)}([\boldsymbol{0}_{k_{\delta }}^{\prime },1]\mathfrak{)}|$.
\end{lemma}

\begin{remark}
\normalfont Result (\ref{PZ|W-PZ}) is based on a (conditional) Slepian-type
inequality, following from Lemma 3.1 in \cite{Chernozhukov_etal2013} %
\citep[cf.][Theorem 2]{Chernozhukov_etal2015}. We ultimately exploit $%
k_{\theta ,n}=o(\sqrt{n}/\ln (n)^{2})$ which is just slightly better than $%
k_{\theta ,n}=o(\sqrt{n}/(\ln (n))^{4})$ in Lemma \ref{lm:expansion}, again
due to the iid Gaussian draw $\eta _{t}$.
\end{remark}

Lemmas \ref{lm:b_boot_expansion} and \ref{lm_ZZ} yield the following
fundamental result.

\begin{lemma}
\label{lm_boot_theta}Let Assumption \ref{assum:param_boot} hold, and let $%
\max_{1\leq i\leq k_{\theta ,n}}|\mathcal{W}_{n,i}$ $-$ $\mathcal{W}_{i}|$ $%
= $ $o_{p}(1/(\ln (n))^{2})$ for non-stochastic $\mathcal{W}_{i}$ $\in $ $%
(0,\infty )$. Let $\{\boldsymbol{\tilde{Z}}_{(i)}(\lambda )\}_{i\in \mathbb{N%
}}$, $\boldsymbol{\tilde{Z}}_{(i)}(\lambda )$ $\sim $ $N(0,\tilde{\sigma}%
_{(i)}^{2}(\lambda ))$, be an independent copy of the Theorem \ref%
{th:max_theta_hat}.a null distribution process $\{\boldsymbol{Z}%
_{(i)}(\lambda )\}_{i\in \mathbb{N}}$, that is independent of the asymptotic
draw $\{x_{t},y_{t}\}_{t=1}^{\infty }$. Then $\max_{1\leq i\leq k_{\theta
,n}}|\sqrt{n}\mathcal{W}_{n,i}\widehat{\tilde{\theta}}_{i}|$ $\Rightarrow
^{p}$ $\max_{i\in \mathbb{N}}|\mathcal{W}_{i}\boldsymbol{\tilde{Z}}_{(i)}([%
\boldsymbol{0}_{k_{\delta }}^{\prime },1]\mathfrak{)}|$ for any monotonic
sequence of positive integers $\{k_{\theta ,n}\}$ that satisfies $k_{\theta
,n}$ $=$ $o(\sqrt{n}/\ln (n)^{2})$.
\end{lemma}

We now have the main result of this section: the approximate p-value $\tilde{%
p}_{n,\mathcal{M}_{n}}$ promotes a correctly sized and consistent test
asymptotically.

\begin{theorem}
\label{th:p_value_boot}Let Assumptions \ref{assum:ident}-\ref%
{assum:param_boot} hold, and let $\max_{1\leq i\leq k_{\theta ,n}}|\mathcal{W%
}_{n,i}$ $-$ $\mathcal{W}_{i}|$ $=$ $o_{p}(1/(\ln (n))^{2})$ for
non-stochastic $\mathcal{W}_{i}$ $\in $ $(0,\infty )$.\medskip \newline
$a.$ Under $H_{0}$, $P(\tilde{p}_{n,\mathcal{M}_{n}}$ $<$ $\alpha )$ $%
\rightarrow $ $\alpha $ for any monotonic sequence of positive integers $%
\{k_{\theta ,n}\}_{n\geq 1}$ with $k_{\theta ,n}$ $=$ $o(\sqrt{n}/(\ln
(n))^{4})$.\medskip \newline
$b.$ Under $H_{1}$, $P(\tilde{p}_{n,\mathcal{M}_{n}}$ $<$ $\alpha )$ $%
\rightarrow $ $1$ if $\theta _{0,i}$ $\neq $ $0$ for some $i$ $\in $ $%
\mathbb{N}$ such that $1$ $\leq $ $i$ $\leq $ $\lim_{n\rightarrow \infty
}k_{\theta ,n}$.
\end{theorem}

\subsubsection{Local Power}

The bootstrapped max-test has non-trivial power against a sequence of $\sqrt{%
n}$-local alternatives $H_{1}^{L}$ $:$ $\theta _{0}$ $=$ $c/\sqrt{n}$ where $%
c$ $=$ $\left[ c_{i}\right] _{i=1}^{k_{\theta }}$ $\in $ $\mathbb{R}%
^{k_{\theta }}$. The theory developed in Sections \ref{sec:first_order} and %
\ref{sec:max_test_dist} easily carries over to $H_{1}^{L}$. A formal proof
of the following is omitted because it is identical to arguments under $%
H_{0} $. Let $k_{\theta ,n}=o(\sqrt{n}/(\ln (n))^{4})$. Under $H_{1}^{L}$\
we have a first order expansion $|\max_{1\leq i\leq k_{\theta ,n}}|\sqrt{n}%
\mathcal{W}_{n,i}(\hat{\theta}_{i}$ $-$ $c_{i}/\sqrt{n})|$ $-$ $\max_{1\leq
i\leq k_{\theta ,n}}|\mathcal{W}_{i}[\boldsymbol{0}_{k_{\delta }}^{\prime
},1]\mathcal{\hat{Z}}_{(i)}||$ $\overset{p}{\rightarrow }$ $0$ and therefore
Gaussian approximation:%
\begin{equation*}
\left\vert P\left( \max_{1\leq i\leq k_{\theta ,n}}\left\vert \sqrt{n}%
\mathcal{W}_{n,i}\left( \hat{\theta}_{i}-c_{i}/\sqrt{n}\right) \right\vert
\leq z\right) -P\left( \max_{1\leq i\leq k_{\theta ,n}}\left\vert 
\boldsymbol{Z}_{n,(i)}([\boldsymbol{0}_{k_{\delta }}^{\prime
},1])\right\vert \leq z\right) \right\vert \rightarrow 0\text{ }\forall
z\geq 0,
\end{equation*}

Now define for arbitrary $d$ $\in $ $\mathbb{R}^{k_{\theta }}$ and sample $%
\mathfrak{S}_{n}$ $=$ $\{x_{t},y_{t}\}_{t=1}^{n}$: 
\begin{equation*}
\mathcal{P}_{n}(d)=P\left( \max_{1\leq i\leq k_{\theta ,n}}\left\vert \sqrt{n%
}\mathcal{W}_{n,i}\widehat{\tilde{\theta}}_{1}\right\vert >\max_{1\leq i\leq
k_{\theta ,n}}\left\vert \sqrt{n}\mathcal{W}_{n,i}\left( \hat{\theta}_{i}-%
\frac{c_{i}}{\sqrt{n}}\right) +d_{i}\right\vert |\mathfrak{S}_{n}\right) .
\end{equation*}%
By the Glivenko-Cantelli Theorem with $\mathcal{M}$ $=$ $\mathcal{M}_{n}$ $%
\rightarrow $ $\infty $, the triangle inequality, and \linebreak $%
\max_{1\leq i\leq k_{\theta ,n}}|\mathcal{W}_{n,i}$ $-$ $\mathcal{W}_{i}|$ $%
= $ $o_{p}(1/(\ln (n))^{2})$, the bootstrapped p-value approximation
satisfies: 
\begin{eqnarray}
\tilde{p}_{n,\mathcal{M}} &=&P\left( \max_{1\leq i\leq k_{\theta
,n}}\left\vert \sqrt{n}\mathcal{W}_{n,i}\widehat{\tilde{\theta}}%
_{1}\right\vert >\max_{1\leq i\leq k_{\theta ,n}}\left\vert \sqrt{n}\mathcal{%
W}_{n,i}\left( \hat{\theta}_{i}-\frac{c_{i}}{\sqrt{n}}\right) +\mathcal{W}%
_{i}c_{i}\right\vert |\mathfrak{S}_{n}\right) +o_{p}(1)  \notag \\
&\equiv &\mathcal{P}_{n}\left( \left\{ \mathcal{W}_{i}c_{i}\right\}
_{i=1}^{k_{\theta ,n}}\right) +o_{p}(1).  \label{PnHL}
\end{eqnarray}%
Now, by the proof of Theorem \ref{th:p_value_boot}, cf. 
\citet[Theorem 2 and its
proof]{Hansen1996}, if $c$ $=$ $\boldsymbol{0}_{k_{\theta }}$ then $\mathcal{%
P}_{n}(\{\mathcal{W}_{i}c_{i}\}_{i=1}^{k_{\theta ,n}})$ $=$ $\mathcal{P}_{n}(%
\boldsymbol{0}_{k_{\theta }})$ is asymptotically uniformly distributed on $%
[0,1]$, hence $P(\mathcal{P}_{n}(\boldsymbol{0}_{k_{\theta }})$ $<$ $\alpha
) $ $\rightarrow $ $\alpha $ provided $k_{\theta ,n}$ $=$ $o(\sqrt{n}/(\ln
(n))^{4})$. Under $H_{1}^{L}$ generally, from (\ref{PnHL}) the bootstrap
test p-value satisfies $\lim_{n\rightarrow \infty }P(\tilde{p}_{n,\mathcal{M}%
_{n}}$ $<$ $\alpha )$ $\nearrow $ $1$ monotonically as the maximum local
drift $\max_{1\leq i\leq k_{\theta }}|c_{i}|$ $\nearrow $ $\ \infty $ since $%
\mathcal{W}_{i}$ $\in $ $(0,\infty )$ $\forall i$. Asymptotic local power
therefore does not depend on the degree, if any, of sparsity.

By comparison, for example, \cite{Zhong_Chen_Xu_2013} consider a high
dimensional mean model $y_{i}$ $=$ $\theta _{0}$ $+$ $\epsilon _{i}$ and
test $H_{0}$ $:$ $\theta _{0}$ $=$ $\boldsymbol{0}_{p}$, where $y_{i}$ $\in $
$\mathbb{R}^{p}$. Under their alternative $p^{1-\xi }$ elements $\theta
_{0,i}$ $\neq $ $0$ for some $\xi $ $\in $ $(1/2,1)$, where $\theta _{0,i}$ $%
=$ $\sqrt{r\ln (p)/n}$ for some $r$ $>$ $0$, and by assumption $\ln (p)$ $=$ 
$o(n^{1/3})$, hence $\theta _{0,i}$ $=$ $o(1/n^{2/3})$. See also \cite%
{Arias-Castro_Candes_Plan_2011} and \cite{Delaigle_Hall_Jin_2011} and their
references. Our max-test is consistent against deviations $\theta _{0}$ $=$ $%
c/n^{\zeta }$ for any $\zeta $ $\in $ $[0,1/2)$ and $c$ $\neq $ $\boldsymbol{%
0}_{k_{\theta }}$, and has non-trivial power when $\zeta $ $=$ $1/2$ and $c$ 
$\neq $ $\boldsymbol{0}$. Further, any positive number of elements $\theta
_{0,i}$ $\neq $ $0$ under $H_{1}$ or $H_{1}^{L}$ (i.e. $1$ $\leq $ $%
\sum_{i=0}^{k_{\theta }}I(\theta _{0,i}$ $\neq $ $0)$ $\leq $ $k_{\theta }$%
). Thus, in the simple setting of \cite{Zhong_Chen_Xu_2013}, alternatives
can be closer to the null than allowed here but sparsity is enforced. Here,
a far more complex parametric setting is allowed, and any degree of
(non)sparsity is permitted under $H_{1}$ or $H_{1}^{L}$ (as long as some $%
\theta _{0,i}$ $\neq $ $0$)$.$

\section{Example: Linear Regression\label{sec:example}}

We verify all assumptions for a linear regression model 
\citep[see][Appendix
D.2, for a logistic regression model]{supp_mat_testmanyzeros}: 
\begin{equation*}
y_{t}=\delta _{0}^{\prime }x_{\delta ,t}+\theta _{0}^{\prime }x_{\theta
,t}+\epsilon _{t}=\beta _{0}^{\prime }x_{t}+\epsilon _{t}.
\end{equation*}%
Let $x_{\delta ,t}$ $\in $ $\mathbb{R}^{k_{\delta }}$, $x_{\theta ,t}$ $\in $
$\mathbb{R}^{k_{\theta }}$, and $E[\epsilon _{t}]$ $=$ $0$. Assume $%
\min_{t\in \mathbb{N}}E[\epsilon _{t}^{2}]$ $>$ $0$ and $\min_{i,t\in 
\mathbb{N}}E[x_{i,t}^{2}]$ $>$ $0$. Assume $\{x_{t},\epsilon _{t}\}$ are
independent across $t$ and mutually independent for unique $\beta _{0}$ $=$ $%
[\delta _{0}^{\prime },\theta _{0}^{\prime }]^{\prime }$. $\delta _{0}$ and $%
\theta _{0,i}$ are interior points of compact $\mathcal{D}$ $\subset $ $%
\mathbb{R}^{k_{\delta }}$, $k_{\delta }$ $\in $ $\mathbb{N}$, and $\Theta
_{i}$ $\subset $ $\mathbb{R}$. Assume $\Theta $ $\subset $ $\{\times
_{i=1}^{k_{\theta }}\Theta _{i}$ $:$ $|\theta |$ $<$ $\infty \}$. In a
linear framework we need $|\theta |$ $<$ $\infty $ to ensure $\sup_{\beta
\in \mathcal{B}}||\beta ^{\prime }x_{t}||_{p}$ $<$ $\infty $ for some $p$ $>$
$4$ under general conditions on $x_{t}$, which is used to verify several
conditions in bootstrap Assumption \ref{assum:param_boot}.d.\footnote{$%
|\theta |$ $<$ $\infty $ is trivial when $k_{\theta }$ $<$ $\infty $, and
covers sparse (i.e. $\sum_{i=1}^{k_{\theta }}I(\theta _{i}$ $\neq $ $0)$ $<$ 
$\infty $) and nonsparse (e.g. $\theta _{i}$ $=$ $O(\rho ^{i})$, $|\rho |$ $<
$ $1$) cases. However, we do not require $|\theta |$ $<$ $\infty $ with
logistic or similar "squash"-like response \citep[see Appendix D.2
in][]{supp_mat_testmanyzeros}.}

Let $\epsilon _{t}$ and $x_{i,t}$\ be uniformly sub-exponential: $\max_{t\in 
\mathbb{N}}P(|\epsilon _{t}|$ $>$ $c)$ $\leq $ $\mathcal{C}\exp \{-\mathcal{K%
}c\}$ and \linebreak $\max_{i,t\in \mathbb{N}}P(|x_{i,t}|$ $>$ $c)$ $\leq $ $%
\mathcal{C}\exp \{-\mathcal{K}c\}$ for some finite $\mathcal{C},\mathcal{K}$ 
$>$ $0$ that may be different in different places. Assume all pairs $%
(x_{i,t},x_{j,t}$ $:$ $i$ $\neq $ $j)$ are conditionally sub-exponential: $%
\max_{i,j,t\in \mathbb{N}}P(\left\vert x_{i,t}\right\vert $ $>$ $c|x_{j,t})$ 
$\leq $ $\max_{j,t\in \mathbb{N}}\mathcal{C}(x_{j,t})\exp \left\{ -\mathcal{K%
}c\right\} $ $a.s.$ for some $\sigma (x_{j,t})$-measurable random variable $%
\mathcal{C}(x_{j,t})$ with $E[\max_{j,t\in \mathbb{N}}\mathcal{C}%
(x_{j,t})^{2}]$ $<$ $\infty $. The latter expedites the tail bounds in
bootstrap Assumption \ref{assum:param_boot}.c, and includes covariates with
bounded support, and mutually independent covariates.

We want to test $H_{0}$ $:$ $\theta _{0}$ $=$ $0$. The parsimonious models
are 
\begin{equation*}
y_{t}=\delta _{(i)}^{\ast \prime }x_{\delta ,t}+\theta _{i}^{\ast }x_{\theta
,i,t}+v_{(i),t}=\beta _{(i)}^{\ast \prime }x_{(i),t}+v_{(i),t}\text{, }%
i=1,...,k_{\theta ,n},
\end{equation*}%
where $E[v_{(i),t}x_{(i),t}]$ $=$ $0$ for unique $\beta _{(i)}^{\ast }$ in
the interior of $\mathcal{B}_{(i)}$ $\equiv $ $\mathcal{D}\times \Theta _{i}$%
. Squared error loss is used: $\mathcal{L}(\beta )$ $=$ $.5E[(y_{t}$ $-$ $%
\beta ^{\prime }x_{t})^{2}]$ and $\mathcal{L}_{(i)}(\beta )$ $=$ $.5E[(y_{t}$
$-$ $\beta _{(i)}^{\prime }x_{(i).t})^{2}]$. We estimate $\beta _{(i)}^{\ast
}$ by least squares with criterion $\mathcal{\hat{L}}_{(i)}(\beta )$ $\equiv 
$ $.5\sum_{t=1}^{n}(y_{t}$ $-$ $\beta _{(i)}^{\prime }x_{(i).t})^{2}$, hence 
$\widehat{\mathcal{G}}_{(i)}$ $=$ $-\sum_{t=1}^{n}v_{(i),t}x_{(i),t},$ $%
\widehat{\mathcal{H}}_{(i)}$ $=$ $\sum_{t=1}^{n}x_{(i),t}x_{(i),t}^{\prime }$%
, \ and therefore $\mathcal{H}_{(i)}$ $=$ $\mathcal{\bar{H}}_{(i)}$ $=$ $%
\mathfrak{H}_{(i)}$ $=$ $\lim_{n\rightarrow \infty
}1/n\sum_{t=1}^{n}E[x_{(i),t}x_{(i),t}^{\prime }]$.

Assume $\inf_{\lambda ^{\prime }\lambda =1}\min_{1\leq i\leq k_{\theta
,n}}\{1/n\sum_{t=1}^{n}(\lambda ^{\prime }x_{(i),t})^{2}\}$ $>$ $0$ $a.s$.,
and $\max_{1\leq i\leq k_{\theta ,n}}||(1/n\sum_{t=1}^{n}x_{(i),t}$ $\times $
$x_{(i),t}^{\prime })^{-1}||$ $=$ $O_{p}(k_{\theta ,n})$. The latter holds,
for example, when there are $k_{\delta }$ $=$ $0$ nuisance parameters and $%
\min_{i\in \mathbb{N}}\lim \inf_{n\rightarrow \infty
}1/n\sum_{t=1}^{n}x_{i,t}^{2}$ $>$ $0$ $a.s.$ Finally, we need to assume:%
\begin{eqnarray*}
&&\max_{1\leq i\leq k_{\theta ,n}}\left\vert \frac{1}{n}\sum_{t=1}^{n}E\left[
\epsilon _{t}^{2}\right] E\left[ x_{i,t}^{2}\right] -\lim_{n\rightarrow
\infty }\frac{1}{n}\sum_{t=1}^{n}E\left[ \epsilon _{t}^{2}\right] E\left[
x_{i,t}^{2}\right] \right\vert =O(k_{\theta ,n}/\sqrt{n}) \\
&&\max_{1\leq i\leq k_{\theta ,n}}\left\vert \frac{1}{n}\sum_{t=1}^{n}E\left[
\left( \beta _{0}^{\prime }x_{t}\right) ^{2}x_{i,t}^{2}\right]
-\lim_{n\rightarrow \infty }\frac{1}{n}\sum_{t=1}^{n}E\left[ \left( \beta
_{0}^{\prime }x_{t}\right) ^{2}x_{i,t}^{2}\right] \right\vert =O(k_{\theta
,n}/\sqrt{n}).
\end{eqnarray*}%
Each restricts heterogeneity, and is trivial when $(\epsilon _{t},x_{i,t})$
are identically distributed across $t$, or have bounded heterogeneity.

See \citet[Appendix D.1]{supp_mat_testmanyzeros} for a proof of the
following result.

\begin{lemma}
\label{lm:ex_linear_reg}Assumptions \ref{assum:ident}-\ref{assum:param_boot}%
\ hold in the setting above.
\end{lemma}

\begin{remark}
\normalfont The result extends with only mild changes to the proof to a
general class of additively nonlinear models, for example, $y_{t}$ $=$ $%
\delta _{0,2}^{\prime }x_{\delta ,t}F(\delta _{0,1}^{\prime }x_{\delta ,t})$ 
$+$ $\theta _{0}^{\prime }x_{\theta ,t}$ $+$ $\epsilon _{t}$. where $F(\cdot
)$ is known, measurable, bounded ("squash"-like) and suitably smooth. We
need both $\delta _{0,i}$ $\neq $ $0$ to ensure identification. Examples
include neural network models and smooth switching and threshold models
(see, e.g., Lee and Porter, \citeyear{LeePorter1984}, and White, %
\citeyear{White1989}, for historical references).
\end{remark}

\section{Monte Carlo Experiments\label{sec:sim}}

We perform a Monte Carlo study based on a linear regression model for
cross-sectional data. We use 10,000 independently drawn samples of sizes $n$ 
$\in $ $\{100,250,500,1000\}$. In each case two max-tests are performed: one
uses a flat weight $\mathcal{W}_{n,i}$ $=$ $1$, and the other $\mathcal{W}%
_{n,i}$ is equal to the inverted standard error of $\hat{\theta}_{n,i}$. We
call these the \textit{max-test} and \textit{max-t-test} respectively. Other
tests are discussed below. The number of parsimonious regression models is $%
k_{\theta ,n}$ $\in $ $\{10,35,[5n^{1/2-\iota }]\}$ where $\iota $ $=$ $%
10^{-10}$. Note $[5n^{1/2-\iota }]$ $\in $ $\{50,79,112,158\}$ for $n$ $\in $
$\{100,250,500,1000\}$. There are $k_{\delta }$ $\in $ $\{0,10\}$ nuisance
parameters. All bootstrapped p-values are based on $\mathcal{M}$ $=$ $1000$
independently drawn samples from $N(0,1).$The significance levels are $%
\alpha $ $\in $ $\{.01,.05,.10\}$.

The following presents our benchmark design, after which we discuss
robustness checks.

\subsection{Test of Linear Regression Zero Restrictions}

\subsubsection{Set-Up}

The DGP is $y_{t}$ $=$ $\delta _{0}^{\prime }x_{\delta ,t}$ $+$ $\theta
_{0}^{\prime }x_{\theta ,t}$ $+$ $\epsilon _{t}$ where $\epsilon _{t}$ is
iid standard normal. We consider three cases for the regressors. In the
first case $[x_{\delta ,t},x_{\theta ,t}]$ are serially and mutually
independent standard normals. In the second case they are block-wise
dependent normals ($x_{\delta ,t}$ and $x_{\theta ,t}$ are mutually
independent). In the third case they are within and across block dependent
normals. Not surprisingly all tests perform better under the first two
cases. We therefore only report results for the more realistic third case.

Covariates in case three are drawn as follows. Define $k$ $\equiv $ $%
k_{\delta }$ $+$ $k_{\theta ,n}$. Combine $x_{t}$ $\equiv $ $[x_{\delta
,t}^{\prime },x_{\theta ,t}^{\prime }]^{\prime }\in $ $\mathbb{R}^{k}$, and
let $(w_{t},v_{t})$\ be mutually independent draws from $N(0,I_{k})$. The
regressors are $x_{t}$ $=$ $Aw_{t}$ $+$ $v_{t}$. We randomly draw each
element of $A\in $ $\mathbb{R}^{k}$ from a uniform distribution on $[-1,1]$.
If $A$ does not have full column rank then we add a randomly drawn $\iota $
from $[0,1]$ to each diagonal component (in every case the resulting $A$ had
full column rank).\medskip

We fix $\delta $ $=$ $\mathbf{1}_{k_{\delta }}$, a $k_{\delta }$ $\times $ $%
1 $ vector of ones. The benchmark models are as follows. Under the null $%
\theta _{0}$ $=$ $[\theta _{0,i}]_{i=1}^{k_{\theta ,n}}$ $=$ $0$. The
alternatives are $(i)$ $\theta _{0,1}$ $=$ $.001$ with $[\theta
_{0,i}]_{i=2}^{k_{\theta ,n}}$ $=$ $0$; $(ii)$ $\theta _{0,i}$ $=$ $%
i/k_{\theta ,n}$ for $i$ $=$ $1,...,k_{\theta ,n}$; and $(iii)$ $\theta
_{0,i}$ $=$ $.001$ for $i$ $=$ $1,...,k_{\theta ,n}$. See Table \ref{tbl_h1}
for reference. Under $(i)$ we have a small deviation from the null. At large
deviations (e.g. $\theta _{0,1}$ $=$ $1/2$) all tests are comparable. Small
deviations like $\theta _{0,1}$ $\in $ $(0,.1]$, however, reveal test
differences. That said, as we will see below the max-test generally
dominates by size, and by power for small deviations.

\begin{table}[h!]
\caption{Alternative Models}
\label{tbl_h1}%
\begin{tabular}{ll}
Alternative $(i):$ & $\theta _{0,1}$ $=$ $.001$ and $[\theta _{0,i}]_{i=2}^{%
\mathring{k}_{\theta ,n}}$ $=$ $0$ \\ 
Alternative $(ii):$ & $\theta _{0,i}$ $=$ $i/\mathring{k}_{\theta ,n}$ for $%
i $ $=$ $1,...,\mathring{k}_{\theta ,n}$ \\ 
Alternative $(iii):$ & $\theta _{0,i}$ $=$ $.001$ for $i$ $=$ $1,...,%
\mathring{k}_{\theta ,n}$%
\end{tabular}%
\end{table}

In case $(iii)$ all $\theta _{0,i}$ have the same value. A max-test might
plausibly be sub-optimal to a bootstrapped Wald test since the latter
combines each estimated parameter, while the max-statistic works with the
single largest estimate of a sequence of identical parameters. We will see
below that contrary to this logic, the max-test works very well, precisely
because it incorporates so little, yet relevant, information. A Wald
statistic for large $k_{\theta ,n}$ adds significant sampling error due to
the inverted variance matrix and due to low degrees of freedom by estimating
the original high dimensional model.

Benchmark results are discussed in Section \ref{sec:sim_bench}. In Section %
\ref{sec:sim_rob} we investigate additional models as robustness checks
against the benchmark cases. We discuss the various tests next.

\subsubsection{Max-Tests}

We estimate $k_{\theta ,n}$ parsimonious models $y_{t}$ $=$ $\delta
_{(i)}^{\ast \prime }x_{\delta ,t}$ $+$ $\theta _{i}^{\ast }x_{\theta ,i,t}$ 
$+$ $v_{(i),t}$ by least squares. Denote by $\mathcal{T}_{n}$ the resulting
max-test or max-t-test statistic.\medskip \textbf{\ }The bootstrapped test
statistic $\mathcal{\tilde{T}}_{n}$ and p-value $\tilde{p}_{n,\mathcal{M}}$
are computed as in (\ref{T_boot}) and (\ref{p_boot}). We reject $H_{0}$ when 
$\tilde{p}_{n,\mathcal{M}}$ $<$ $\alpha $.

\subsubsection{Wald Tests}

We also perform asymptotic Wald and parametric bootstrapped Wald tests, and
asymptotic and bootstrapped normalized Wald tests. In all cases the model $%
y_{t}$ $=$ $\delta _{0}^{\prime }x_{\delta ,t}$ $+$ $\theta _{0}^{\prime
}x_{\theta ,t}$ $+$ $\epsilon _{t}$ is estimated by least squares, where $%
x_{\theta ,t}$ has dimension $k_{\theta ,n}$. The asymptotic Wald test is
based on the $\chi ^{2}(k_{\theta ,n})$ distribution, which is only valid
(asymptotically) for fixed $k_{\theta ,n}$ $=$ $k$\ or $k_{\theta ,n}$ $%
\rightarrow $ $(0,\infty )$.

The bootstrap Wald test is based on the fixed-design parametric wild
bootstrap in \cite{Goncalves_Kilian_2004}. The algorithm is as follows. Let $%
\mathcal{W}_{n}$ be the Wald statistic for a test of $H_{0}$ $:$ $\theta
_{0} $ $=$ $0$. Write $y_{t}$ $=$ $\delta _{0}^{\prime }x_{\delta ,t}$ $+$ $%
\theta _{0}^{\prime }x_{\theta ,t}$ $+$ $\epsilon _{t}$ $=$ $\beta
_{0}^{\prime }x_{t}$ $+$ $\epsilon _{t}$ and $\epsilon _{t}(\beta )$ $\equiv 
$ $y_{t}$ $-$ $\beta ^{\prime }x_{t}$, and let $\hat{\beta}_{n}$ be the
unrestricted least squares estimator of $\beta _{0}$. Generate a sample $\{%
\hat{y}_{t}^{\ast }\}_{t=1}^{n}$ with the null imposed: $\hat{y}_{t}^{\ast }$
$=$ $\hat{\delta}_{n}^{\prime }x_{\delta ,t}$ $+$ $\epsilon _{t}(\hat{\beta}%
_{n})\eta _{t}$ where $\eta _{t}$ is iid $N(0,1)$. Construct the model $\hat{%
y}_{t}^{\ast }$ $=$ $b^{\prime }x_{t}$ $+$ $\upsilon _{t}^{\ast }$, estimate 
$b$ by least squares and compute a Wald statistic $\mathcal{W}_{n}^{\ast }$
for the null hypothesis. Repeat $\mathcal{M}$ times resulting in the
sequence of bootstrapped Wald statistics $\{\mathcal{W}_{n,i}^{\ast
}\}_{i=1}^{\mathcal{M}}$. The p-value approximation is $p_{n}$ $=$ $%
1/M\sum_{i=1}^{\mathcal{M}}I(\mathcal{W}_{n}$ $>$ $\mathcal{W}_{n,i}^{\ast
}) $, and we reject $H_{0}$ when $p_{n}$ $<$ $\alpha $.

The normalized Wald statistic is $\mathcal{W}_{n}^{s}$ $\equiv $ $(\mathcal{W%
}_{n}$ $-$ $k_{\theta ,n})/\sqrt{2k_{\theta ,n}}$. Under the null $\mathcal{W%
}_{n}^{s}$ $\overset{d}{\rightarrow }$ $N(0,1)$ as $k_{\theta ,n}$ $%
\rightarrow $ $\infty $, and as long as $k_{\theta ,n}/n$ $\rightarrow $ $0$
then $\mathcal{W}_{n}^{s}$ $\overset{p}{\rightarrow }$ $\infty $ under $%
H_{1} $. The asymptotic test is therefore a one-sided test which rejects the
null when $\mathcal{W}_{n}^{s}$ $>$ $Z_{\alpha }$, where $Z_{\alpha }$ is
the standard normal upper tail $\alpha $-level critical value. The
bootstrapped test is performed as above: we compute $\mathcal{W}%
_{n,i}^{s\ast }$ $\equiv $ $(\mathcal{W}_{n,i}^{\ast }$ $-$ $k_{\theta ,n})/%
\sqrt{2k_{\theta ,n}}$ with p-value approximation $p_{n}$ $=$ $%
1/M\sum_{i=1}^{\mathcal{M}}I(\mathcal{W}_{n}^{s}$ $>$ $\mathcal{W}%
_{n,i}^{s\ast })$. Trivially $\mathcal{W}_{n}^{s}$ $>$ $\mathcal{W}%
_{n,i}^{s\ast }$ \textit{if and only if} $\mathcal{W}_{n,i}^{\ast }$ $>$ $%
\mathcal{W}_{n}$, hence the bootstrapped normalized Wald test and
bootstrapped Wald test are identical. We therefore only discuss the
bootstrapped Wald test.

\subsection{Simulation Results}

\subsubsection{Benchmark Results\label{sec:sim_bench}}

\paragraph{Empirical Size}

Test results under the three covariate cases are similar, although all tests
perform slightly less well under correlated regressors within and across
blocks (case three) which we focus on.

The asymptotic Wald test and asymptotic normalized Wald test are generally
severely over-sized due to the magnitude of $k_{\delta }$ $+$ $k_{\theta ,n}$%
. The bootstrapped Wald test is strongly under-sized when $n$ is small, but
also when $k_{\delta }$ $=$ $10$ and/or when $k_{\theta ,n}$ is large. The
max-tests typically lead to qualitatively similar results with empirical
size close to nominal size (see Table \ref{table:max_wald_rejH0}).

\paragraph{Empirical Power}

Under alternative $(i)$\ only one parameter deviates from the null, $\theta
_{0,1}$ $=$ $.001$. The difficulty in detecting this deviation is clearly
apparent at smaller sample sizes. Across sample sizes and size of $k_{\delta
}$, moreover, the max-tests dominate the bootstrapped Wald test. The
max-tests can achieve a power gain of as much as $40\%$ over the Wald tests.
For example, at $n$ $=$ $100$, $k_{\delta }$ $=$ $10$ and $k_{\theta ,n}$ $=$
$50$ the Wald test has zero power above nominal size while the max-test and
max-t-test have power $\{.039,.149,.266\}$ and $\{.099,.267,.427\}$ for
sizes $\{1\%,5\%,10\%\}$. At $n$ $=$ $250$, $k_{\delta }$ $=$ $10$ and $%
k_{\theta ,n}$ $=$ $79$, max-test and max-t-test power are $%
\{.514,.742,.838\}$ and $\{.617,.799,.873\}$, yet Wald test power is only $%
\{.008,.096,.260\}$. A similar result occurs at $n$ $=$ $500.$

When there are zero nuisance covariates $k_{\delta }$ $=$ $0$, the
differences are just as stark. Max-test and max-t-test power when $n$ $=$ $%
100$ and $k_{\theta ,n}$ $=$ $35$, for example, are $\{.451,.762,.870\}$ and 
$\{.531,.795,.890\}$, while Wald test power is $\{.020,.133,.251\}$. Only
when $n$ $=$ $1000$ does power across tests become comparable. Thus, the
max-tests dominate in terms of detecting a singular small deviation from the
null.

Under alternative $(ii)$ we have $\theta _{0,i}$ $=$ $i/k_{\theta ,n}$ for $%
i $ $=$ $1,...,k_{\theta ,n}$: all parameters deviate from zero at values
between $(0,1]$. In this case the max- and max-t-tests, and the Wald test,
are generally similar. This is not surprising since the (bootstrapped) Wald
test becomes competitive under $(i)$ when $\theta _{0,1}$ has values farther
from zero.

Under $(iii)$, where $\theta _{0,i}$ $=$ $.001$ for $i$ $=$ $1,...,k_{\theta
,n}$, all tests are naturally better able to detect a deviation from the
null than under $(i)$. The max-tests, however, dominate the remaining tests
at $n$ $\leq $ $250$, in particular when there are nuisance parameters $%
k_{\delta }$ $>$ $0$. See Tables \ref{table:max_wald_rejH11}-\ref%
{table:max_wald_rejH13}. As with alternative $(i)$, under $(iii)$ all tests
are comparable when $\theta _{0,i}$ is larger, e.g. $\theta _{0,i}$ $=$ $1/2$
for each $i$ (not shown).

\subsubsection{Robustness Checks\label{sec:sim_rob}}

We consider various robustness checks against the above design:
heterogeneous covariates, and large(r) $k_{\delta }$.

\paragraph{a. Regressor Dispersion}

While the covariates $x_{\theta ,i,t}$ are not homogeneous in all benchmark
cases, their variances are fairly similar or identical depending on case.
The dispersion of $x_{\theta ,i,t}$, however, has a direct impact on the
dispersion of the key parsimonious model estimators $\hat{\theta}_{i}$.

We therefore add three covariate cases allowing for different variances
across $x_{\theta ,t}$ $=$ $[x_{\theta ,i,t}]$ $\sim $ $N(0,\Psi )$. We
bypass covariate dependence and let $[x_{\delta ,t},x_{\theta ,t}]$ be
within and across block iid, in order to focus on the pure effects of
covariate dispersion. When present, each nuisance $x_{\delta ,i,t}$ is iid $%
N(0,1)$. In each case we use a diagonal variance matrix $\Psi $ $=$ $[\Psi
_{i,j}]$. In case $(a)$ $\Psi _{i,i}$ $=$ $1$ $+$ $100(i$ $-$ $1)/k_{\theta
,n}$ $\in $ $[1,100]$ hence $x_{\theta ,i,t}$ has a monotonically larger
variance as $i$ increases. In this setting a max-t-test may perform better
due to self-scaling. In the remaining two cases diagonal $\Psi $ has entries 
$(b)$ $\Psi _{1,1}$ $=$ $10$ or $(c)$ $\Psi _{1,1}$ $=$ $100,$ and all other 
$\Psi _{i,i}$ $=$ $1$, hence there is a single potentially influential
regressor.

We only inspect test performance under the null. See Tables 5-7 in %
\citet[Appendix E]{supp_mat_testmanyzeros}. Covariate dispersion matters
more when there are large differences. The first case with increasing
variance from $1$ to $100$ does not lead to major size distortions
generally. The second and third cases, however, where there is one
influential covariate, the Wald test can suffer major size distortions,
while the max-test and max-t-test are similar, or in some cases the
max-t-test has slightly better empirical size as anticipated. The max-tests
are mildly oversized. The slight size distortions still exist when $n$ $=$ $%
1000$, although they are not as attenuated as at smaller sample sizes. The
Wald test suffers major size distortions even at $n$ $=$ $1000$ when $%
k_{\theta ,n}$ is large (in particular when $k_{\delta }$ $>$ $0$).

\paragraph{b. Large $k_{\protect\delta }$}

Finally, we study two cases $k_{\delta }$ $\in $ $\{20,40\}$ under the null
where $k_{\delta }$ is larger than in the baseline cases $\{0,10\}$. The
covariates are within and across group iid standard normal. In view of the
large number of parameters estimated, computation time is heavy, hence we
only report results for $n$ $\in $ $\{100,250,500\}$. Results are contained
in Table 8 in \citet[Appendix E]{supp_mat_testmanyzeros}. The max-tests are
slightly over-sized, in particular when $k_{\delta }$ $=$ $40$, while the
bootstrapped Wald test has greater size distortions. At lower sample sizes
the max-test dominates, but the max-test and max-t-test are comparable at $n$
$=$ $500$.

\section{Conclusion\label{sec:conclude}}

We present a class of max-tests of zero restrictions for a high dimensional
parameter. In a general loss and extremum estimator setting we prove that
the hypotheses can be identified by using many potentially vastly lower
dimension parameterizations. We use a test statistic that is the maximum in
absolute value of the weighted key estimated parameters across the many low
dimension settings. The lower dimension helps improve estimation accuracy,
while a max-statistic alleviates the need for a multivariate normalization
used in Wald and score statistics. Thus, we avoid inverting a potentially
large dimension variance estimator that may be a poor proxy for the true
sampling dispersion. Asymptotic theory sidesteps traditional extreme value
theoretic arguments. We instead focus on an approximate p-value computed by
parametric wild bootstrap for a nonlinear regression model and independent
data. In simulation experiments the max-tests generally have good or sharp
size, and dominate a bootstrapped (and normalized) Wald test in terms of
size and power. Future work should ($i$) lean toward allowing for a broader
class of semi-nonparametric models, including models with a high dimensional
parameter $\theta _{0}$ \textit{and} infinite dimensional unknown function $%
h $; and ($ii$) verify key exponential moment bounds for first order
summands in a general time series setting.

\setcounter{equation}{0} \renewcommand{\theequation}{{\thesection}.%
\arabic{equation}} \appendix

\section{Appendix: Proofs\label{app:proofs}}

\setstretch{1.2} Let $\boldsymbol{0}_{d}$ be a $d$ $\times $ $1$ vector of
zeros, and write $(\partial /\partial x)f(x_{0})$ $=$ $[(\partial /\partial
x_{i})f(x)|_{x=x_{0}}]$. Throughout $O_{p}(1)$ and $o_{p}(1)$ are not
functions of model counter $i$. $\{k_{\theta ,n},k_{n}\}$ are a
monotonically increasing sequences of positive integers. Supporting lemmas
are proved in \citet[Appendix
C]{supp_mat_testmanyzeros}.\medskip

The proof of asymptotic approximation Lemma \ref{lm:expansion} requires a
general stochastic domination theory for the maximum of possibly
non-Gaussian random variables. This is presented as three general lemmas,
plus a high dimensional rate of convergence bound.

\begin{lemma}
\label{lm:max_p-converg}Let $\{\mathcal{X}_{n,i}\}_{n\geq 1}$ be random
variables on a common probability space, $i$ $=$ $1,2,...$ Then $\max_{1\leq
i\leq k_{n}}|\mathcal{X}_{n,i}|$ $=$ $O_{p}(\ln \left( k_{n}\right) \ln (n))$
for any sequence of positive integers $\{k_{n}\}$ provided $\max_{1\leq
i\leq k_{n}}E[\exp \left\{ \vartheta \left\vert \mathcal{X}_{n,i}\right\vert
\right\} ]$ $=$ $O(n^{\xi \ln (k_{n})})$ for some $\vartheta ,\xi $ $>$ $0$.
\end{lemma}

The following related lemma will play a role in the bootstrap asymptotic
theory.

\begin{lemma}
\label{lm:max_p-converg:eta}Let $\{k_{n}\}_{n\geq 1}$ be any sequence of
positive integers with $k_{n}$ $=$ $O(n)$. Let $\{w_{i,t}\}_{t=1}^{\infty }$
be random variables on a common probability space, independent across $t$,
for $i$ $=$ $1,2,...$.\medskip \newline
$a.$ Assume either ($i$) $\max_{i\in \mathbb{N},1\leq t\leq n}|w_{i,t}|$ $%
\leq $ $\mathcal{M}_{n}$ $a.s.$ for non-random positive $\mathcal{M}_{n}$ $=$
$O(\ln \left( n\right) )$; or ($ii$) $\max_{i,t\in \mathbb{N}}P(|w_{i,t}|$ $%
> $ $c)$ $\leq $ $C\exp \{-\mathcal{K}c\}$ for some finite $C,\mathcal{K}$ $%
> $ $0$. Let $\{\eta _{t}\}_{t=1}^{n}$ be iid $N(0,1)$ random variables,
independent of $\{w_{i,t}\}_{t=1}^{n}$, $i$ $=$ $1,2..$. Then $\max_{1\leq
i\leq k_{n}}|1/\sqrt{n}\sum_{t=1}^{n}\eta _{t}w_{i,t}|$ $=$ $O_{p}(\ln
(k_{n})\ln (n))$.$\medskip $\newline
$b.$ If $\max_{i,t\in \mathbb{N}}E[w_{i,t}^{2}]$ $<$ $\infty $ then $%
\max_{1\leq i\leq k_{n}}|1/n\sum_{t=1}^{n}(\eta _{t}^{2}w_{i,t}$ $-$ $%
E[w_{i,t}])|$ $=$ $O_{p}(k_{n}/\sqrt{n})$.
\end{lemma}

In bootstrap cases where the Lemma \ref{lm:max_p-converg} condition $%
\max_{1\leq i\leq k_{n}}E[\exp \left\{ \vartheta \left\vert \mathcal{X}%
_{n,i}\right\vert \right\} ]$ $=$ $O(n^{\xi \ln (k_{n})})$ need not hold
(and is not strictly required), we have the following more restrictive
result for partial sums.

\begin{lemma}
\label{lm:max_p-converg2}Let $\{x_{i,t}\}_{t=1}^{\infty }$ be random
variables on a common probability space, $i$ $=$ $1,2,...$ Assume $x_{i,t}$
are independent across $t$, $E[x_{i,t}]$ $=$ $0$ $\forall i,t$, and $%
1/n\sum_{t=1}^{n}\max_{i\in \mathbb{N}}\{E[x_{i,t}^{2}]\}$ $=$ $O(1)$. Then $%
\max_{1\leq i\leq k_{n}}|1/\sqrt{n}\sum_{t=1}^{n}x_{i,t}|$ $=$ $O_{p}(k_{n})$
for any sequence of positive integers $\{k_{n}\}$.
\end{lemma}

We also need the (uniform) rate of convergence of $\hat{\beta}_{(i)}$.

\begin{lemma}
\label{lm:beta_hat_rate}Under Assumptions \ref{assum:consist_suff} and \ref%
{assum:expand} $\max_{1\leq i\leq k_{\theta ,n}}|\sqrt{n}(\hat{\beta}_{(i)}$ 
$-$ $\beta _{(i)}^{\ast })|$ $=$ $O_{p}(\left( \ln \left( n\right) \right)
^{2})$ provided $k_{\theta ,n}$ $=$ $O(\sqrt{n}/(\ln (n))^{2})$.
\end{lemma}

\noindent \textbf{Proof of Lemma \ref{lm:expansion}.} Recall $\{\varpi
_{n}\}_{n\in \mathbb{N}}$ are non-stochastic positive numbers, $\varpi _{n}$ 
$\rightarrow $ $0$, and may be different in different places, and 
\begin{equation*}
\mathcal{B}_{n,(i)}=\mathcal{B}_{n,(i)}(\varpi _{n})\equiv \left\{ \beta
_{(i)}\in \mathcal{B}_{(i)}:\left\Vert \beta _{(i)}-\beta _{(i)}^{\ast
}\right\Vert \leq \varpi _{n}\right\} .
\end{equation*}

Recall $\mathcal{\hat{Z}}_{(i)}$ $\equiv $ $-\mathfrak{H}_{(i)}^{-1}\widehat{%
\mathcal{G}}_{(i)}/\sqrt{n}$, $\widehat{\mathfrak{H}}_{(i)}(\cdot )$ $\equiv 
$ $\widehat{\mathcal{H}}_{(i)}(\cdot )/n$, $\mathfrak{H}_{(i)}(\cdot )$ $%
\equiv $ $\plim_{n\rightarrow \infty }\widehat{\mathcal{H}}_{(i)}(\cdot )/n$
and $\mathfrak{H}_{(i)}$ $=$ $\mathfrak{H}_{(i)}(\beta _{(i)}^{\ast })$.
Assume $k_{\theta ,n}$ $=$ $o(\sqrt{n}/(\ln (n))^{4})$. By multiple uses of
the triangle inequality:%
\begin{eqnarray*}
&&\left\vert \max_{1\leq i\leq k_{\theta ,n}}\left\vert \sqrt{n}\mathcal{W}%
_{n,i}\hat{\theta}_{i}\right\vert -\max_{1\leq i\leq k_{\theta
,n}}\left\vert \mathcal{W}_{i}[\boldsymbol{0}_{k_{\delta }}^{\prime },1]%
\mathcal{\hat{Z}}_{(i)}\right\vert \right\vert \\
&&\text{ \ \ \ \ \ \ \ \ \ \ \ \ \ \ \ }\leq \max_{1\leq i\leq k_{\theta
,n}}\left\vert \sqrt{n}\mathcal{W}_{n,i}\hat{\theta}_{i}-\mathcal{W}_{i}[%
\boldsymbol{0}_{k_{\delta }}^{\prime },1]\mathcal{\hat{Z}}_{(i)}\right\vert
\\
&&\text{ \ \ \ \ \ \ \ \ \ \ \ \ \ \ \ }\leq \max_{i\in \mathbb{N}%
}\left\vert \mathcal{W}_{i}\right\vert \times \max_{1\leq i\leq k_{\theta
,n}}\left\vert \sqrt{n}\hat{\theta}_{i}-[\boldsymbol{0}_{k_{\delta
}}^{\prime },1]\mathcal{\hat{Z}}_{(i)}\right\vert \\
&&\text{ \ \ \ \ \ \ \ \ \ \ \ \ \ \ \ \ \ \ \ \ \ \ \ \ \ \ }+\max_{1\leq
i\leq k_{\theta ,n}}\left\vert \mathcal{W}_{n,i}-\mathcal{W}_{i}\right\vert
\times \max_{1\leq i\leq k_{\theta ,n}}\left\vert \sqrt{n}\hat{\theta}_{i}-[%
\boldsymbol{0}_{k_{\delta }}^{\prime },1]\mathcal{\hat{Z}}_{(i)}\right\vert
\\
&&\text{ \ \ \ \ \ \ \ \ \ \ \ \ \ \ \ \ \ \ \ \ \ \ \ \ \ \ }+\max_{1\leq
i\leq k_{\theta ,n}}\left\vert [\boldsymbol{0}_{k_{\delta }}^{\prime },1]%
\mathcal{\hat{Z}}_{(i)}\right\vert \times \max_{1\leq i\leq k_{\theta
,n}}\left\vert \mathcal{W}_{n,i}-\mathcal{W}_{i}\right\vert .
\end{eqnarray*}%
By assumption $\max_{1\leq i\leq k_{\theta ,n}}|\mathcal{W}_{n,i}$ $-$ $%
\mathcal{W}_{i}|$ $=$ $o_{p}(1/(\ln (n))^{2})$ and $\max_{i\in \mathbb{N}}|%
\mathcal{W}_{i}|$ $<$ $\infty $. Furthermore $\max_{i\in \mathbb{N}}||%
\mathfrak{H}_{(i)}^{-1}||$ $<$ $\infty $ under Assumption \ref{assum:expand}%
.b(ii). Assumption \ref{assum:expand}.c and Lemma \ref{lm:max_p-converg}
yield $\max_{1\leq i\leq k_{\theta ,n}}|\widehat{\mathcal{G}}_{(i)}/\sqrt{n}%
| $ $=$ $O_{p}(\ln \left( k_{\theta ,n}\right) \ln (n))$. Hence since $%
k_{\theta ,n}$ $=$ $O(n)$, $\max_{1\leq i\leq k_{\theta ,n}}|[\boldsymbol{0}%
_{k_{\delta }}^{\prime },1]\mathcal{\hat{Z}}_{(i)}|$ $\times $ $\max_{1\leq
i\leq k_{\theta ,n}}|\mathcal{W}_{n,i}$ $-$ $\mathcal{W}_{i}|$ $\overset{p}{%
\rightarrow }$ $0$. It therefore suffices to prove 
\begin{equation}
\max_{1\leq i\leq k_{\theta ,n}}\left\vert \sqrt{n}\hat{\theta}_{i}-[%
\boldsymbol{0}_{k_{\delta }}^{\prime },1]\mathcal{\hat{Z}}_{(i)}\right\vert 
\overset{p}{\rightarrow }0.  \label{maxtheta_Z}
\end{equation}

By (\ref{1st_expand}) and uniform positive definiteness of $\widehat{%
\mathfrak{H}}_{(i)}(\cdot )$ under Assumption \ref{assum:expand}.b(i), $%
\sqrt{n}(\hat{\beta}_{(i)}$ $-$ $\beta _{(i)}^{\ast })$ $=$ $-\widehat{%
\mathfrak{H}}_{(i)}^{-1}(\ddot{\beta}_{(i)})\widehat{\mathcal{G}}_{(i)}/%
\sqrt{n}$. The mean value theorem yields for some $\{\ddot{\beta}_{(i)}\}$, $%
||\ddot{\beta}_{(i)}$ $-$ $\beta _{(i)}^{\ast }||$ $\leq $ $||\hat{\beta}%
_{(i)}$ $-$ $\beta _{(i)}^{\ast }||$, that may be different in different
places:%
\begin{eqnarray*}
&&\max_{1\leq i\leq k_{\theta ,n}}\left\vert \sqrt{n}\lambda ^{\prime
}\left( \hat{\beta}_{(i)}-\beta _{(i)}^{\ast }\right) -\lambda ^{\prime }%
\mathcal{\hat{Z}}_{(i)}\right\vert \\
&&\text{ \ \ \ \ \ \ \ \ \ \ \ }\leq \left\{ \max_{1\leq i\leq k_{\theta
,n}}\left\vert \lambda ^{\prime }\left\{ \widehat{\mathfrak{H}}_{(i)}^{-1}(%
\ddot{\beta}_{(i)})-\mathfrak{H}_{(i)}^{-1}(\ddot{\beta}_{(i)})\right\}
\right\vert \right. \\
&&\text{ \ \ \ \ \ \ \ \ \ \ \ \ \ \ \ \ \ \ \ \ \ \ }\left. +\max_{1\leq
i\leq k_{\theta ,n}}\lambda ^{\prime }\left\{ \mathfrak{H}_{(i)}^{-1}(\ddot{%
\beta}_{(i)})-\mathfrak{H}_{(i)}^{-1}(\beta _{(i)}^{\ast })\right\} \right\}
\times \max_{1\leq i\leq k_{\theta ,n}}\left\vert \frac{1}{\sqrt{n}}\widehat{%
\mathcal{G}}_{(i)}\right\vert .
\end{eqnarray*}%
Local Lipschitz and Hessian bounds Assumption \ref{assum:expand}.a,b yield
for finite $\mathcal{K}$ $>$ $0$:%
\begin{eqnarray*}
\max_{1\leq i\leq k_{\theta ,n}}\left\vert \mathfrak{H}_{(i)}^{-1}(\ddot{%
\beta}_{(i)})-\mathfrak{H}_{(i)}^{-1}(\beta _{(i)}^{\ast })\right\vert &\leq
&\max_{1\leq i\leq k_{\theta ,n}}\left\vert \mathfrak{H}_{(i)}^{-1}(\ddot{%
\beta}_{(i)})\right\vert \times \left\vert \mathfrak{H}_{(i)}(\ddot{\beta}%
_{(i)})-\mathfrak{H}_{(i)}^{-1}(\beta _{(i)}^{\ast })\right\vert \times
\left\vert \mathfrak{H}_{(i)}(\beta _{(i)}^{\ast })\right\vert \\
&\leq &\mathcal{K}\times \left\vert \ddot{\beta}_{(i)}-\beta _{(i)}^{\ast
}\right\vert .
\end{eqnarray*}%
Now use Lemma \ref{lm:beta_hat_rate}, $||\ddot{\beta}_{(i)}$ $-$ $\beta
_{(i)}^{\ast }||$ $\leq $ $||\hat{\beta}_{(i)}$ $-$ $\beta _{(i)}^{\ast }||$%
, and $k_{\theta ,n}$ $=$ $o(\sqrt{n}/(\ln (n))^{4})$, to yield: 
\begin{equation*}
\max_{1\leq i\leq k_{\theta ,n}}\left\vert \mathfrak{H}_{(i)}^{-1}(\ddot{%
\beta}_{(i)})-\mathfrak{H}_{(i)}^{-1}(\beta _{(i)}^{\ast })\right\vert
=O_{p}(\ln \left( k_{\theta ,n}\right) /\sqrt{n}).
\end{equation*}%
Finally, $\max_{1\leq i\leq k_{\theta ,n}}|\widehat{\mathcal{G}}_{(i)}/\sqrt{%
n}|$ $=$ $O_{p}(\ln \left( k_{\theta ,n}\right) \ln (n))$ from above, and
from the proof of Lemma \ref{lm:beta_hat_rate} we have 
\citep[see][eq.
(C.12)]{supp_mat_testmanyzeros}:%
\begin{equation*}
\max_{1\leq i\leq k_{\theta ,n}}\left\{ \sup_{\beta _{(i)}\in \mathcal{B}%
_{n,(i)}}\left\vert \widehat{\mathfrak{H}}_{(i)}^{-1}(\beta _{(i)})-%
\mathfrak{H}_{(i)}^{-1}(\beta _{(i)})\right\vert \right\} =O_{p}\left( \frac{%
k_{\theta ,n}\ln \left( k_{\theta ,n}\right) \ln (n)}{\sqrt{n}}\right) .
\end{equation*}%
Hence, for any $\{k_{\theta ,n}\}$, $k_{\theta ,n}$ $=$ $o(\sqrt{n}/(\ln
(n))^{4})$:%
\begin{eqnarray*}
&&\max_{1\leq i\leq k_{\theta ,n}}\left\vert \sqrt{n}\lambda ^{\prime
}\left( \hat{\beta}_{(i)}-\beta _{(i)}^{\ast }\right) -\lambda ^{\prime }%
\mathcal{\hat{Z}}_{(i)}\right\vert \\
&&\text{ \ \ \ \ \ \ \ \ \ \ \ \ \ \ \ }=\left\{ O_{p}\left( \frac{k_{\theta
,n}\ln \left( k_{\theta ,n}\right) \ln (n)}{\sqrt{n}}\right) +O_{p}\left( 
\frac{\ln \left( k_{\theta ,n}\right) }{\sqrt{n}}\right) \right\} \times
O_{p}\left( \ln \left( k_{\theta ,n}\right) \ln (n)\right) \\
&&\text{ \ \ \ \ \ \ \ \ \ \ \ \ \ \ \ }=O_{p}\left( \frac{k_{\theta
,n}\left( \ln \left( k_{\theta ,n}\right) \right) ^{2}\left( \ln (n)\right)
^{2}}{\sqrt{n}}\right) =o_{p}(1)
\end{eqnarray*}%
Now invoke $H_{0}$ $:$ $\theta _{0}$ $=$ $0$ and put $\lambda $ $=$ $[%
\boldsymbol{0}_{k_{\delta }}^{\prime },1]$\ to complete the proof. $\mathcal{%
QED}.\bigskip $\newline
\textbf{Proof of Lemma \ref{lm:max_dist}.} In view of Assumption \ref%
{assum:max_dist} and $B_{n}^{2}(\ln (k_{\theta ,n}n))^{7}$ $=$ $O(n^{1-\xi
}) $ the claim follows from Corollary 2.1 of \cite{Chernozhukov_etal2013}. $%
\mathcal{QED}$.\bigskip \newline
\textbf{Proof of Theorem \ref{th:max_theta_hat}.} Convergence in cdfs (\ref%
{PP0}) under ($a$) follows from the imposed properties on $\mathcal{W}_{n,i}$%
, Lemma \ref{lm:expansion}, and Corollary \ref{cor:max_dist} with $\lambda $ 
$=$ $[\boldsymbol{0}_{k_{\delta }}^{\prime },1]^{\prime }$. Then $%
\max_{1\leq i\leq k_{\theta ,n}}|\sqrt{n}\mathcal{W}_{n,i}\hat{\theta}_{i}|$ 
$\overset{d}{\rightarrow }$ $\max_{i\in \mathbb{N}}|\mathcal{W}_{i}%
\boldsymbol{Z}_{(i)}([\boldsymbol{0}_{k_{\delta }}^{\prime },1])|$ follows
immediately given $\boldsymbol{Z}_{(i)}(\lambda )$ $\sim $ $%
N(0,\lim_{n\rightarrow \infty }\sigma _{n,(i)}^{2}(\lambda ))$, and Gaussian
distributions are fully characterized by their means and variances. ($b$)
follows from Theorem \ref{th:consist} consistency of $\hat{\theta}_{i}$, and 
$\mathcal{W}_{n,i}$ $\overset{p}{\rightarrow }$ $\mathcal{W}_{i}$ with $%
\mathcal{W}_{i}$ $\in $ $(0,\infty )$. $\mathcal{QED}$.\medskip

In order to prove the bootstrap expansion Lemma \ref{lm:b_boot_expansion} we
utilize several supporting results. Recall $\hat{\beta}^{(0)}$ $\equiv $ $[%
\hat{\delta}^{(0)\prime },\boldsymbol{0}_{k_{\theta }}^{\prime }]^{\prime }$
where $\sum_{t=1}^{n}\{y_{t}$ $-$ $f(x_{t},[\hat{\delta}^{(0)},\boldsymbol{0}%
_{k_{\theta }}])\}^{2}$ $\leq $ $\sum_{t=1}^{n}\{y_{t}$ $-$ $f(x_{t},[\delta
,\boldsymbol{0}_{k_{\theta }}])\}^{2}$ $\forall \delta $ $\in $ $\{\mathcal{D%
}$ $:$ $||\hat{\delta}^{(0)}$ $-$ $\delta ||$ $>$ $0\}$, and write $\hat{%
\beta}_{(i)}^{(0)}$ $\equiv $ $[\hat{\delta}^{(0)\prime },0]^{\prime }$. Let 
$\delta ^{(0)}$ minimize $E[(y_{t}$ $-$ $f(x_{t}[\delta ,\boldsymbol{0}%
_{k_{\theta }}]))^{2}]$ on $\mathcal{D}$, which is unique for each $t$ by
Assumption \ref{assum:param_boot}.b; and write $\beta ^{(0)}$ $=$ $[\delta
^{(0)\prime },\boldsymbol{0}_{k_{\theta }}^{\prime }]^{\prime }$. Recall $%
y_{n,t}^{\ast }$ $\equiv $ $f(x_{t},\hat{\beta}^{(0)})$ $+$ $\epsilon
_{n,t}^{(0)}\eta _{t}$ and $\epsilon _{n,t}^{(0)}$ $\equiv $ $y_{t}$ $-$ $%
f(x_{t},\hat{\beta}^{(0)})$, and note by construction $f_{(i)}(x_{t},\hat{%
\beta}_{(i)}^{(0)})$ $=$ $f(x_{t},\hat{\beta}^{(0)})$.

We first require a standard uniform law of large numbers

\begin{lemma}
\label{lm:ulln_gen}Let $w_{t}$ $:$ $\mathcal{A}$ $\rightarrow $ $\mathbb{R}$%
, $t$ $\in $ $\mathbb{N}$, be well defined random functions on a common
probability measure space, with compact $\mathcal{A}$ $\subset $ $\mathbb{R}%
^{k}$, $k$ $\in $ $\mathbb{N}$. Assume $w_{t}(\alpha )$ is almost surely
continuously differentiable on $\mathcal{A}$ for all $t$, independent, and $%
\max_{t\in \mathbb{N}}E[\sup_{\alpha \in \mathcal{A}}|(\partial /\partial
\alpha )^{j}w_{t}(\alpha )|^{2}]$ $<$ $\infty $ for $j$ $=$ $0,1$. Then $%
\sup_{\alpha \in \mathcal{A}}|1/n\sum_{t=1}^{n}\{w_{t}(\alpha )$ $-$ $%
E[w_{t}(\alpha )]\}|$ $\overset{p}{\rightarrow }$ $0$.
\end{lemma}

\begin{remark}
\label{rm:ulln}\normalfont An essentially identical argument yields $%
\sup_{\alpha ,\tilde{\alpha}\in \mathcal{A}}|1/n\sum_{t=1}^{n}%
\{w_{t}(a)w_{t}(\tilde{\alpha})$ $-$ $E[w_{t}(a)w_{t}(\tilde{\alpha}]\}|$ $%
\overset{p}{\rightarrow }$ $0$ when $\max_{t\in \mathbb{N}}E[\sup_{\alpha
\in \mathcal{A}}|(\partial /\partial \alpha )^{j}w_{t}(\alpha )|^{2}]$ $<$ $%
\infty $ for $j$ $=$ $0,1$.
\end{remark}

Next, the (uniform) rate of convergence for $\hat{\beta}^{(0)}$ and the
bootstrapped $\widehat{\tilde{\beta}}_{(i)}$ $-$ $\hat{\beta}_{(i)}^{(0)}$.

\begin{lemma}
\label{lm:bootbb}Under Assumption \ref{assum:param_boot}: $(a)$ $\hat{\beta}%
^{(0)}$ $\overset{p}{\rightarrow }$ $\beta ^{(0)}$; $(b)$ $\hat{\beta}^{(0)}$
$=$ $\beta ^{(0)}$ $+$ $O_{p}(1/\sqrt{n})$; and $(c)$ $\max_{1\leq i\leq
k_{\theta ,n}}||\widehat{\tilde{\beta}}_{(i)}$ $-$ $\hat{\beta}%
_{(i)}^{(0)}|| $ $=$ $O_{p}((\ln (n))^{2}/\sqrt{n})$ for any monotonically
increasing sequence of positive integers $\{k_{\theta ,n}\}$, $k_{\theta ,n}$
$=$ $o(\sqrt{n})$.
\end{lemma}

\noindent \textbf{Proof of Lemma \ref{lm:b_boot_expansion}}. Define $%
g_{(i)}(x_{t},\beta _{(i)})$ $\equiv $ $(\partial /\partial \beta
_{(i)})f(x_{t},\beta _{(i)})$, $g_{(i),t}^{(0)}$ $\equiv $ $%
g_{(i)}(x_{t},\beta _{(i)}^{(0)})$, and $h_{(i)}(x_{t},\beta _{(i)})$ $%
\equiv $ $(\partial /\partial \beta _{(i)})^{2}f(x_{t},\beta _{(i)})$. Recall%
\begin{equation*}
\mathcal{\tilde{Z}}_{(i)}^{(0)}\equiv -\mathcal{\bar{H}}_{n.(i)}^{(0)-1}%
\frac{1}{\sqrt{n}}\sum_{t=1}^{n}\eta _{t}(y_{t}-f(x_{t},\beta
^{(0)})g_{(i),t}^{(0)}
\end{equation*}%
where $\mathcal{\bar{H}}_{n.(i)}^{(0)}$ $\equiv $ $1/n%
\sum_{t=1}^{n}E[g_{(i),t}^{(0)}g_{(i),t}^{(0)\prime }]$. Define $\mathcal{%
\bar{H}}_{(i)}^{(0)}$ $\equiv $ $\lim_{n\rightarrow \infty }\mathcal{\bar{H}}%
_{n.(i)}^{(0)}$.

We assume $\max_{i\in \mathbb{N}}|\mathcal{W}_{i}|$ $<$ $\infty $ and $%
\max_{1\leq i\leq k_{\theta ,n}}|\mathcal{W}_{n,i}$ $-$ $\mathcal{W}_{i}|$ $%
= $ $o_{p}(1/(\ln (n))^{2})$; and $\max_{1\leq i\leq k_{\theta ,n}}||%
\mathcal{\bar{H}}_{n.(i)}^{(0)-1}||$ $=$ $O(1)$ holds under Assumption \ref%
{assum:param_boot}.d. Under Assumption \ref{assum:param_boot}.c and the
proof of Lemma B.1 in \cite{supp_mat_testmanyzeros}, $(y_{t}$ $-$ $%
f(x_{t},\beta ^{(0)})g_{(i),t}^{(0)}$ satisfies the conditions of Lemma \ref%
{lm:max_p-converg:eta}.a. Coupled with $k_{\theta ,n}$ $=$ $o(\sqrt{n})$, it
therefore follows $\max_{1\leq i\leq k_{\theta ,n}}||1/\sqrt{n}%
\sum_{t=1}^{n}\eta _{t}(y_{t}$ $-$ $f(x_{t},\beta ^{(0)})g_{(i),t}^{(0)}||$ $%
=$ $O_{p}((\ln (n))^{2})$. Hence $\max_{1\leq i\leq k_{\theta ,n}}|[%
\boldsymbol{0}_{k_{\delta }}^{\prime },1]\mathcal{\tilde{Z}}_{(i)}^{(0)}|$ $%
\times $ $\max_{1\leq i\leq k_{\theta ,n}}|\mathcal{W}_{n,i}$ $-$ $\mathcal{W%
}_{i}|$ $\overset{p}{\rightarrow }$ $0$. Now use the argument leading to (%
\ref{maxtheta_Z}) to deduce it suffices to prove 
\begin{equation}
\max_{1\leq i\leq k_{\theta ,n}}\left\vert \sqrt{n}\widehat{\tilde{\theta}}%
_{i}-[\boldsymbol{0}_{k_{\delta }}^{\prime },1]\mathcal{\tilde{Z}}%
_{(i)}^{(0)}\right\vert \overset{p}{\rightarrow }0.  \label{maxWTWZ}
\end{equation}

In the following $O_{p}(\cdot )$ and $o_{p}(\cdot )$ are not functions of $i$%
, and $\{\beta _{(i)}^{\ast (0)}\}$ satisfies $||\beta _{(i)}^{\ast (0)}$ $-$
$\beta _{(i)}^{(0)}||$ $\leq $ $||\hat{\beta}_{(i)}^{(0)}$ $-$ $\beta
_{(i)}^{(0)}||$ and may differ line to line. Recall $f_{(i)}(x_{t},\beta
_{(i)})$ $=$ $f(x_{t},\delta ,[0,...,\theta _{i},0,,,]^{\prime })$. Define%
\begin{eqnarray*}
&&\widehat{\mathfrak{G}}_{(i)}\equiv \frac{1}{n}\sum_{t=1}^{n}g_{(i)}(x_{t},%
\hat{\beta}_{(i)}^{(0)})g_{(i)}(x_{t},\beta _{(i)}^{\ast (0)}) \\
&&\widehat{\mathfrak{Y}}_{(i)}\equiv \frac{1}{n}\sum_{t=1}^{n}\eta
_{t}\left( y_{t}-f_{(i)}(x_{t},\hat{\beta}_{(i)}^{(0)})\right)
h_{(i)}(x_{t},\beta _{(i)}^{\ast (0)}) \\
&&\widehat{\mathfrak{W}}_{(i)}\equiv \frac{1}{n}\sum_{t=1}^{n}\eta
_{t}\left( y_{t}-f(x_{t},\beta _{(i)}^{(0)})\right) g_{(i)}(x_{t},\hat{\beta}%
_{(i)}^{(0)}).
\end{eqnarray*}%
Arguments in 
\citet[proof of Lemma A.6, see
(C.17)-(C.19)]{supp_mat_testmanyzeros} yield:%
\begin{equation}
\widehat{\mathfrak{W}}_{(i)}=\left\{ \widehat{\mathfrak{G}}_{(i)}-\widehat{%
\mathfrak{Y}}_{(i)}\right\} \times \left( \widehat{\tilde{\beta}}_{(i)}-\hat{%
\beta}_{(i)}^{(0)}\right) +O_{p}\left( 1/\sqrt{n}\right) .  \label{FOC}
\end{equation}

Consider $\widehat{\mathfrak{G}}_{(i)}$. Expanding terms with $\hat{\beta}%
_{(i)}^{(0)}$ around $\beta _{(i)}^{(0)}$ yields:%
\begin{eqnarray*}
&&\max_{1\leq i\leq k_{\theta ,n}}\left\Vert \widehat{\mathfrak{G}}_{(i)}-%
\mathcal{\bar{H}}_{n.(i)}^{(0)}\right\Vert \\
&&\text{ \ \ \ }\leq \max_{1\leq i\leq k_{\theta ,n}}\left\Vert \frac{1}{n}%
\sum_{t=1}^{n}g_{(i),t}^{(0)}g_{(i),t}^{(0)\prime }-\mathcal{\bar{H}}%
_{n.(i)}^{(0)}\right\Vert \\
&&\text{ \ \ \ \ \ \ \ }+\max_{1\leq i\leq k_{\theta ,n}}\left\Vert \frac{1}{%
n}\sum_{t=1}^{n}\left( g_{(i)}(x_{t},\hat{\beta}_{(i)}^{(0)})g_{(i)}(x_{t},%
\beta _{(i)}^{\ast (0)})^{\prime }-g_{(i),t}^{(0)}g_{(i),t}^{(0)\prime
}\right) \right\Vert =\mathfrak{A}_{n}+\mathfrak{B}_{n}.
\end{eqnarray*}%
By Lemma \ref{lm:max_p-converg2} and Assumption \ref{assum:param_boot}.d, $%
\mathfrak{A}_{n}$ $\overset{p}{\rightarrow }$ $0$.

If we prove $\mathfrak{B}_{n}\overset{p}{\rightarrow }$ $0$ then it follows%
\begin{equation}
\max_{1\leq i\leq k_{\theta ,n}}\left\Vert \widehat{\mathfrak{G}}_{(i)}-%
\mathcal{\bar{H}}_{n.(i)}^{(0)}\right\Vert \rightarrow 0.  \label{gg}
\end{equation}%
For any pair $(l,m)$ $\in $ $\{1,..,k_{\delta }$ $+$ $1\}$ write%
\begin{equation*}
a_{i,t}\equiv \sup_{\beta _{(i)}\in \mathcal{B}_{(i)}}\left\Vert
g_{m,(i)}(x_{t},\beta _{(i)})\right\Vert \sup_{\beta _{(i)}\in \mathcal{B}%
_{(i)}}\left\Vert \frac{\partial }{\partial \beta _{(i)}}g_{l,(i)}(x_{t},%
\beta _{(i)})\right\Vert .
\end{equation*}%
The mean-value theorem, and adding and subtracting $1/n%
\sum_{t=1}^{n}E[a_{i,t}]$, yield: 
\begin{eqnarray*}
&&\max_{1\leq i\leq k_{\theta ,n}}\left\vert \frac{1}{n}%
\sum_{t=1}^{n}g_{l,(i)}(x_{t},\hat{\beta}_{(i)}^{(0)})g_{m,(i)}(x_{t},\beta
_{(i)}^{\ast (0)})-\frac{1}{n}%
\sum_{t=1}^{n}g_{l,(i),t}^{(0)}g_{m,(i),t}^{(0)}\right\vert \\
&&\text{ \ \ \ \ \ \ \ \ \ \ \ \ \ }\leq 2\max_{1\leq i\leq k_{\theta
,n}}\left\{ \left\vert \frac{1}{n}\sum_{t=1}^{n}\left( a_{i,t}-E\left[
a_{i,t}\right] \right) \right\vert \right\} \times \max_{1\leq i\leq
k_{\theta ,n}}\left\Vert \hat{\beta}_{(i)}^{(0)}-\beta
_{(i)}^{(0)}\right\Vert \\
&&\text{ \ \ \ \ \ \ \ \ \ \ \ \ \ \ \ \ \ \ \ \ \ }+2\max_{1\leq i\leq
k_{\theta ,n}}\left\{ \left\vert \frac{1}{n}\sum_{t=1}^{n}E\left[ a_{i,t}%
\right] \right\vert \right\} \times \max_{1\leq i\leq k_{\theta
,n}}\left\Vert \hat{\beta}_{(i)}^{(0)}-\beta _{(i)}^{(0)}\right\Vert .
\end{eqnarray*}%
Under Assumption \ref{assum:param_boot}.d and the inequality $cd$ $\leq $ $%
c^{2}/2+d^{2}/2$, $a_{i,t}$ $-$ $E[a_{i,t}]$ satisfies the conditions of
Lemma \ref{lm:max_p-converg2}, and $\max_{1\leq i\leq k_{\theta
,n}}|1/n\sum_{t=1}^{n}E[a_{i,t}]|$ $\leq $ $\max_{i,t\in \mathbb{N}%
}|E[a_{i,t}]|$ $<$ $\infty $. Now invoke Lemma \ref{lm:bootbb}, the facts
that $(\hat{\beta}_{(i)}^{(0)},\beta _{(i)}^{(0)})$ are constant over $i$,
and $(l,m)$ are arbitrary, and $k_{\theta ,n}$ $=$ $o(\sqrt{n})$ by
assumption, to conclude 
\begin{equation*}
\mathfrak{B}_{n}=O_{p}(k_{\theta ,n}/\sqrt{n})\times O_{p}(1/\sqrt{n}%
)+O(1)\times O_{p}(1/\sqrt{n})=O_{p}(1/\sqrt{n}).
\end{equation*}

Next $\widehat{\mathfrak{Y}}_{(i)}$. Under Assumption \ref{assum:param_boot}%
.c and the proof of Lemma B.1 in \cite{supp_mat_testmanyzeros}, $\sup_{\beta
_{(i)}\in \mathcal{B}_{(i)}}|(y_{t}$ $-$ $f(x_{t},[\delta ,\boldsymbol{0}%
_{k_{\theta }}]))|$ $\times $ $\sup_{\delta \in \mathcal{D}}|h_{(i)}(x_{t},%
\tilde{\beta}_{(i)})$ satisfies the requirements of Lemma \ref%
{lm:max_p-converg:eta}.a. Then $k_{\theta ,n}$ $=$ $o(\sqrt{n})$ by
supposition yields:%
\begin{equation}
\max_{1\leq i\leq k_{\theta ,n}}\left\vert \widehat{\mathfrak{Y}}%
_{(i)}\right\vert =O_{p}\left( \frac{\ln \left( k_{\theta ,n}\right) \ln
\left( n\right) }{\sqrt{n}}\right) =o_{p}\left( \frac{\left( \ln \left(
n\right) \right) ^{2}}{\sqrt{n}}\right) .  \label{yfh}
\end{equation}

Now $\widehat{\mathfrak{W}}_{(i)}$. The mean value theorem, (\ref{yfh}) and
Lemma \ref{lm:bootbb}.b combined yield: 
\begin{eqnarray}
&&\max_{1\leq i\leq k_{\theta ,n}}\left\vert \frac{1}{n}\sum_{t=1}^{n}\eta
_{t}\left( y_{t}-f(x_{t},\beta ^{(0)})\right) \left\{ g_{(i)}(x_{t},\hat{%
\beta}_{(i)}^{(0)})-g_{(i),t}^{(0)}\right\} \right\vert  \label{etayfgg} \\
&&\text{ \ \ \ \ \ \ \ }\leq \max_{1\leq i\leq k_{\theta ,n}}\left\vert 
\frac{1}{n}\sum_{t=1}^{n}\eta _{t}\left( y_{t}-f_{(i)}(x_{t},\hat{\beta}%
_{(i)}^{(0)})\right) h_{(i)}(x_{t},\beta _{(i)}^{\ast (0)})\right\vert
\max_{1\leq i\leq k_{\theta ,n}}\left\vert \hat{\beta}_{(i)}^{(0)}-\beta
_{(i)}^{(0)}\right\vert  \notag \\
&&\text{ \ \ \ \ \ \ \ }=O_{p}\left( \frac{\ln \left( k_{\theta ,n}\right)
\ln \left( n\right) }{\sqrt{n}}\right) \times O_{p}\left( \frac{1}{\sqrt{n}}%
\right) =o_{p}\left( \frac{\left( \ln \left( n\right) \right) ^{2}}{n}\right)
\notag
\end{eqnarray}

Combine (\ref{FOC})-(\ref{etayfgg}) to deduce:%
\begin{eqnarray*}
&&\frac{1}{n}\sum_{t=1}^{n}\eta _{t}\left( y_{t}-f(x_{t},\beta
^{(0)})\right) g_{(i),t}^{(0)} \\
&&\text{ \ \ \ \ \ \ \ \ \ \ \ }=\left\{ \mathcal{\bar{H}}%
_{n.(i)}^{(0)}+o_{p}\left( \frac{\left( \ln \left( n\right) \right) ^{2}}{%
\sqrt{n}}\right) \right\} \times \left( \widehat{\tilde{\beta}}_{(i)}-\hat{%
\beta}_{(i)}^{(0)}\right) +O_{p}\left( 1/\sqrt{n}\right) .
\end{eqnarray*}%
Under Assumption \ref{assum:param_boot}.d $\mathcal{\bar{H}}_{n,(i)}^{(0)}$
is uniformly positive definite, and $\max_{i\in \mathbb{N}}||\mathcal{\bar{H}%
}_{n,(i)}^{(0)-1}||$ $=$ $O(1)$. We may therefore write: 
\begin{equation*}
\widehat{\tilde{\beta}}_{(i)}-\hat{\beta}_{(i)}^{(0)}=\left( \mathcal{\bar{H}%
}_{n,(i)}^{(0)-1}+o_{p}\left( \frac{\left( \ln \left( n\right) \right) ^{2}}{%
\sqrt{n}}\right) \right) \frac{1}{n}\sum_{t=1}^{n}\eta _{t}\left(
y_{t}-f(x_{t},\beta ^{(0)})\right) g_{(i),t}^{(0)}+O_{p}\left( 1/\sqrt{n}%
\right) .
\end{equation*}%
Arguments in the proof of Lemma \ref{lm:bootbb}.c yield $\max_{1\leq i\leq
k_{\theta ,n}}||1/n\sum_{t=1}^{n}\eta _{t}(y_{t}$ $-$ $f(x_{t},\beta
^{(0)}))g_{(i),t}^{(0)}||$ $=$ $O_{p}((\ln (n))^{2}/\sqrt{n})$ %
\citep[see][(C.17) and subsequent derivations]{supp_mat_testmanyzeros}.
Therefore%
\begin{equation*}
\max_{1\leq i\leq k_{\theta ,n}}\left\Vert \sqrt{n}\left( \widehat{\tilde{%
\beta}}_{(i)}-\hat{\beta}_{(i)}^{(0)}\right) -\mathcal{\bar{H}}%
_{n,(i)}^{(0)-1}\frac{1}{\sqrt{n}}\sum_{t=1}^{n}\eta _{t}\left(
y_{t}-f(x_{t},\beta ^{(0)})\right) g_{(i),t}^{(0)}\right\Vert +o_{p}\left( 
\frac{\left( \ln \left( n\right) \right) ^{4}}{\sqrt{n}}\right) ,
\end{equation*}%
This proves (\ref{maxWTWZ}) since $[\boldsymbol{0}_{k_{\delta }}^{\prime
},1](\widehat{\tilde{\beta}}_{(i)}$ $-$ $\hat{\beta}_{(i)}^{(0)})$ $=$ $%
\widehat{\tilde{\theta}}_{i}$. $\mathcal{QED}$.\bigskip\ \newline
\textbf{Proof of Lemma \ref{lm_ZZ}.\medskip }\newline
\textbf{Step 1.}\qquad Consider (\ref{PZ|W-PZ}). Recall $\mathcal{\tilde{Z}}%
_{(i)}^{(0)}$ $=$ $-\mathcal{\bar{H}}_{n,(i)}^{(0)-1}n^{-1/2}\sum_{t=1}^{n}%
\eta _{t}G_{i,t}^{(0)}$ where $G_{i,t}^{(0)}$ $\equiv $ $(y_{t}$ $-$ $%
f(x_{t},\beta ^{(0)})g_{(i),t}^{(0)}$ with $g_{(i),t}^{(0)}$ $\equiv $ $%
g_{(i)}(x_{t},\beta _{(i)}^{(0)})$, $\mathcal{\bar{H}}_{n,(i)}^{(0)}$ $%
\equiv $ $1/n\sum_{t=1}^{n}E[g_{(i),t}^{(0)}g_{(i),t}^{(0)\prime }]$, and $%
\eta _{t}$ is iid $N(0,1)$ and independent of the sample $\mathfrak{S}_{n}$ $%
\equiv $ $\{x_{t},y_{t}\}_{t=1}^{n}$. As always we only care about $\lambda $
$=$ $[\boldsymbol{0}_{k_{\delta }}^{\prime },1]^{\prime }$, hence uniformity
over $\lambda $ is unnecessary, and we may assume $||\lambda ||$ $=$ $1$.
Write 
\begin{eqnarray*}
&&\hat{\sigma}_{n,(ij)}^{2}(\lambda )\equiv \lambda ^{\prime }E\left[ 
\mathcal{\tilde{Z}}_{(i)}^{(0)}\mathcal{\tilde{Z}}_{(j)}^{(0)}|\mathfrak{S}%
_{n}\right] \lambda =\lambda ^{\prime }\mathcal{\bar{H}}_{n,(i)}^{(0)-1}%
\frac{1}{n}\sum_{t=1}^{n}\eta _{t}^{2}G_{i,t}^{(0)}G_{j,t}^{(0)\prime }%
\mathcal{\bar{H}}_{n,(j)}^{(0)-1}\lambda \\
&&\tilde{\sigma}_{n,(ij)}^{2}(\lambda )\equiv \lambda ^{\prime }\mathcal{%
\bar{H}}_{n,(i)}^{(0)-1}\frac{1}{n}\sum_{t=1}^{n}E\left[
G_{i,t}^{(0)}G_{j,t}^{(0)\prime }\right] \mathcal{\bar{H}}%
_{n,(j)}^{(0)-1}\lambda .
\end{eqnarray*}%
By construction $\tilde{\sigma}_{n,(ii)}^{2}(\lambda )$ $=$ $\tilde{\sigma}%
_{n,(i)}^{2}(\lambda )$ $\rightarrow $ $\tilde{\sigma}_{(i)}^{2}(\lambda )$,
cf. (\ref{Vni0}).

By construction $\lambda ^{\prime }\mathcal{\tilde{Z}}_{(i)}^{(0)}|\mathfrak{%
S}_{n}$ $\sim $ $N(0,\hat{\sigma}_{n,(ii)}^{2}(\lambda ))$. Recall $\{%
\boldsymbol{\tilde{Z}}_{n,(i)}(\lambda )\}_{n\in \mathbb{N}}$ are sequences
with $\boldsymbol{\tilde{Z}}_{n,(i)}(\lambda )$ $\sim $ $N(0,\tilde{\sigma}%
_{n,(ii)}^{2}(\lambda ))$ independent of $\mathfrak{S}_{n}$. Define 
\begin{equation*}
\Delta _{k_{\theta ,n}}(\lambda )\equiv \max_{1\leq i,j\leq k_{\theta
,n}}\left\vert \hat{\sigma}_{n,(ij)}^{2}(\lambda )-\tilde{\sigma}%
_{n,(ij)}^{2}(\lambda )\right\vert .
\end{equation*}%
Lemma 3.1 in \cite{Chernozhukov_etal2013}\ therefore yields the following
conditional Slepian-type inequality 
\citep[cf.][Theorem 2, Proposition
1]{Chernozhukov_etal2015}:%
\begin{eqnarray*}
\mathcal{E}_{k_{\theta ,n}} &\equiv &\sup_{z\geq 0}\left\vert P\left(
\max_{1\leq i\leq k_{\theta ,n}}\left\vert [\boldsymbol{0}_{k_{\delta
}}^{\prime },1]\mathcal{\tilde{Z}}_{(i)}^{(0)}\right\vert \leq z|\mathfrak{S}%
_{n}\right) -P\left( \max_{1\leq i\leq k_{\theta ,n}}\left\vert \boldsymbol{%
\tilde{Z}}_{n,(i)}([\boldsymbol{0}_{k_{\delta }}^{\prime },1]^{\prime
})\right\vert \leq z\right) \right\vert \\
&=&O_{p}\left( \Delta _{k_{\theta ,n}}^{1/3}(\lambda )\max \left\{ 1,\ln
\left( k_{\theta ,n}/\Delta _{k_{\theta ,n}}(\lambda )\right) \right\}
^{2/3}\right) .
\end{eqnarray*}

The proof of (\ref{PZ|W-PZ}) is complete if we show $\mathcal{E}_{k_{\theta
,n}}$ $=$ $o_{p}(1)$. We need to bound $\Delta _{k_{\theta ,n}}(\lambda )$.
By adding and subtracting terms, and the triangle inequality:%
\begin{equation*}
\Delta _{k_{\theta ,n}}(\lambda )\leq \max_{i\in \mathbb{N}}\left\Vert 
\mathcal{\bar{H}}_{n,(i)}^{(0)-1}\right\Vert ^{2}\times \max_{1\leq i,j\leq
k_{\theta ,n}}\left\vert \frac{1}{n}\sum_{t=1}^{n}\left( \eta
_{t}^{2}G_{i,t}^{(0)}G_{j,t}^{(0)\prime }-E\left[ G_{i,t}^{(0)}G_{j,t}^{(0)%
\prime }\right] \right) \right\vert .
\end{equation*}%
Under Assumption \ref{assum:param_boot}.d $\max_{i\in \mathbb{N}}||\mathcal{%
\bar{H}}_{n,(i)}^{(0)-1}||$ $<$ $\infty $, and $E[\eta
_{t}^{2}G_{i,t}^{(0)}G_{j,t}^{(0)\prime }]$ $=$ $E[G_{i,t}^{(0)}G_{j,t}^{(0)%
\prime }]$. Furthermore, Assumption \ref{assum:param_boot}.c and the proof
of Lemma B.1 in \cite{supp_mat_testmanyzeros} imply Lemma \ref%
{lm:max_p-converg:eta}.b applies to $\eta
_{t}^{2}G_{i,t}^{(0)}G_{j,t}^{(0)\prime }$: 
\begin{equation*}
\max_{1\leq i,j\leq k_{\theta ,n}}\left\vert \frac{1}{n}\sum_{t=1}^{n}\left(
\eta _{t}^{2}G_{i,t}^{(0)}G_{j,t}^{(0)\prime }-E\left[
G_{i,t}^{(0)}G_{j,t}^{(0)\prime }\right] \right) \right\vert
=O_{p}(k_{\theta ,n}/\sqrt{n}).
\end{equation*}%
Hence $\Delta _{k_{\theta ,n}}(\lambda )$ $=$ $O_{p}(k_{\theta ,n}/\sqrt{n})$%
. Now use $k_{\theta ,n}=o(\sqrt{n}/\ln (n)^{2})$ to conclude%
\begin{equation*}
\mathcal{E}_{k_{\theta ,n}}=O_{p}\left( \Delta _{k_{\theta
,n}}^{1/3}(\lambda )\max \left\{ 1,\ln \left( k_{\theta ,n}/\Delta
_{k_{\theta ,n}}(\lambda )\right) \right\} ^{2/3}\right) =O_{p}\left( \left( 
\frac{k_{\theta ,n}}{\sqrt{n}}\right) ^{1/3}\ln (n)^{2/3}\right) =o_{p}(1)
\end{equation*}%
\textbf{Step 2.}\qquad Now let $\{\boldsymbol{\tilde{Z}}_{(i)}(\lambda
)\}_{i\in \mathbb{N}}$, $\boldsymbol{\tilde{Z}}_{(i)}(\lambda )$ $\sim $ $%
N(0,\tilde{\sigma}_{(i)}^{2}(\lambda ))$ and $\tilde{\sigma}%
_{(i)}^{2}(\lambda )$ $=$ $\lim_{n\rightarrow \infty }\tilde{\sigma}%
_{n,(i)}^{2}(\lambda )$, be an independent copy of the Theorem \ref%
{th:max_theta_hat}.a null process $\{\boldsymbol{Z}_{(i)}(\lambda )\}_{i\in 
\mathbb{N}}$, independent of the asymptotic draw $\{x_{t},y_{t}\}_{t=1}^{%
\infty }$. It remains to prove $\max_{1\leq i\leq k_{\theta ,n}}|[%
\boldsymbol{0}_{k_{\delta }}^{\prime },1]\mathcal{\tilde{Z}}_{(i)}^{(0)}|$ $%
\Rightarrow ^{p}$ $\max_{i\in \mathbb{N}}|\boldsymbol{\tilde{Z}}_{(i)}([%
\boldsymbol{0}_{k_{\delta }}^{\prime },1])|$.

Define $\mathring{\Delta}_{k_{\theta ,n}}(\lambda )$ $\equiv $ $\max_{1\leq
i,j\leq k_{\theta ,n}}|\tilde{\sigma}_{n,(ij)}^{2}(\lambda )$ $-$ $\tilde{%
\sigma}_{(ij)}^{2}(\lambda )|$ and $\mathcal{\bar{H}}_{(i)}^{(0)}$ $\equiv $ 
$\lim_{n\rightarrow \infty }\mathcal{\bar{H}}_{n,(i)}^{(0)}$, where 
\begin{equation*}
\tilde{\sigma}_{(ij)}^{2}(\lambda )\equiv \lim_{n\rightarrow \infty }\tilde{%
\sigma}_{n,(ij)}^{2}(\lambda )=\lambda ^{\prime }\mathcal{\bar{H}}%
_{(i)}^{(0)-1}\lim_{n\rightarrow \infty }\frac{1}{n}%
\sum_{t=1}^{n}E[G_{i,t}^{(0)}G_{j,t}^{(0)\prime }]\mathcal{\bar{H}}%
_{(j)}^{(0)-1}\lambda .
\end{equation*}%
By the argument above, if we show $\mathring{\Delta}_{k_{\theta
,n}}^{1/3}(\lambda )$ $\times $ $\max \{1,\ln (k_{\theta ,n}/\mathring{\Delta%
}_{k_{\theta ,n}}(\lambda ))\}^{2/3}$ $=$ $o(1)$ then%
\begin{equation}
\sup_{z\geq 0}\left\vert P\left( \max_{1\leq i\leq k_{\theta ,n}}\left\vert 
\boldsymbol{\tilde{Z}}_{n,(i)}(\lambda )\right\vert \leq z\right) -P\left(
\max_{1\leq i\leq k_{\theta ,n}}\left\vert \boldsymbol{\tilde{Z}}%
_{(i)}(\lambda )\right\vert \leq z\right) \right\vert \rightarrow 0.
\label{ZZ}
\end{equation}%
Adding and subtracting terms, and using Assumption \ref{assum:param_boot}.d,
yield:%
\begin{eqnarray*}
\mathring{\Delta}_{k_{\theta ,n}}(\lambda ) &\leq &\max_{1\leq i,j\leq
k_{\theta ,n}}\left\vert \lambda ^{\prime }\mathcal{\bar{H}}%
_{n,(i)}^{(0)-1}\left( \frac{1}{n}\sum_{t=1}^{n}E\left[
G_{i,t}^{(0)}G_{j,t}^{(0)\prime }\right] -\lim_{n\rightarrow \infty }\frac{1%
}{n}\sum_{t=1}^{n}E[G_{i,t}^{(0)}G_{j,t}^{(0)\prime }]\right) \mathcal{\bar{H%
}}_{n,(j)}^{(0)-1}\lambda \right\vert \\
&&+\max_{1\leq i,j\leq k_{\theta ,n}}\left\vert \lambda ^{\prime }\left( 
\mathcal{\bar{H}}_{n,(i)}^{(0)-1}-\mathcal{\bar{H}}_{(i)}^{(0)-1}\right)
\lim_{n\rightarrow \infty }\frac{1}{n}\sum_{t=1}^{n}E\left[
G_{i,t}^{(0)}G_{j,t}^{(0)\prime }\right] \left( \mathcal{\bar{H}}%
_{n,(j)}^{(0)-1}-\mathcal{\bar{H}}_{(j)}^{(0)-1}\right) \lambda \right\vert
\\
&&+\max_{1\leq i,j\leq k_{\theta ,n}}\left\vert \lambda ^{\prime }\left( 
\mathcal{\bar{H}}_{n,(i)}^{(0)-1}-\mathcal{\bar{H}}_{(i)}^{(0)-1}\right)
\lim_{n\rightarrow \infty }\frac{1}{n}\sum_{t=1}^{n}E\left[
G_{i,t}^{(0)}G_{j,t}^{(0)\prime }\right] \mathcal{\bar{H}}%
_{(j)}^{(0)-1}\lambda \right\vert \\
&&+\max_{1\leq i,j\leq k_{\theta ,n}}\left\vert \lambda ^{\prime }\mathcal{%
\bar{H}}_{(i)}^{(0)-1}\lim_{n\rightarrow \infty }\frac{1}{n}\sum_{t=1}^{n}E%
\left[ G_{i,t}^{(0)}G_{j,t}^{(0)\prime }\right] \left( \mathcal{\bar{H}}%
_{n,(j)}^{(0)-1}-\mathcal{\bar{H}}_{(j)}^{(0)-1}\right) \lambda \right\vert
\\
&& \\
&\leq &K\max_{1\leq i\leq k_{\theta ,n}}\left\Vert \mathcal{\bar{H}}%
_{n,(i)}^{(0)-1}\right\Vert ^{2}\max_{1\leq i,j\leq k_{\theta ,n}}\left\Vert 
\frac{1}{n}\sum_{t=1}^{n}E\left[ G_{i,t}^{(0)}G_{j,t}^{(0)\prime }\right]
-\lim_{n\rightarrow \infty }\frac{1}{n}%
\sum_{t=1}^{n}E[G_{i,t}^{(0)}G_{j,t}^{(0)\prime }]\right\Vert \\
&&+K\max_{1\leq i\leq k_{\theta ,n}}\left\Vert \mathcal{\bar{H}}%
_{n,(i)}^{(0)-1}-\mathcal{\bar{H}}_{(i)}^{(0)-1}\right\Vert ^{2}\max_{1\leq
i,j\leq k_{\theta ,n}}\left\Vert \lim_{n\rightarrow \infty }\frac{1}{n}%
\sum_{t=1}^{n}E\left[ G_{i,t}^{(0)}G_{j,t}^{(0)\prime }\right] \right\Vert \\
&&+K\max_{1\leq i\leq k_{\theta ,n}}\left\Vert \mathcal{\bar{H}}%
_{n,(i)}^{(0)-1}-\mathcal{\bar{H}}_{(i)}^{(0)-1}\right\Vert \max_{i\in 
\mathbb{N}}\left\Vert \mathcal{\bar{H}}_{(i)}^{(0)-1}\right\Vert \max_{1\leq
i,j\leq k_{\theta ,n}}\left\Vert \lim_{n\rightarrow \infty }\frac{1}{n}%
\sum_{t=1}^{n}E\left[ G_{i,t}^{(0)}G_{j,t}^{(0)\prime }\right] \right\Vert \\
&& \\
&=&O\left( 1\right) \times O\left( k_{\theta ,n}/\sqrt{n}\right) +O\left(
k_{\theta ,n}/\sqrt{n}\right) \times O(1)=O\left( k_{\theta ,n}/\sqrt{n}%
\right)
\end{eqnarray*}%
The proof of (\ref{ZZ}) now follows from the Step 1 argument.

Together (\ref{PZ|W-PZ}) and (\ref{ZZ}) yield convergence in Kolmogorov
distance: 
\begin{equation*}
\sup_{z\geq 0}\left\vert P\left( \max_{1\leq i\leq k_{\theta ,n}}\left\vert
\lambda ^{\prime }\mathcal{\tilde{Z}}_{(i)}^{(0)}\right\vert \leq z|%
\mathfrak{S}_{n}\right) -P\left( \max_{1\leq i\leq k_{\theta ,n}}\left\vert 
\boldsymbol{\tilde{Z}}_{(i)}(\lambda )\right\vert \leq z\right) \right\vert 
\overset{p}{\rightarrow }0.
\end{equation*}%
Hence $\max_{1\leq i\leq k_{\theta ,n}}|[\boldsymbol{0}_{k_{\delta
}}^{\prime },1]\mathcal{\tilde{Z}}_{(i)}^{(0)}|$ $|\mathfrak{S}_{n}$ $%
\overset{d}{\rightarrow }$ $\max_{i\in \mathbb{N}}|\boldsymbol{\tilde{Z}}%
_{(i)}([\boldsymbol{0}_{k_{\delta }}^{\prime },1])|$ $awp1$ with respect to $%
\{x_{t},y_{t}\}_{t=1}^{\infty }$. Therefore $\max_{1\leq i\leq k_{\theta
,n}}|[\boldsymbol{0}_{k_{\delta }}^{\prime },1]\mathcal{\tilde{Z}}%
_{(i)}^{(0)}|$ $\Rightarrow ^{p}$ $\max_{i\in \mathbb{N}}|\boldsymbol{\tilde{%
Z}}_{(i)}([\boldsymbol{0}_{k_{\delta }}^{\prime },1])|$, cf. 
\citet[Section
3]{GineZinn1990}. $\mathcal{QED}.\bigskip $\newline
\textbf{Proof of Lemma \ref{lm_boot_theta}.}\qquad Lemmas \ref%
{lm:b_boot_expansion} and \ref{lm_ZZ}, coupled with $\max_{1\leq i\leq
k_{\theta ,n}}|\mathcal{W}_{n,i}$ $-$ $\mathcal{W}_{i}|$ $=$ $o_{p}(1/(\ln
(n))^{2})$ for non-stochastic $\mathcal{W}_{i}$ $\in $ $(0,\infty )$, prove
the claim. $\mathcal{QED}$.\bigskip \newline
\textbf{Proof of Theorem \ref{th:p_value_boot}}.\medskip \newline
\textbf{Claim (a).}\qquad\ Assumption \ref{assum:param_boot}.c,d imply
Assumption \ref{assum:max_dist}$^{\ast }$ holds. Therefore $\max_{1\leq
i\leq k_{\theta ,n}}|\sqrt{n}\mathcal{W}_{n,i}\hat{\theta}_{i}|$ $\overset{d}%
{\rightarrow }$ $\max_{i\in \mathbb{N}}|\mathcal{W}_{i}\boldsymbol{Z}_{(i)}([%
\boldsymbol{0}_{k_{\delta }}^{\prime },1])|$ by Lemma \ref%
{lm:suff_Assum4abiic} and Theorem \ref{th:max_theta_hat}.a. Furthermore by
Lemma \ref{lm_boot_theta} $\max_{1\leq i\leq k_{\theta ,n}}\mathcal{\tilde{T}%
}_{n}$ $\Rightarrow ^{p}$ $\max_{i\in \mathbb{N}}|\mathcal{W}_{i}\boldsymbol{%
\tilde{Z}}_{(i)}([\boldsymbol{0}_{k_{\delta }}^{\prime },1])|$, where $\{%
\boldsymbol{\tilde{Z}}_{(i)}(\lambda )\}_{i\in \mathbb{N}}$ is an
independent copy of \linebreak $\{\boldsymbol{Z}_{(i)}(\lambda )\}_{i\in 
\mathbb{N}}$, independent of the asymptotic draw $\{x_{t},y_{t}\}_{t=1}^{%
\infty }$. The claim now follows from arguments in \citet[p. 427]{Hansen1996}%
.\medskip \newline
\textbf{Claim (b).}\qquad Let $H_{1}$ hold, and define $\mathfrak{S}_{n}$ $%
\equiv $ $\{x_{t},y_{t}\}_{t=1}^{n}$. In view of $\max_{1\leq i\leq
k_{\theta ,n}}\mathcal{\tilde{T}}_{n}$ $\Rightarrow ^{p}$ $\max_{i\in 
\mathbb{N}}|\mathcal{W}_{i}\boldsymbol{\tilde{Z}}_{(i)}([\boldsymbol{0}%
_{k_{\delta }}^{\prime },1])|$ we have \citep[eq.
(3.4)]{GineZinn1990} 
\begin{equation*}
\sup_{z>0}\left\vert \mathcal{P}\left( \mathcal{\tilde{T}}_{n}\leq z|%
\mathfrak{S}_{n}\right) -P\left( \max\nolimits_{1\leq i\leq k_{\theta
,n}}\left\vert \boldsymbol{\mathring{Z}}_{(i)}([\boldsymbol{0}_{k_{\delta
}}^{\prime },1])\right\vert \leq z\right) \right\vert \overset{p}{%
\rightarrow }0,
\end{equation*}%
where $\{\boldsymbol{\mathring{Z}}_{(i)}(\cdot )\}_{i\in \mathbb{N}}$ is an
independent copy of $\{\boldsymbol{\tilde{Z}}_{(i)}(\cdot )\}_{i\in \mathbb{N%
}}$. Furthermore, $\tilde{p}_{n,\mathcal{M}_{n}}$ $=$ $\mathcal{P}(\mathcal{%
\tilde{T}}_{n,1}$ $>$ $\mathcal{T}_{n}|\mathfrak{S}_{n})$ $+$ $o_{p}(1)$ by
the Glivenko-Cantelli theorem and independence across bootstrap samples.
Also, $\mathcal{T}_{n}$ $\overset{p}{\rightarrow }$ $\infty $ by Theorem \ref%
{th:max_theta_hat}.b. Therefore $\mathcal{P}_{\mathfrak{S}_{n}}(\mathcal{%
\tilde{T}}_{n,1}$ $>$ $\mathcal{T}_{n})$ $\overset{p}{\rightarrow }$ $0$
hence $\tilde{p}_{n,\mathcal{M}_{n}}$ $=$ $\mathcal{P}_{\mathfrak{S}_{n}}(%
\mathcal{\tilde{T}}_{n,1}$ $\geq $ $\mathcal{T}_{n})$ $+$ $o_{p}(1)$ $%
\overset{p}{\rightarrow }$ $0$, and therefore $P(\tilde{p}_{n,\mathcal{M}%
_{n}}$ $<$ $\alpha )$ $\rightarrow $ $1$. $\mathcal{QED}$.

\setstretch{.5} 
\bibliographystyle{econometrica}
\bibliography{refs_maxtest_many_zeros}

\singlespacing\setstretch{1} \clearpage
\newpage

\begin{sidewaystable}[tbp]
\caption{Rejection Frequencies under $H_{0}: \theta = 0$} \label%
{table:max_wald_rejH0}

\begin{center}
{\small
\resizebox{\textwidth}{!}{
\begin{tabular}{l|lll|lll|lll|lll|lll|lll}
\hline\hline
\multicolumn{10}{c|}{$k_{\delta }=0$} & \multicolumn{9}{|c}{$k_{\delta }=10$}
\\ \hline
\multicolumn{19}{c}{$n=100$} \\ \hline\hline
& \multicolumn{3}{|c|}{$\mathring{k}_{\theta ,n}=10$} & \multicolumn{3}{|l|}{%
$\mathring{k}_{\theta ,n}=35$} & \multicolumn{3}{|l|}{$\mathring{k}_{\theta
,n}=50$} & \multicolumn{3}{|c|}{$\mathring{k}_{\theta ,n}=10$} &
\multicolumn{3}{|l|}{$\mathring{k}_{\theta ,n}=35$} & \multicolumn{3}{|l}{$%
\mathring{k}_{\theta ,n}=50$} \\ \hline
\multicolumn{1}{r|}{Test / Size} & $\mathbf{1\%}$ & $\mathbf{5\%}$ & $%
\mathbf{10\%}$ & $\mathbf{1\%}$ & $\mathbf{5\%}$ & $\mathbf{10\%}$ & $%
\mathbf{1\%}$ & $\mathbf{5\%}$ & $\mathbf{10\%}$ & $\mathbf{1\%}$ & $\mathbf{%
5\%}$ & $\mathbf{10\%}$ & $\mathbf{1\%}$ & $\mathbf{5\%}$ & $\mathbf{10\%}$
& $\mathbf{1\%}$ & $\mathbf{5\%}$ & $\mathbf{10\%}$ \\ \hline
Max-Test & .008 & .051 & .101 & .006 & .051 & .107 & .008 & .041 & .086 &
.013 & .090 & .118 & .006 & .057 & .129 & .007 & .053 & .118 \\
Max-t-Test & .010 & .048 & .110 & .010 & .055 & .114 & .003 & .042 & .089
& .018 & .072 & .142 & .019 & .092 & .180 & .018 & .107 & .183 \\
Wald & .003 & .035 & .076 & .000 & .001 & .016 & .000 & .000 & .005 & .002
& .057 & .122 & .000 & .018 & .081 & .000 & .003 & .040 \\ \hline\hline
\multicolumn{19}{c}{$n=250$} \\ \hline
& \multicolumn{3}{|c|}{$\mathring{k}_{\theta ,n}=10$} & \multicolumn{3}{|l|}{%
$\mathring{k}_{\theta ,n}=35$} & \multicolumn{3}{|l|}{$\mathring{k}_{\theta
,n}=79$} & \multicolumn{3}{|c|}{$\mathring{k}_{\theta ,n}=10$} &
\multicolumn{3}{|l|}{$\mathring{k}_{\theta ,n}=35$} & \multicolumn{3}{|l}{$%
\mathring{k}_{\theta ,n}=79$} \\ \hline
\multicolumn{1}{r|}{Test / Size} & $\mathbf{1\%}$ & $\mathbf{5\%}$ & $%
\mathbf{10\%}$ & $\mathbf{1\%}$ & $\mathbf{5\%}$ & $\mathbf{10\%}$ & $%
\mathbf{1\%}$ & $\mathbf{5\%}$ & $\mathbf{10\%}$ & $\mathbf{1\%}$ & $\mathbf{%
5\%}$ & $\mathbf{10\%}$ & $\mathbf{1\%}$ & $\mathbf{5\%}$ & $\mathbf{10\%}$
& $\mathbf{1\%}$ & $\mathbf{5\%}$ & $\mathbf{10\%}$ \\ \hline
Max-Test & .003 & .046 & .094 & .010 & .046 & .095 & .007 & .052 & .111 &
.010 & .064 & 116 & .007 & .053 & .110 & .006 & .049 & .101 \\
Max-t-Test & .005 & .043 & .090 & .008 & .052 & .094 & .006 & .040 & .082
& .014 & .064 & .136 & .009 & .057 & .101 & .008 & .069 & .120 \\
Wald & .005 & .041 & .091 & .001 & .019 & .059 & .001 & .003 & 022 & .010
& .053 & .113 & .003 & .035 & .083 & .000 & .005 & .041 \\ \hline\hline
\multicolumn{19}{c}{$n=500$} \\ \hline
& \multicolumn{3}{|c|}{$\mathring{k}_{\theta ,n}=10$} & \multicolumn{3}{|l|}{%
$\mathring{k}_{\theta ,n}=35$} & \multicolumn{3}{|l|}{$\mathring{k}_{\theta
,n}=112$} & \multicolumn{3}{|c|}{$\mathring{k}_{\theta ,n}=10$} &
\multicolumn{3}{|l|}{$\mathring{k}_{\theta ,n}=35$} & \multicolumn{3}{|l}{$%
\mathring{k}_{\theta ,n}=112$} \\ \hline
\multicolumn{1}{r|}{Test / Size} & $\mathbf{1\%}$ & $\mathbf{5\%}$ & $%
\mathbf{10\%}$ & $\mathbf{1\%}$ & $\mathbf{5\%}$ & $\mathbf{10\%}$ & $%
\mathbf{1\%}$ & $\mathbf{5\%}$ & $\mathbf{10\%}$ & $\mathbf{1\%}$ & $\mathbf{%
5\%}$ & $\mathbf{10\%}$ & $\mathbf{1\%}$ & $\mathbf{5\%}$ & $\mathbf{10\%}$
& $\mathbf{1\%}$ & $\mathbf{5\%}$ & $\mathbf{10\%}$ \\ \hline
Max-Test & .015 & .049 & .099 & .006 & .054 & .107 & .007 & .045 & .098 &
.009 & .041 & .103 & .009 & .070 & .117 & .004 & .027 & .085 \\
Max-t-Test & .013 & .041 & .093 & .008 & .048 & .104 & .006 & .042 & .105
& .008 & .046 & .103 & .011 & .061 & .103 & .005 & .034 & .094 \\
Wald & .005 & .034 & .082 & .003 & .028 & .070 & .000 & .004 & .023 & .005
& .041 & .099 & .004 & .029 & .089 & .002 & .015 & .063 \\ \hline\hline
\multicolumn{19}{c}{$n=1000$} \\ \hline
& \multicolumn{3}{|c|}{$\mathring{k}_{\theta ,n}=10$} & \multicolumn{3}{|l|}{%
$\mathring{k}_{\theta ,n}=35$} & \multicolumn{3}{|l|}{$\mathring{k}_{\theta
,n}=158$} & \multicolumn{3}{|c|}{$\mathring{k}_{\theta ,n}=10$} &
\multicolumn{3}{|l|}{$\mathring{k}_{\theta ,n}=35$} & \multicolumn{3}{|l}{$%
\mathring{k}_{\theta ,n}=158$} \\ \hline
\multicolumn{1}{r|}{Test / Size} & $\mathbf{1\%}$ & $\mathbf{5\%}$ & $%
\mathbf{10\%}$ & $\mathbf{1\%}$ & $\mathbf{5\%}$ & $\mathbf{10\%}$ & $%
\mathbf{1\%}$ & $\mathbf{5\%}$ & $\mathbf{10\%}$ & $\mathbf{1\%}$ & $\mathbf{%
5\%}$ & $\mathbf{10\%}$ & $\mathbf{1\%}$ & $\mathbf{5\%}$ & $\mathbf{10\%}$
& $\mathbf{1\%}$ & $\mathbf{5\%}$ & $\mathbf{10\%}$ \\ \hline
Max-Test & .006 & .064 & .115 & .005 & .,044 & .094 & .004 & .005 & .101
& .015 & .063 & .108 & .013 & 045 & .096 & .010 & .051 & .089 \\
Max-t-Test & 012 & .056 & .116 & .010 & .045 & .089 & .006 & .050 & .102
& .010 & .061 & .115 & .014 & .055 & .107 & .012 & .042 & .087 \\
Wald & .015 & .045 & .099 & .007 & 034 & .080 & .000 & .017 & .056 & .010
& .052 & .116 & .007 & .034 & .076 & .002 & .028 & .052 \\ \hline\hline
\end{tabular}
}}
\end{center}

{\small All p-vales are bootstrapped, based on 1,000 independently drawn samples.}

\end{sidewaystable}

\clearpage

\begin{sidewaystable}[tbp]
\caption{Linear Regression: Rejection Frequencies under $H_{1}:\theta_{1} =
.001$ and $\theta_{i}= 0$, $i\ge 2$]} \label{table:max_wald_rejH11}

\begin{center}
{\small
\resizebox{\textwidth}{!}{
\begin{tabular}{l|lll|lll|lll|lll|lll|lll}
\hline\hline
\multicolumn{10}{c|}{$k_{\delta }=0$} & \multicolumn{9}{|c}{$k_{\delta }=10$}
\\ \hline
\multicolumn{19}{c}{$n=100$} \\ \hline\hline
& \multicolumn{3}{|c|}{$\mathring{k}_{\theta ,n}=10$} & \multicolumn{3}{|l|}{%
$\mathring{k}_{\theta ,n}=35$} & \multicolumn{3}{|l|}{$\mathring{k}_{\theta
,n}=50$} & \multicolumn{3}{|c|}{$\mathring{k}_{\theta ,n}=10$} &
\multicolumn{3}{|l|}{$\mathring{k}_{\theta ,n}=35$} & \multicolumn{3}{|l}{$%
\mathring{k}_{\theta ,n}=50$} \\ \hline
\multicolumn{1}{r|}{Test / Size} & $\mathbf{1\%}$ & $\mathbf{5\%}$ & $%
\mathbf{10\%}$ & $\mathbf{1\%}$ & $\mathbf{5\%}$ & $\mathbf{10\%}$ & $%
\mathbf{1\%}$ & $\mathbf{5\%}$ & $\mathbf{10\%}$ & $\mathbf{1\%}$ & $\mathbf{%
5\%}$ & $\mathbf{10\%}$ & $\mathbf{1\%}$ & $\mathbf{5\%}$ & $\mathbf{10\%}$
& $\mathbf{1\%}$ & $\mathbf{5\%}$ & $\mathbf{10\%}$ \\ \hline
Max-Test & .101 & .307 & .466 & .451 & .762 & .870 & .605 & .862 & .938 &
.011 & .058 & .125 & .028 & .119 & .218 & .039 & .149 & .266 \\
Max-t-Test & .142 & .388 & .536 & .531 & .795 & .890 & .673 & .876 & .937
& .019 & .083 & .190 & .053 & .164 & .275 & .099 & .267 & .427 \\
Wald & .020 & .133 & .251 & .002 & .033 & .131 & .000 & .001 & .038 & .004
& .060 & .156 & .000 & .031 & .115 & .000 & .009 & .060 \\ \hline\hline
\multicolumn{19}{c}{$n=250$} \\ \hline
& \multicolumn{3}{|c|}{$\mathring{k}_{\theta ,n}=10$} & \multicolumn{3}{|l|}{%
$\mathring{k}_{\theta ,n}=35$} & \multicolumn{3}{|l|}{$\mathring{k}_{\theta
,n}=79$} & \multicolumn{3}{|c|}{$\mathring{k}_{\theta ,n}=10$} &
\multicolumn{3}{|l|}{$\mathring{k}_{\theta ,n}=35$} & \multicolumn{3}{|l}{$%
\mathring{k}_{\theta ,n}=79$} \\ \hline
\multicolumn{1}{r|}{Test / Size} & $\mathbf{1\%}$ & $\mathbf{5\%}$ & $%
\mathbf{10\%}$ & $\mathbf{1\%}$ & $\mathbf{5\%}$ & $\mathbf{10\%}$ & $%
\mathbf{1\%}$ & $\mathbf{5\%}$ & $\mathbf{10\%}$ & $\mathbf{1\%}$ & $\mathbf{%
5\%}$ & $\mathbf{10\%}$ & $\mathbf{1\%}$ & $\mathbf{5\%}$ & $\mathbf{10\%}$
& $\mathbf{1\%}$ & $\mathbf{5\%}$ & $\mathbf{10\%}$ \\ \hline
Max-Test & .248 & .460 & .600 & .924 & .978 & .992 & 1.00 & 1.00 & 1.00 &
.027 & .127 & .217 & .148 & .347 & .485 & .514 & .742 & .838 \\
Max-t-Test & .219 & .462 & .584 & .914 & .975 & .990 & 1.00 & 1.00 & 1.00
& .037 & .149 & .245 & .281 & .500 & .616 & .617 & .799 & .873 \\
Wald & .068 & .224 & .338 & .212 & .522 & .685 & .309 & .850 & .963 & .031
& .127 & .217 & .034 & .190 & .341 & .008 & .096 & .260 \\ \hline\hline
\multicolumn{19}{c}{$n=500$} \\ \hline
& \multicolumn{3}{|c|}{$\mathring{k}_{\theta ,n}=10$} & \multicolumn{3}{|l|}{%
$\mathring{k}_{\theta ,n}=35$} & \multicolumn{3}{|l|}{$\mathring{k}_{\theta
,n}=112$} & \multicolumn{3}{|c|}{$\mathring{k}_{\theta ,n}=10$} &
\multicolumn{3}{|l|}{$\mathring{k}_{\theta ,n}=35$} & \multicolumn{3}{|l}{$%
\mathring{k}_{\theta ,n}=112$} \\ \hline
\multicolumn{1}{r|}{Test / Size} & $\mathbf{1\%}$ & $\mathbf{5\%}$ & $%
\mathbf{10\%}$ & $\mathbf{1\%}$ & $\mathbf{5\%}$ & $\mathbf{10\%}$ & $%
\mathbf{1\%}$ & $\mathbf{5\%}$ & $\mathbf{10\%}$ & $\mathbf{1\%}$ & $\mathbf{%
5\%}$ & $\mathbf{10\%}$ & $\mathbf{1\%}$ & $\mathbf{5\%}$ & $\mathbf{10\%}$
& $\mathbf{1\%}$ & $\mathbf{5\%}$ & $\mathbf{10\%}$ \\ \hline
Max-Test & .833 & .963 & .987 & 1.00 & 1.00 & 1.00 & 1.00 & 1.00 & 1.00 &
.020 & .122 & .235 & .490 & .739 & .836 & 1.00 & 1.00 & 1.00 \\
Max-t-Test & .931 & .985 & .994 & 1.00 & 1.00 & 1.00 & 1.00 & 1.00 & 1.00
& .063 & .199 & .291 & .680 & .829 & .893 & 1.00 & 1.00 & 1.00 \\
Wald & .693 & .876 & .939 & .971 & .994 & .98 & 1.00 & 1.00 & 1.00 & .029
& .142 & .263 & .199 & .483 & .631 & .297 & .680 & .857 \\ \hline\hline
\multicolumn{19}{c}{$n=1000$} \\ \hline
& \multicolumn{3}{|c|}{$\mathring{k}_{\theta ,n}=10$} & \multicolumn{3}{|l|}{%
$\mathring{k}_{\theta ,n}=35$} & \multicolumn{3}{|l|}{$\mathring{k}_{\theta
,n}=158$} & \multicolumn{3}{|c|}{$\mathring{k}_{\theta ,n}=10$} &
\multicolumn{3}{|l|}{$\mathring{k}_{\theta ,n}=35$} & \multicolumn{3}{|l}{$%
\mathring{k}_{\theta ,n}=158$} \\ \hline
\multicolumn{1}{r|}{Test / Size} & $\mathbf{1\%}$ & $\mathbf{5\%}$ & $%
\mathbf{10\%}$ & $\mathbf{1\%}$ & $\mathbf{5\%}$ & $\mathbf{10\%}$ & $%
\mathbf{1\%}$ & $\mathbf{5\%}$ & $\mathbf{10\%}$ & $\mathbf{1\%}$ & $\mathbf{%
5\%}$ & $\mathbf{10\%}$ & $\mathbf{1\%}$ & $\mathbf{5\%}$ & $\mathbf{10\%}$
& $\mathbf{1\%}$ & $\mathbf{5\%}$ & $\mathbf{10\%}$ \\ \hline
Max-Test & 1.00 & 1.00 & 1.00 & 1.00 & 1.00 & 1.00 & 1.00 & 1.00 & 1.00 &
.463 & .682 & .768 & .979 & .997 & .998 & 1.00 & 1.00 & 1.00 \\
Max-t-Test & 1.00 & 1.00 & 1.00 & 1.00 & 1.00 & 1.00 & 1.00 & 1.00 & 1.00
& .363 & .595 & .710 & .997 & 1.00 & 1.00 & 1.00 & 1.00 & 1.00 \\
Wald & .987 & .999 & 1.00 & 1.00 & 1.00 & 1.00 & 1.00 & 1.00 & 1.00 &
.1940 & .395 & .537 & .821 & .945 & .968 & .998 & 1.00 & 1.00 \\ \hline\hline
\end{tabular}
}}
\end{center}

{\small All p-vales are bootstrapped, based on 1,000 independently drawn samples.}

\end{sidewaystable}

\clearpage

\begin{sidewaystable}[tbp]
\caption{Linear Regression: Rejection Frequencies under $H_{1}:\theta_{i} = i/2$ for
$i=1,...,10$ and $\theta_{i} = 0$ $i\ge 11$]} \label{table:max_wald_rejH12}

\begin{center}
{\small
\resizebox{\textwidth}{!}{
\begin{tabular}{l|lll|lll|lll|lll|lll|lll}
\hline\hline
\multicolumn{10}{c|}{$k_{\delta }=0$} & \multicolumn{9}{|c}{$k_{\delta }=10$}
\\ \hline
\multicolumn{19}{c}{$n=100$} \\ \hline\hline
& \multicolumn{3}{|c|}{$\mathring{k}_{\theta ,n}=10$} & \multicolumn{3}{|l|}{%
$\mathring{k}_{\theta ,n}=35$} & \multicolumn{3}{|l|}{$\mathring{k}_{\theta
,n}=50$} & \multicolumn{3}{|c|}{$\mathring{k}_{\theta ,n}=10$} &
\multicolumn{3}{|l|}{$\mathring{k}_{\theta ,n}=35$} & \multicolumn{3}{|l}{$%
\mathring{k}_{\theta ,n}=50$} \\ \hline
\multicolumn{1}{r|}{Test / Size} & $\mathbf{1\%}$ & $\mathbf{5\%}$ & $%
\mathbf{10\%}$ & $\mathbf{1\%}$ & $\mathbf{5\%}$ & $\mathbf{10\%}$ & $%
\mathbf{1\%}$ & $\mathbf{5\%}$ & $\mathbf{10\%}$ & $\mathbf{1\%}$ & $\mathbf{%
5\%}$ & $\mathbf{10\%}$ & $\mathbf{1\%}$ & $\mathbf{5\%}$ & $\mathbf{10\%}$
& $\mathbf{1\%}$ & $\mathbf{5\%}$ & $\mathbf{10\%}$ \\ \hline
Max-Test & 1.00 & 1.00 & 1.00 & 1.00 & 1.00 & 1.00 & 1.00 & 1.00 & 1.00 &
1.00 & 1.00 & 1.00 & 1.00 & 1.00 & 1.00 & 1.00 & 1.00 & 1.00 \\
Max-t-Test & 1.00 & 1.00 & 1.00 & 1.00 & 1.00 & 1.00 & 1.00 & 1.00 & 1.00
& 1.00 & 1.00 & 1.00 & 1.00 & 1.00 & 1.00 & 1.00 & 1.00 & 1.00 \\
Wald & 1.00 & 1.00 & 1.00 & 1.00 & 1.00 & 1.00 & .869 & 1.00 & 1.00 & 1.00
& 1.00 & 1.00 & 1.00 & 1.00 & 1.00 & .868 & 1.00 & 1.00 \\ \hline\hline
\multicolumn{19}{c}{$n=250$} \\ \hline
& \multicolumn{3}{|c|}{$\mathring{k}_{\theta ,n}=10$} & \multicolumn{3}{|l|}{%
$\mathring{k}_{\theta ,n}=35$} & \multicolumn{3}{|l|}{$\mathring{k}_{\theta
,n}=79$} & \multicolumn{3}{|c|}{$\mathring{k}_{\theta ,n}=10$} &
\multicolumn{3}{|l|}{$\mathring{k}_{\theta ,n}=35$} & \multicolumn{3}{|l}{$%
\mathring{k}_{\theta ,n}=79$} \\ \hline
\multicolumn{1}{r|}{Test / Size} & $\mathbf{1\%}$ & $\mathbf{5\%}$ & $%
\mathbf{10\%}$ & $\mathbf{1\%}$ & $\mathbf{5\%}$ & $\mathbf{10\%}$ & $%
\mathbf{1\%}$ & $\mathbf{5\%}$ & $\mathbf{10\%}$ & $\mathbf{1\%}$ & $\mathbf{%
5\%}$ & $\mathbf{10\%}$ & $\mathbf{1\%}$ & $\mathbf{5\%}$ & $\mathbf{10\%}$
& $\mathbf{1\%}$ & $\mathbf{5\%}$ & $\mathbf{10\%}$ \\ \hline
Max-Test & 1.00 & 1.00 & 1.00 & 1.00 & 1.00 & 1.00 & 1.00 & 1.00 & 1.00 &
1.00 & 1.00 & 1.00 & 1.00 & 1.00 & 1.00 & 1.00 & 1.00 & 1.00 \\
Max-t-Test & 1.00 & 1.00 & 1.00 & 1.00 & 1.00 & 1.00 & 1.00 & 1.00 & 1.00
& 1.00 & 1.00 & 1.00 & 1.00 & 1.00 & 1.00 & 1.00 & 1.00 & 1.00 \\
Wald & 1.00 & 1.00 & 1.00 & 1.00 & 1.00 & 1.00 & 1.00 & 1.00 & 1.00 & 1.00
& 1.00 & 1.00 & 1.00 & 1.00 & 1.00 & 1.00 & 1.00 & 1.00 \\ \hline\hline
\multicolumn{19}{c}{$n=500$} \\ \hline
& \multicolumn{3}{|c|}{$\mathring{k}_{\theta ,n}=10$} & \multicolumn{3}{|l|}{%
$\mathring{k}_{\theta ,n}=35$} & \multicolumn{3}{|l|}{$\mathring{k}_{\theta
,n}=112$} & \multicolumn{3}{|c|}{$\mathring{k}_{\theta ,n}=10$} &
\multicolumn{3}{|l|}{$\mathring{k}_{\theta ,n}=35$} & \multicolumn{3}{|l}{$%
\mathring{k}_{\theta ,n}=112$} \\ \hline
\multicolumn{1}{r|}{Test / Size} & $\mathbf{1\%}$ & $\mathbf{5\%}$ & $%
\mathbf{10\%}$ & $\mathbf{1\%}$ & $\mathbf{5\%}$ & $\mathbf{10\%}$ & $%
\mathbf{1\%}$ & $\mathbf{5\%}$ & $\mathbf{10\%}$ & $\mathbf{1\%}$ & $\mathbf{%
5\%}$ & $\mathbf{10\%}$ & $\mathbf{1\%}$ & $\mathbf{5\%}$ & $\mathbf{10\%}$
& $\mathbf{1\%}$ & $\mathbf{5\%}$ & $\mathbf{10\%}$ \\ \hline
Max-Test & 1.00 & 1.00 & 1.00 & 1.00 & 1.00 & 1.00 & 1.00 & 1.00 & 1.00 &
1.00 & 1.00 & 1.00 & 1.00 & 1.00 & 1.00 & 1.00 & 1.00 & 1.00 \\
Max-t-Test & 1.00 & 1.00 & 1.00 & 1.00 & 1.00 & 1.00 & 1.00 & 1.00 & 1.00
& 1.00 & 1.00 & 1.00 & 1.00 & 1.00 & 1.00 & 1.00 & 1.00 & 1.00 \\
Wald & 1.00 & 1.00 & 1.00 & 1.00 & 1.00 & 1.00 & 1.00 & 1.00 & 1.00 & 1.00
& 1.00 & 1.00 & 1.00 & 1.00 & 1.00 & 1.00 & 1.00 & 1.00 \\ \hline\hline
\multicolumn{19}{c}{$n=1000$} \\ \hline
& \multicolumn{3}{|c|}{$\mathring{k}_{\theta ,n}=10$} & \multicolumn{3}{|l|}{%
$\mathring{k}_{\theta ,n}=35$} & \multicolumn{3}{|l|}{$\mathring{k}_{\theta
,n}=158$} & \multicolumn{3}{|c|}{$\mathring{k}_{\theta ,n}=10$} &
\multicolumn{3}{|l|}{$\mathring{k}_{\theta ,n}=35$} & \multicolumn{3}{|l}{$%
\mathring{k}_{\theta ,n}=158$} \\ \hline
\multicolumn{1}{r|}{Test / Size} & $\mathbf{1\%}$ & $\mathbf{5\%}$ & $%
\mathbf{10\%}$ & $\mathbf{1\%}$ & $\mathbf{5\%}$ & $\mathbf{10\%}$ & $%
\mathbf{1\%}$ & $\mathbf{5\%}$ & $\mathbf{10\%}$ & $\mathbf{1\%}$ & $\mathbf{%
5\%}$ & $\mathbf{10\%}$ & $\mathbf{1\%}$ & $\mathbf{5\%}$ & $\mathbf{10\%}$
& $\mathbf{1\%}$ & $\mathbf{5\%}$ & $\mathbf{10\%}$ \\ \hline
Max-Test & 1.00 & 1.00 & 1.00 & 1.00 & 1.00 & 1.00 & 1.00 & 1.00 & 1.00 &
1.00 & 1.00 & 1.00 & 1.00 & 1.00 & 1.00 & 1.00 & 1.00 & 1.00 \\
Max-t-Test & 1.00 & 1.00 & 1.00 & 1.00 & 1.00 & 1.00 & 1.00 & 1.00 & 1.00
& 1.00 & 1.00 & 1.00 & 1.00 & 1.00 & 1.00 & 1.00 & 1.00 & 1.00 \\
Wald & 1.00 & 1.00 & 1.00 & 1.00 & 1.00 & 1.00 & 1.00 & 1.00 & 1.00 & 1.00
& 1.00 & 1.00 & 1.00 & 1.00 & 1.00 & 1.00 & 1.00 & 1.00 \\ \hline\hline
\end{tabular}
}}
\end{center}

{\small All p-vales are bootstrapped, based on 1,000 independently drawn samples.}

\end{sidewaystable}

\clearpage

\begin{sidewaystable}[tbp]
\caption{Linear Regression: Rejection Frequencies under $H_{1}:\theta_{i} =
.001$ for each $i$} \label{table:max_wald_rejH13}

\begin{center}
{\small
\resizebox{\textwidth}{!}{
\begin{tabular}{l|lll|lll|lll|lll|lll|lll}
\hline\hline
\multicolumn{10}{c|}{$k_{\delta }=0$} & \multicolumn{9}{|c}{$k_{\delta }=10$}
\\ \hline
\multicolumn{19}{c}{$n=100$} \\ \hline\hline
& \multicolumn{3}{|c|}{$\mathring{k}_{\theta ,n}=10$} & \multicolumn{3}{|l|}{%
$\mathring{k}_{\theta ,n}=35$} & \multicolumn{3}{|l|}{$\mathring{k}_{\theta
,n}=50$} & \multicolumn{3}{|c|}{$\mathring{k}_{\theta ,n}=10$} &
\multicolumn{3}{|l|}{$\mathring{k}_{\theta ,n}=35$} & \multicolumn{3}{|l}{$%
\mathring{k}_{\theta ,n}=50$} \\ \hline
\multicolumn{1}{r|}{Test / Size} & $\mathbf{1\%}$ & $\mathbf{5\%}$ & $%
\mathbf{10\%}$ & $\mathbf{1\%}$ & $\mathbf{5\%}$ & $\mathbf{10\%}$ & $%
\mathbf{1\%}$ & $\mathbf{5\%}$ & $\mathbf{10\%}$ & $\mathbf{1\%}$ & $\mathbf{%
5\%}$ & $\mathbf{10\%}$ & $\mathbf{1\%}$ & $\mathbf{5\%}$ & $\mathbf{10\%}$
& $\mathbf{1\%}$ & $\mathbf{5\%}$ & $\mathbf{10\%}$ \\ \hline
Max-Test & 1.00 & 1.00 & 1.00 & 1.00 & 1.00 & 1.00 & 1.00 & 1.00 & 1.00 &
.729 & .917 & .965 & .995 & 1.00 & 1.00 & 1.00 & 1.00 & 1.00 \\
Max-t-Test & 1.00 & 1.00 & 1.00 & 1.00 & 1.00 & 1.00 & 1.00 & 1.00 & 1.00
& .706 & .898 & .951 & 1.00 & 1.00 & 1.00 & 1.00 & 1.00 & 1.00 \\
Wald & 1.00 & 1.00 & 1.00 & .996 & 1.00 & 1.00 & .885 & 1.00 & 1.00 & .499
& .812 & .909 & .991 & 1.00 & 1.00 & .747 & 1.00 & 1.00 \\ \hline\hline
\multicolumn{19}{c}{$n=250$} \\ \hline
& \multicolumn{3}{|c|}{$\mathring{k}_{\theta ,n}=10$} & \multicolumn{3}{|l|}{%
$\mathring{k}_{\theta ,n}=35$} & \multicolumn{3}{|l|}{$\mathring{k}_{\theta
,n}=79$} & \multicolumn{3}{|c|}{$\mathring{k}_{\theta ,n}=10$} &
\multicolumn{3}{|l|}{$\mathring{k}_{\theta ,n}=35$} & \multicolumn{3}{|l}{$%
\mathring{k}_{\theta ,n}=79$} \\ \hline
\multicolumn{1}{r|}{Test / Size} & $\mathbf{1\%}$ & $\mathbf{5\%}$ & $%
\mathbf{10\%}$ & $\mathbf{1\%}$ & $\mathbf{5\%}$ & $\mathbf{10\%}$ & $%
\mathbf{1\%}$ & $\mathbf{5\%}$ & $\mathbf{10\%}$ & $\mathbf{1\%}$ & $\mathbf{%
5\%}$ & $\mathbf{10\%}$ & $\mathbf{1\%}$ & $\mathbf{5\%}$ & $\mathbf{10\%}$
& $\mathbf{1\%}$ & $\mathbf{5\%}$ & $\mathbf{10\%}$ \\ \hline
Max-Test & 1.00 & 1.00 & 1.00 & 1.00 & 1.00 & 1.00 & 1.00 & 1.00 & 1.00 &
.583 & .879 & .947 & 1.00 & 1.00 & 1.00 & 1.00 & 1.00 & 1.00 \\
Max-t-Test & 1.00 & 1.00 & 1.00 & 1.00 & 1.00 & 1.00 & 1.00 & 1.00 & 1.00
& .883 & .964 & .991 & 1.00 & 1.00 & 1.00 & 1.00 & 1.00 & 1.00 \\
Wald & 1.00 & 1.00 & 1.00 & 1.00 & 1.00 & 1.00 & 1.00 & 1.00 & 1.00 & .756
& .926 & .965 & 1.00 & 1.00 & 1.00 & 1.00 & 1.00 & 1.00 \\ \hline\hline
\multicolumn{19}{c}{$n=500$} \\ \hline
& \multicolumn{3}{|c|}{$\mathring{k}_{\theta ,n}=10$} & \multicolumn{3}{|l|}{%
$\mathring{k}_{\theta ,n}=35$} & \multicolumn{3}{|l|}{$\mathring{k}_{\theta
,n}=112$} & \multicolumn{3}{|c|}{$\mathring{k}_{\theta ,n}=10$} &
\multicolumn{3}{|l|}{$\mathring{k}_{\theta ,n}=35$} & \multicolumn{3}{|l}{$%
\mathring{k}_{\theta ,n}=112$} \\ \hline
\multicolumn{1}{r|}{Test / Size} & $\mathbf{1\%}$ & $\mathbf{5\%}$ & $%
\mathbf{10\%}$ & $\mathbf{1\%}$ & $\mathbf{5\%}$ & $\mathbf{10\%}$ & $%
\mathbf{1\%}$ & $\mathbf{5\%}$ & $\mathbf{10\%}$ & $\mathbf{1\%}$ & $\mathbf{%
5\%}$ & $\mathbf{10\%}$ & $\mathbf{1\%}$ & $\mathbf{5\%}$ & $\mathbf{10\%}$
& $\mathbf{1\%}$ & $\mathbf{5\%}$ & $\mathbf{10\%}$ \\ \hline
Max-Test & 1.00 & 1.00 & 1.00 & 1.00 & 1.00 & 1.00 & 1.00 & 1.00 & 1.00 &
1.00 & 1.00 & 1.00 & 1.00 & 1.00 & 1.00 & 1.00 & 1.00 & 1.00 \\
Max-t-Test & 1.00 & 1.00 & 1.00 & 1.00 & 1.00 & 1.00 & 1.00 & 1.00 & 1.00
& 1.00 & 1.00 & 1.00 & 1.00 & 1.00 & 1.00 & 1.00 & 1.00 & 1.00 \\
Wald & 1.00 & 1.00 & 1.00 & 1.00 & 1.00 & 1.00 & 1.00 & 1.00 & 1.00 & 1.00
& 1.00 & 1.00 & 1.00 & 1.00 & 1.00 & 1.00 & 1.00 & 1.00 \\ \hline\hline
\multicolumn{19}{c}{$n=1000$} \\ \hline
& \multicolumn{3}{|c|}{$\mathring{k}_{\theta ,n}=10$} & \multicolumn{3}{|l|}{%
$\mathring{k}_{\theta ,n}=35$} & \multicolumn{3}{|l|}{$\mathring{k}_{\theta
,n}=158$} & \multicolumn{3}{|c|}{$\mathring{k}_{\theta ,n}=10$} &
\multicolumn{3}{|l|}{$\mathring{k}_{\theta ,n}=35$} & \multicolumn{3}{|l}{$%
\mathring{k}_{\theta ,n}=158$} \\ \hline
\multicolumn{1}{r|}{Test / Size} & $\mathbf{1\%}$ & $\mathbf{5\%}$ & $%
\mathbf{10\%}$ & $\mathbf{1\%}$ & $\mathbf{5\%}$ & $\mathbf{10\%}$ & $%
\mathbf{1\%}$ & $\mathbf{5\%}$ & $\mathbf{10\%}$ & $\mathbf{1\%}$ & $\mathbf{%
5\%}$ & $\mathbf{10\%}$ & $\mathbf{1\%}$ & $\mathbf{5\%}$ & $\mathbf{10\%}$
& $\mathbf{1\%}$ & $\mathbf{5\%}$ & $\mathbf{10\%}$ \\ \hline
Max-Test & 1.00 & 1.00 & 1.00 & 1.00 & 1.00 & 1.00 & 1.00 & 1.00 & 1.00 &
1.00 & 1.00 & 1.00 & 1.00 & 1.00 & 1.00 & 1.00 & 1.00 & 1.00 \\
Max-t-Test & 1.00 & 1.00 & 1.00 & 1.00 & 1.00 & 1.00 & 1.00 & 1.00 & 1.00
& 1.00 & 1.00 & 1.00 & 1.00 & 1.00 & 1.00 & 1.00 & 1.00 & 1.00 \\
Wald & 1.00 & 1.00 & 1.00 & 1.00 & 1.00 & 1.00 & 1.00 & 1.00 & 1.00 & 1.00
& 1.00 & 1.00 & 1.00 & 1.00 & 1.00 & 1.00 & 1.00 & 1.00 \\ \hline\hline
\end{tabular}%
}}
\end{center}

{\small All p-vales are bootstrapped, based on 1,000 independently drawn samples.}

\end{sidewaystable}

\clearpage

\end{document}